\documentclass[12pt]{article}

\usepackage{epic}
\usepackage{eepic}
\usepackage{epsf} 
\usepackage{ amssymb, amscd}
\usepackage{amsmath}
\usepackage{color}
\numberwithin{equation}{section}

\newcommand{\lori}{\longrightarrow}

\newcommand\as[1]  {[\alpha_{{#1}}]}
\newcommand\asx[2] {[\alpha_{{#1}}^{({#2})}]}

\def\cA            {{\mathcal{A}}}

\def\cE            {{\mathcal{E}}}

\def\cH            {{\mathcal{H}}}

\def\cK            {{\mathcal{K}}}

\def\cO            {{\mathcal{O}}}

\def\cS            {{\mathcal{S}}}

\def\bbA           {\mathbb{A}}
\def\bbC           {\mathbb{C}}
\def\bbD           {\mathbb{D}}
\def\bbE           {\mathbb{E}}  

\def\bbN           {\mathbb{N}}

\def\bbP           {\mathbb{P}}
\def\bbQ           {\mathbb{Q}}
\def\bbR           {\mathbb{R}}
\def\bbT           {\mathbb{T}}
\def\bbZ           {\mathbb{Z}}
\def\SUz           {{\mathit{SU}}(2)}
\def\SUn           {{\mathit{SU}}(n)}

\def\MXN           {{}_M {\mathcal{X}}_N}
\def\MXM           {{}_M {\mathcal{X}}_M}

\def\MXMo          {{}_M^{} {\mathcal{X}}_M^0}

\def\MXMp          {{}_M^{} {\mathcal{X}}_M^+}
\def\MXMm          {{}_M^{} {\mathcal{X}}_M^-}
\def\MXMpm         {{}_M^{} {\mathcal{X}}_M^\pm}

\def\NXN           {{}_N {\mathcal{X}}_N}

\def\NXM           {{}_N {\mathcal{X}}_M}

\def\a             {\alpha}
\def\la            {\lambda}

\def\rmc           {{\mathrm{c}}}

\def\SUd           {{\mathit{SU}}(3)}

\def\SUn           {{\mathit{SU}}(n)}

\def\SUz           {{\mathit{SU}}(2)}

\def\SUf           {{\mathit{SU}}(4)}

\def\tr            {{\mathrm{tr}}}

\def\Mat           {{\mathrm{Mat}}}
\def\i{{\rm i}}
\def\sdprod{{\times\!\vrule height5pt depth0pt width0.4pt\,}}

\newcommand\asp[1] {[\alpha^+_{{#1}}]}
\newcommand\asm[1] {[\alpha^-_{{#1}}]}
\newcommand\aspm[1]{[\alpha^\pm_{{#1}}]}
\newcommand\asprod[2]{[\alpha^+_{{#1}}\circ\alpha^-_{{#2}}]}

\def\LSUn          {{\mathit{LSU}}(n)}
\def\LSUz          {{\mathit{LSU}}(2)}
\def\LISUn         {{\mathit{L}}_I{\mathit{SU}}(n)}

\def\LIcSUn        {{\mathit{L}}_{I^\rmc}{\mathit{SU}}(n)}
\def\id            {\mathrm{id}}
\def\NXMt          {{}_{N} {\mathcal{X}}_{\tilde{M}}}

\title{Modular Invariants and Twisted Equivariant K-theory}
\author{
{\sc David E.\ Evans}\\
 {\footnotesize School of Mathematics, Cardiff University,}\\
 {\footnotesize Senghennydd Road, Cardiff CF24 4AG, Wales, U.K.}\\
 {\footnotesize e-mail: {\tt EvansDE@cf.ac.uk}}\\ \\
 {\sc Terry  Gannon }\\
 {\footnotesize Department of Mathematics, University of Alberta,}\\
{\footnotesize Edmonton, Alberta, Canada T6G 2G1}\\
{\footnotesize e-mail: {\tt tgannon@math.ualberta.ca}} }

\begin{document}
\maketitle

\begin{abstract} Freed-Hopkins-Teleman expressed the Verlinde algebra as
twisted equivariant $K$-theory. We study how to recover the full system
(fusion algebra of defect lines), nimrep (cylindrical partition function),
etc of modular invariant partition functions of conformal field theories 
associated to loop
groups. We work out several examples corresponding to conformal embeddings 
and orbifolds. We identify a new aspect of the A-D-E pattern of $SU(2)$ 
modular invariants.
\end{abstract}

{\footnotesize
\tableofcontents
}

\section{Introduction}

Let $G$ be  a compact connected simply connected Lie group. The equivalence 
classes of its finite-dimensional
representations under direct sum and tensor product form the {\it
representation ring} $R_G$. This ring can be realised as the equivariant
(topological) $K$-group $K_G^{0}(pt)$ of $G$ acting trivially on a point $pt$. 
 From Bott periodicity, $K^{2+*}_G(X)=
K^{*}_G(X)$; the other $K$-group is $K_G^{1}(pt)=0$.
Note that these $K$-groups depend on the Lie group $G$ and not on the
Lie algebra $\mathfrak{g}(G)$: replacing $G$ with $G/Z$ for a central
subgroup $Z$ changes the representation ring but not the Lie algebra.
This sensitivity to the Lie group will persist throughout this
paper, and is fundamental to what follows.

Let $LG$ be the loop group $\{f:S^1\rightarrow G\}$ associated to $G$.
The most interesting representations of $LG$ are projective; the
corresponding central extensions $LG_k$ of $LG$ by $S^1$ are parametrised by
the {\it level} $k\in\bbZ$. The
 loop group analogue of the ring $R_G$ is the {\it Verlinde algebra} $Ver_k(G)$,
spanned by the (equivalence classes of) positive energy representations of 
$LG_k$, with operations direct sum and fusion product.  Physically, $Ver_k(G)$ is the
fusion algebra of the Wess-Zumino-Witten conformal field theory corresponding
to $G$, with a central charge determined by $k$.
Freed-Hopkins-Teleman \cite{FHT,FHTi,FHTii,FHTlg} identify $Ver_k(G)$
with the twisted equivariant $K$-group
${}^\tau K^{{\rm dim}(G)}_G(G)$ for some twist $\tau\in 
H^3_G(G;\bbZ)$ ($\simeq\bbZ$ for $G$ simple) depending on $k$, where $G$ acts on itself by
conjugation. The other $K$-group, namely
${}^\tau K^{1+{\rm dim}(G)}_G(G)$, again is 0. The dimension shift here, by
dim$(G)$, is due to an implicit application of Poincar\'e duality, and is
a hint that things here are more naturally expressed in (the equivalent 
language of) $K$-homology. The ring structure of $Ver_k(G)$ is recovered
from the push-forward of group multiplication $m:G\times G\rightarrow G$, whereas its
$R_G$-module structure comes from the push-forward of
the inclusion $1\hookrightarrow G$ of the identity. Strictly speaking, $\tau\in
H^3_G(G;\bbZ)$ only determines the $K$-theory $^\tau K_G^*(G)$ as an additive
group; the ring structure comes from a choice of lift (if it exists) of
the 3-cocycle $\tau$ to $H^4(BG;\bbZ)$. But for $G$ compact connected simply
connected, transgression identifies $H^4(BG;\bbZ)$ with $H^3_G(G;\bbZ)$ so we 
can (and will) let $\tau$ parametrise the full ring structure.
Note that, as in the preceding paragraph, the $K$-groups depend on $G$ and
not its Lie algebra -- eg. the Verlinde algebra $^\tau K_{SO3}^1(SO(3))$
involves only nonspinors and a fixed point resolution arises, as one would
expect with Wess-Zumino-Witten on $SO(3)$.

Those authors were helped to their loop group theory, through considering a finite
group toy model. But in \cite{Ev}, the finite group story is developed much
more completely, using the braided subfactor approach. Let a finite group
$G$ act on itself by conjugation; then the transgression
 $H^4(BG;\bbZ)\rightarrow H^3_G(G;\bbZ)$, $\sigma\mapsto\tau(\sigma)$,
will in general be neither surjective nor injective. For twist $\sigma\in H^4(BG;\bbZ)$,
${}^{\tau(\sigma)} K^0_G(G)$ is as a ring
the Verlinde algebra of the $\sigma$-twisted quantum double
${\cal D}^\sigma(G)$ of $G$, and ${}^{\tau(\sigma)} K^1_G(G)$ is again 0.
For example, for $G=\bbZ_n$, $H^3_G(G;\bbZ)= 0$ while $H^4(BG;\bbZ)\simeq
\bbZ_n$; as an additive group, $^{\tau(\sigma)}K_G^0(G)\simeq \bbZ^{n^2}$ is
independent of $\sigma$, but as a ring, $^{\tau(\sigma)}K^0_G(G)$ is the
group ring $\bbZ[\bbZ_d\times\bbZ_{n^2/d}]$ where $d={\rm gcd}(2n,\sigma)$.

However \cite{Ev} argues that much more is possible. The viable modular invariants
for ${\cal D}^\tau(G)$ are parametrised by pairs $(H,\psi)$ for a subgroup $H$ 
of $G\times G$ and $\psi\in H^2(BH;\bbT)$ \cite{EP1,Ost}. Let $H\times H$ act 
on $G\times G$ on the left and right: $(h_L,h_R).(g,g')=h_L(g,g')h_R$. Then 
$^\tau K^0_{H\times H}(G\times G)$ can be identified with the {\it full system}
(the fusion algebra of defect lines -- in section 1.3 we define this using 
subfactors and explain its physical meaning),
and again $^\tau K^1_{H\times H}(G\times G)=0$. Choosing $H$ to be the diagonal
subgroup isomorphic to $G$ recovers the Verlinde algebra. We review this in
subsection 1.4.

Finite groups are much simpler than loop groups -- for example there is no
direct analogue of the $(H,\psi)$ parametrisation of viable modular
invariants --  but similar extensions of Freed-Hopkins-Teleman can be expected.
This paper explores these extensions. 
Atiyah \cite{Ati} writes: ``The K-theory approach [to the Verlinde algebra]
is totally new and much more direct than most other ways. It remains to be
thoroughly explored.'' In \cite{FHTlg} v.1, second paragraph, we read that
the relation of twisted K-theory to 
``Chern-Simons (3D TFT) structure ... is, at present, not understood.'' 
This provides the context for our work. 

For example, given a {\it conformal
embedding} (see section 1.3) of $LH$ level $k$ in $LG$ level $\ell$ and
appropriate choice of twist $\tau$, we would expect 
$^\tau K^{{\rm dim}({H})}_{{H}}(G)$ to be related to
the full system for the modular invariant of $H$ level $k$ coming from the
diagonal modular invariant of $G$ level $\ell$. The group ${H}$
 acts adjointly on $G$ as usual. 
For the special case of $H=G$ and $k=\ell$, this
construction recovers that of \cite{FHT}.
We find however that for $H\ne G$, this relation is not as direct.
For example, even in
some of the simplest examples, the full system can appear here with a
multiplicity, and
$^\tau K^{1+{\rm dim}({H})}_{{H}}(G)$ may not vanish.
%$^\tau K^{{\rm dim}(\tilde{H})}_{\widetilde{H}}(G)$ to be related to
%the full system for the modular invariant of $H$ level $k$ coming from the
%diagonal modular invariant of $G$ level $\ell$.
%We write  $\widetilde{H}$ here for
%the universal covering  group of the closed subgroup $H$ of $G$; $\widetilde{H}$
% acts adjointly on $G$
%by first projecting onto $H$. 
%For the special case of $H=G$ and $k=\ell$, this
%construction recovers that of \cite{FHT}.
%We find however that for $H\ne G$, this relation is not as direct.
%For example, even in
%some of the simplest examples, the full system can appear here with a
%%multiplicity, and
%$^\tau K^{1+{\rm dim}(\tilde{H})}_{\tilde{H}}(G)$ may not vanish.

We give several examples, most importantly $D_4$ and $E_6$ on 
Cappelli-Itzykson-Zuber's A-D-E list of $G=SU(2)$ 
modular invariants \cite{CIZ}. It is intriguing
that the quaternionic and tetrahedral groups $Q_4$ and $BA_4$  play
fundamental roles in this $K$-homological analysis for $D_4$ and $E_6$, 
respectively, as the largest finite stabilisers ($Q_4$ and $BA_4$ correspond
to $\bbD_4$ and $\bbE_6$ in McKay's A-D-E list \cite{McK} of finite 
subgroups of $SU(2)$). We show in section 1.4 that the analogue of conformal 
embeddings  for finite groups works perfectly.

The orbifold construction also seems tractable from this point-of-view. 
In particular, in section 4 we compare
the  Verlinde algebra of the $\pi$-permutation orbifold of $G$ level $k$
(for any permutation group $\pi<S_n$) to
the twisted equivariant $K$-homology of $G^n\sdprod \pi$ acting on $G^n$,
where $G^n$ acts adjointly on itself and $\pi$ acts by permuting the
factors. Again, the analogue for finite groups seems to work perfectly
(see the beginning of section 4).

The most important source of modular invariants for loop groups
are the simple current invariants. As these correspond to strings living on
the nonsimply connected groups $G/Z$ (for $Z$ a subgroup of the centre
of $G$), we would expect the full system to be given by the $K$-homology of
$G\times G$ acting diagonally on $(G/Z)\times (G/Z)$.
We shall investigate this in the sequel to this paper.
By contrast, $^\tau K_G^{{\rm dim}(G)}(G/Z)$ for $G$ acting adjointly on $G/Z$ 
should be
the associated nimrep (and again vanish for dim$(G)+1$). The example of
$^\tau K^*_{SU2}(SO(3))$ was worked out in \cite{BS}. It would be very interesting
to understand the $K$-homology capturing the Verlinde algebra of the
 Goddard-Kent-Olive coset construction; in \cite{SN2} the `chiral algebra'
 (a subalgebra of the Verlinde algebra) of some $N=2$ coset models was
identified with a $K$-group.

These considerations suffice to handle every $SU(2)$ modular invariant,
except for the one called $E_7$. This can be obtained from a $\bbZ_2$-orbifold
of the $D_{10}$ modular invariant; the associated $K$-homology will be 
worked out in the sequel.

We explore these
natural constructions and extensions with a detailed study of several simple
but representative examples.
We construct the relevant (twisted graded) bundles --
Meinrenken \cite{M} recently found an independent construction, elegant but
less general, for some of these bundles, and we compare his to ours at
the end of subsection 2.2. The bulk of the paper consists of the detailed
calculations; in the concluding section we interpret these in the 
context of conformal field theory. To
keep this paper relatively self-contained, we begin
with some background material from $K$-theory/-homology and conformal field theory.

\subsection{K-theoretic preliminaries}

The standard references for $K$-theory, $K$-homology, and their twisted
versions, are \cite{AS,BB,CW,DK,Kar1,Ros1,Ros2}.

$K$-theory or $K$-cohomology on a  compact Hausdorff space $X$ looks at the 
vector bundles over $X$. In the operator algebraic formulation, this can be 
equivalently pictured via finitely generated projective modules over the 
$C^{*}$-algebra $C(X)$, the space of complex valued continuous functions
on $X$. This gives the abelian group $K^{0}(X)$, as the Grothendieck group
or completion of  the semigroup of vector bundles or modules. 
For locally compact spaces, we need to be a bit careful,
with inserting and removing one-point compactifications 
or $K$-theory with compact support.
  More precisely, if $X^{+}$ is the one point compactification of locally
 compact space $X$, then
$K^{0}(X)$ is identified with the kernel of the natural map $K^{0}(X^{+})
\rightarrow K^{0}({\rm point})$. The group $K^{1}(X)$
can be defined via
suspensions as $K^{0}(\bbR \times X)$, or through unitaries modulo the
connected component of the identity
in matrices over $C(X^{+})$.  We  then write
$K^{*}(X) =  K_{*}(C_{0}(X))$, for the $C^{*}$-algebra $C_{0}(X)$ of complex valued
continuous functions on $X$. We  can  identify this with 
$K_{*}(C_{0}(X ; \cK))$, the $K$-cohomology of the $C^{*}$-algebra of the space of 
$\cK$-valued functions on $X$, vanishing at infinity, where $\cK$ is the compact operators
on a separable infinite dimensional Hilbert space $\cH$.

The $K$-homology of a compact Hausdorff space $X$ can be understood as classifying extensions of the form
\begin{equation}\label{extE}
 0 \rightarrow \cK \rightarrow \cE \rightarrow C(X) \rightarrow 0\,.
 \end{equation}
More precisely, the degree-one $K$-homology group $K_{1}(X)$ classifies these extensions, and  again
we can define $K_{1}(X) = K_{1}(X^{+})$  and $K_{0}(X) = K_{1}(\bbR \times X)$, for
a locally compact space $X$, using
one-point compactifications and suspensions. Again, we  then write 
$K_{*}(X) =  K^{*}(C_{0}(X))$, and   identify this with 
$K^{*}(C_{0}(X ; \cK))\,.$

The $C^{*}$-algebra $C_{0}(X ; \cK)$ can be twisted, in the sense of twisting 
this space of sections of the trivial bundle, with fibres the compacts $\cK$, 
over $X$, by taking a non trivial bundle $\cK_{\tau}$ and the corresponding 
space of sections $C_{0}(X ; \cK_{\tau})$ \cite{RW}. These algebras are 
locally Morita equivalent to the the trivial algebras $C_{0}(U ; \cK)$, for 
small open subsets $U \subset X$, but the gluing together of these trivial 
algebras is classified by a $\rm{\check Cech}$ cohomology class of $X$, the 
Dixmier-Douady invariant $H^{3}(X; \bbZ)$. We take a cover $\{ U_{i} : i \}$ 
of $X$ by open sets, with a gluing described by a matching on intersections 
$U_{ij} = U_{i} \cap U_{j}$ given by automorphisms $\mu_{ij}$  of $\cK$, 
where $\mu_{ij}\mu_{jk} = \mu_{ik}$ on triple intersections $U_{ijk} = 
U_{i} \cap U_{j} \cap U_{k}$. This gives an element of $H^{1}(X; \mathrm{Aut}(\cK)) 
\simeq H^{1}(X; PU(\cH))$, where $PU(\cH)$ denotes the projective unitary group
$U(\cH)/\bbT$. The latter group $H^{1}(X; PU(\cH))$ is identified with
$H^{2}(X; \bbT) \simeq H^{3}(X; \bbZ)$, the Dixmier-Douady invariant, by taking
the cocycle $\mu_{ij} = Ad(g_{ij})$ to the $\bbT$-valued cocycle 
$g_{ij}g_{jk}g_{ik}^{*}$.

If we wish to include a grading on space of  sections  then there is an 
additional degree of freedom given by  $H^{1}(X; \bbZ_{2})$. First, we decompose
$\cH \simeq \cH\oplus \cH$, with grading self adjoint unitary $\sigma$ which interchanges these components,
and corresponding grading automorphism $Ad(\sigma)$ on the compacts $\cK$.
The graded automorphisms
Aut$^{gr}(\cK) $ are those that commute with $Ad(\sigma)$, and are
implemented by the graded unitaries $U_{gr}(\cH)$.
The transition functions are now given by
graded automorphims $g_{ij}$  on intersections, yielding an element of 
$$H^{1}(X; \mathrm{Aut}^{gr}(\cK)) \simeq
H^{1}(X; PU_{gr}(\cH)) \simeq H^{1}(X; \bbZ_{2}) \oplus H^{3}(X; \bbZ)\,.$$
The graded automorphisms of $\cK$ are identified with the projective
graded unitaries $PU_{gr}(\cH) = U_{gr}(\cH)/\bbT$.
Then  $H^{1}(X; PU_{gr}(\cH)) \rightarrow
H^{3}(X; \bbZ)$ is obtained by ignoring the grading, while the projection
$H^{1}(X; PU_{gr}(\cH)) \rightarrow  H^{1}(X; \bbZ_{2})$ is
through the degree map $U_{gr} \rightarrow \bbZ_{2}$.

The graded Morita equivalence classes of separable $\bbZ_{2}$-graded
continuous trace algebras, with spectrum $X$ are classified by the graded 
Brauer group of $X$, namely: \cite{EMP1,EMP2}
\begin{equation}\label{h1plush3}
 H^{1}(X; \bbZ_{2}) \oplus H^{3}(X; \bbZ)\,.
\end{equation}
Such an element then defines a graded bundle $\cK_{\tau}$ on $X$,  so that we can form the
graded $C^{*}$-algebra of sections $C_{0}(X ; \cK_{\tau})$, and take the corresponding
(graded Kasparov \cite{GK1,GK2,GK3,BB}) $K$-theory:
\begin{eqnarray}
^{\tau}K^{*}(X)& =& K_{*}(C_{0}(X ; \cK_{\tau}))\,,\nonumber\\
^{\tau}K_{*}(X)& =& K^{*}(C_{0}(X ; \cK_{\tau}))\,.\nonumber
\end{eqnarray}

The graded $K$-theory is understood as follows \cite{Tro, HG}. Let $\cS$ denote $C_{0}(\bbR)$ with
the grading induced by the flip $x \mapsto -x$ on $\bbR$. If $A$ is a unital graded $C^{*}$-algebra, then
the graded $K$-theory is defined by $K_{0}(A) = [[\cS , A \hat\otimes \cK]]$, the space of graded homotopy classes
of graded $*$-homomorphisms from $\cS$ into the graded tensor product 
$A  \hat\otimes \cK$. Non-unital algebras, or locally 
compact spaces are handled by unitalisation or one-point compactification as before and $K_{1}$ is handled by suspension.
This suspension can be realised via Clifford algebras: $K_{1}(A) \simeq K_0(A \hat\otimes {\cal C}_{0,1})\,.$
Here $\hat\otimes$ is the graded tensor product and ${\cal C}_{p,q}$ is the graded Clifford algebra for the quadratic form $(-1_p,1_q)$ on $\bbR ^{p+q}$, so that ${\cal C}_{0,1} \simeq \bbC\oplus\bbC$ with the nontrivial grading.

We need equivariant versions of this twisted $K$-theory, using equivariant $\rm{\check Cech}$ cohomology
to describe these twistings. If a group $G$ acts on our space $X$, we define 
equivariant cohomology by the Borel construction 
$H^{*}_{G}(X) = H^{*}((X \times E{G})/G)$, where $E{G}$ is a contractible
space on which $G$ acts
freely, and the quotient is taken for the diagonal action \cite{Bot}. In particular,
$H^{*}_{G}({\rm point}) = H^{*}(B{G})$, where $B{G}$ is the classifying space $E{G}/G$.

The equivariant graded Morita equivalence classes of separable 
$\bbZ_{2}$-graded $G$-equivariant continuous trace algebras, with spectrum $X$ 
are classified by the equivariant graded Brauer group of $X$, namely:
\begin{equation}\label{Gh1plush3}
 H^{1}_{G}(X; \bbZ_{2}) \oplus H^{3}_{G}(X; \bbZ)
 \,.
\end{equation}

We can form a product on these structures. If $\xi$ and $\eta$ are $G$-equivariant graded bundles of
compact operators on $X$, we can form the product $\xi \hat\otimes_{X} \eta$ as $G$-equivariant graded bundle of
compacts on $X$, so that 
$$C_{0}(X; \cK_{\xi \hat\otimes_{X} \eta}) \simeq  C_{0}(X; \cK_{\xi}) \hat\otimes_{C_{0}(X)}
 C_{0}(X; \cK_{\eta})\,,$$ 
where $\hat\otimes_{C_{0}(X)}$ denotes the graded tensor product of ${C_{0}(X)}$-module graded $C^{*}$-algebras.
In terms of the decomposition (\ref{Gh1plush3}), then $\xi \hat\otimes_{X} \eta$ is identified with
\begin{equation}
\xi \hat\otimes_{X} \eta=
(\xi_{1} + \eta_{1},\xi_{3} + \eta_{3} +\beta(\xi_{1}\eta_{1}))\ ,\label{bock} 
\end{equation}
where  $\xi_{1}\eta_{1} \in H^{2}_{G}(X; \bbZ/2)$
is the cup product and $\beta : H^{2}_{G}(X; \bbZ/2) \rightarrow  H^{3}_{G}(X; \bbZ)$ is the Bockstein homomorphism.
For simplification, we prefer to  write this as $\xi + \eta$. In particular we can identify $\xi$ with
$(\xi_{1}, 0) + (0, \xi_3)$ so that in some sense we can treat the $H^1$ and $H^3$ twists independently.

Twisted equivariant $K$-theory is then defined as:
\begin{eqnarray}
^{\tau}K^{*}_{G}(X)& =& K_{*}(C_{0}(X ; \cK_{\tau})\rtimes G)\,,\\
^{\tau}K_{*}^{G}(X)& =& K^{*}(C_{0}(X ; \cK_{\tau}) \rtimes G)\,,
\label{defkhom}\end{eqnarray}
where `$\rtimes$' here denotes the crossed-product construction.

If $\alpha \in H^1_G(X;  \bbZ/2)$ is a $G$-equivariant real line bundle $L = L_\alpha$ on $X$, with
projection $\pi : L \rightarrow X$, and if $\xi$ is $G$-equivariant graded bundle of compacts on
$X$, with $\xi_1 = \alpha$, then $\pi^*(\xi_3)$ is a $G$-equivariant (ungraded) bundle of compacts on $L$.
Moreover 
$$ C_{0}(X; \cK_{\xi}) \simeq C_{0}(X; \cK_{\xi_1}) \hat\otimes_{C_{0}(X)}
 C_{0}(X; \cK_{\xi_3}) \simeq C_{0}(L; \cK_{\pi^*(\xi_3)})\hat\otimes\, {\cal C}_{0,1}\, . $$ 
This expresses the graded $K$-theory in terms of ungraded $K$-theory:
\begin{equation} \label{lineb}
^{\xi}K^i_G(X) \simeq  {}^{\pi^*(\xi_3)}K^{i+1}_G(L_{\xi_{1}})\,.\end{equation}
See \cite{AH,EMP1}, and Theorem 4.4 and section 5.7 of \cite{Kar1} for details.
Compare Definition 4.12 of \cite{FHT}.
%(Terry - I am not sure I have got this interpretation of Karoubi correct - but it is pretty close)

These $K$-groups $^\tau K^G_*(X)$, $^\tau K_G^*(X)$ possess a natural
$R_G$-module structure, coming from the map $X\rightarrow pt$ \cite{Se} --
as mentioned next section, $K^G_0(pt)\simeq H_G^0(pt)\simeq R_G$.
When the twist $\tau\in H^3_G(X;\bbZ)$ is transgressed from $H^4(X;\bbZ)$,
the $K$-groups
$^\tau K^G_*(X)$ and $^\tau K_G^*(X)$ carry graded ring structures \cite{TX}, 
coming from the external Kasparov product in equivariant $KK$-theory.

\subsection{Assorted practicalities}

By $R_G$ we mean the representation
ring of the group $G$. Some of these we'll need are
\begin{eqnarray}
R_{SU2}&=&\bbZ[\sigma]=\mathrm{Span}\{\sigma_1,\sigma_2,\ldots\}\, ,\nonumber\\
R_{O2}&=&\bbZ[\delta,\kappa]/(\delta\kappa=\kappa,\delta^2=1)=\mathrm{Span}\{1,\delta,
\kappa_1,\kappa_2,\ldots\}\, ,\nonumber\\
R_{\bbT}&=&\bbZ[a^{\pm 1}]=\mathrm{Span}\{1,a,a^{-1},a^2,a^{-2},\ldots\}\, ,
\nonumber\end{eqnarray}
where $\sigma_i$ is the $i$-dimensional $SU(2)$-representation (so $\sigma=
\sigma_2$ is the defining representation), $\delta=det$, $\kappa_i$ is the
two-dimensional $O(2)$-representation with winding number $i$ (so $\kappa
=\kappa_1$ is the defining representation), and $a^i$  is the one-dimensional
representation for the circle $SO(2)=U(1)=S^1=\mathbb{T}$ with winding number $i$.

Restriction makes both $R_{O2}$ and $R_{\bbT}$ into $R_{SU2}$-modules;
the generator $\sigma\in R_{SU2}$ restricts to $\kappa_1\in R_{O2}$ and
to $a+a^{-1}\in R_{\bbT}$. Induction from $R_{\bbT}$ to $R_{O2}$ takes 1 to
$1+\delta$ and $a^i$ to $\kappa_{|i|}$; {\it Dirac} induction from $R_{\bbT}$
(resp. $R_{O2}$) to $R_{SU2}$ is discussed later in this subsection (resp.
in section 2.1).

To keep our calculations under some control, we will usually act with
$G=SU(2)$. Hence its finite subgroups will often arise as stabilisers.
As is well-known (see eg.\ \cite{McK,IN}), these fall into an A-D-E pattern:
they are the cyclic groups $\mathbb{A}_n=\bbZ_{n+1}=\langle 1,1,n\rangle$, 
double-covers $\mathbb{D}_n
=BD_{n-2}=\langle 2,2,n\rangle$ of the
dihedral groups $D_{n-2}$, as well as the binary tetrahedral $\mathbb{E}_6
=BA_4=\langle 2,3,3\rangle$, binary
octahedral $\mathbb{E}_7=BS_4=\langle 2,3,4\rangle$ and
binary icosahedral $\mathbb{E}_8=BA_5=\langle 2,3,5\rangle$
groups, where the {\it binary polyhedral group} $\langle \ell,m,n\rangle$
is defined by \cite{CM}
\begin{equation}\label{coxmos}
\langle \ell,m,n\rangle=\langle a,b,c\,|a^\ell=b^m=c^n=abc\rangle\,.
\end{equation}
 The vertices
of the corresponding extended ${A}$, ${D}$, or ${E}$ diagram are 
labelled with the irreducible representations of the finite subgroup $H$; the embedding
$H\hookrightarrow G$ defines a two-dimensional representation $\rho=\mathrm{Res}_H^G\sigma$
of $H$, and decomposing into irreducibles the tensor of
$\rho$ with the irreducible ones recovers the
edges of that diagram. In this way the Dynkin diagram encodes the
$R_{SU2}$-module structure of $R_H$.
We give what seems to be a new aspect of this
A-D-E correspondence, in section 2.1 below.

Of these, the ones we will need in this paper are
\begin{eqnarray}
R_{\mathbb{A}1}&=&\mathrm{Span}\{r_1'',r_{-1}''\}\,,\nonumber\\
R_{\mathbb{A}3}&=&\mathrm{Span}\{r_1',r_{-1}',r_{\i}',r_{-\i}'\}\,,\nonumber\\
R_{\mathbb{A}5}&=&\mathrm{Span}\{r_1,r_{-1},r_{\omega},r_{-\omega},r_{\omega^2},r_{-\omega^2}\}\,,\nonumber\\
R_{\mathbb{D}4}&=&\mathrm{Span}\{s_0,s_1,s_2,s_3,t\}\,,\nonumber\\
R_{\mathbb{D}5}&=&\mathrm{Span}\{s_0',s_1',s_2',s_3',t',t''\}\,,\nonumber\\
R_{\mathbb{E}6}&=&\mathrm{Span}\{x,x',x'',y,y',y'',z\}\,,\nonumber
\end{eqnarray}
where the notation should be clear (see Figure 1). We write $\omega$ for
$e^{2\pi\i/3}$. The representation Res$^{SU2}_H\sigma$ is $2r_{-1}''$, $r'_{\i}+r_{-\i}'$, $r_{-\omega}
+r_{-\omega^2}$, $t$, $t'$, $y$, respectively. All inductions between  these
finite groups are obtained from:
\begin{itemize}
\item  Ind$_{\mathbb{A}1}^{\mathbb{A}3}$, with $r_1''\mapsto r_1'+r'_{-1}$
and $r''_{-1}\mapsto r'_{\i}+r'_{-\i}\,$;
\item  Ind$_{\mathbb{A}1}^{\mathbb{A}5}$, with $r''_{\pm 1}\mapsto r_{\pm 1}
+r_{\pm \omega}+r_{\pm\omega^2}\,$;
\item  Ind$_{\mathbb{A}3}^{\mathbb{D}4}$, with $r_1'\mapsto s_0+s_{2}$,
$r_{-1}'\mapsto s_1+s_3$, and $r_{\pm \i}'\mapsto t\,$;
\item  Ind$_{\mathbb{A}3}^{\mathbb{E}6}$, with $r_1'\mapsto x+x'+x''+z$, 
$r_{-1}'\mapsto 2z$, and $r'_{\pm \i}\mapsto y+y'+y''\,$;
\item  Ind$_{\mathbb{D}4}^{\mathbb{E}6}$, with $s_0\mapsto x+x'+x''$,
$s_1,s_2,s_3\mapsto z$, and $t\mapsto y+y'+y''\,$;
\item  Ind$_{\mathbb{A}5}^{\mathbb{D}5}$, with $r_1\mapsto s_0'+s_1'$, $r_{-1}
\mapsto s_2'+s_3'$, $r_{\omega},r_{\omega^2}\mapsto t''$, and $r_{-\omega},
r_{-\omega^2}\mapsto t'\,$;
\item  Ind$_{\mathbb{A}5}^{\mathbb{E}6}$, with $r_1\mapsto x+z$, $r_{-1}
\mapsto y'+y''$, $r_{\omega}\mapsto x''+z$, $r_{\omega^2}\mapsto x'+z$,
$r_{-\omega}\mapsto y+y'$, and $r_{-\omega^2}\mapsto y+y''\,$.
\end{itemize}

\bigskip 
\begin{figure}[tb]
\begin{center}
\epsfysize=2.5in \centerline{\epsffile{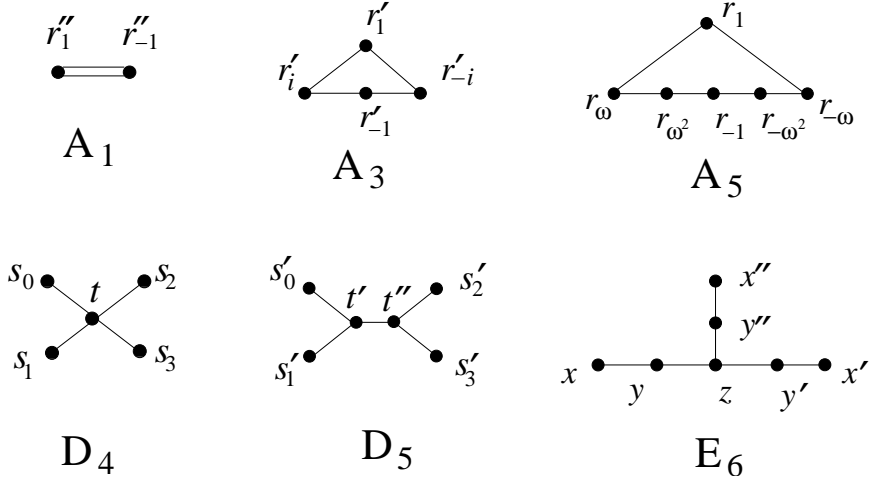}}
\caption{McKay's graphs for some finite subgroups of $SU(2)$}
\end{center}
\end{figure}

The maps between $K$-homology groups tend to be easier to
identify than between $K$-cohomology groups. Also, the answers suggest that
$K$-homology is more natural here (eg. compare $Ver_k(G)={}^\tau K_G^{{\rm dim}\,
G}(G)={}^\tau K^G_0(G)$). For these reasons, we prefer to calculate in
$K$-homology. When the space $X$ is not compact, we must distinguish between
$^{\tau}K_*^G(X)$ and $^{\tau}K_{*,cs}^G$ ($K$-homology with 
compact support): eg. compare
 $^{\tau}K_*^G(\bbR\times X)={}^{\tau}K_{*+1}^G(X)$ with
$^{\tau}K_{*,cs}^G(\bbR\times X)={}^{\tau}K_{*,cs}^G(X)$. Poincar\'e
duality \cite{Ren} relates $K^*$ to $K_{*,cs}$. This yields
two independent ways to compute the $K$-groups. We primarily
use $K_*$, since it permits us to use the six-term exact sequence
(\ref{six-term}).
On the other hand, the $K_{*,cs}$-groups are homotopy invariants.

We have  two main tools for computing twisted equivariant $K$-homology. The first
is obtained by considering the ideal obtained from
a $G$-invariant open subset $X$. Suppose that $\tau^{\prime}$ and
$\tau^{\prime\prime}$ are the restrictions
of $\tau$ on $X$ to $U$ and $X/U$ respectively. Then we have the six-term exact
 sequence for $K_*$:

\begin{equation}\label{six-term}
\begin{array}{ccccc}
    ^{\tau^{\prime}}K^{G}_0(U) & 
    {\longleftarrow} & ^{\tau}K^{G}_{0}(X) &
    {\longleftarrow} & ^{\tau^{\prime\prime}}K^{G}_{0}(X/U)
    \\[3pt]
    \downarrow  & & & & \uparrow
    \\[2pt]
    ^{\tau^{\prime\prime}}K^{G}_1(X/U) & {\lori} &
    ^{\tau}K^{G}_1(X) &{\lori} &^{\tau^{\prime}}K^{G}_1(U)
    \end{array}
    \end{equation}
For $K$-homology with compact support, this fails (consider eg.\ $\mathbb{T}$ with
one point removed). The maps in (\ref{six-term}) are $R_G$-module maps. 

Suppose $X$ is covered by two $G$-invariant open sets, $U_{1}$, and $U_{2}$, and that
$\tau$ restricts to $\tau_{1}$, $\tau_{2}$, and $\tau_{12}$ on $U_{1}$, $U_{2}$ and
$U_{1} \cap U_{2}$ respectively. Then there is the exact Mayer-Vietoris sequence  for $K_*$:

\begin{equation}\label{MV}
\begin{array}{ccccc}
    ^{\tau}K^{G}_0(U_{1} \cup U_{2}) & \stackrel{i^{1}_{*} \ominus i^{2}_{*}
    \,}{\longrightarrow} & ^{\tau_{1}}K^{G}_{0}(U_{1}) \oplus {^{\tau_{2}}K^{G}_{0}(U_{2})}&
    \stackrel{j^{1}_{*} \oplus j^{2}_{*} \,}{\longrightarrow} & ^{\tau_{12}}K^{G}_{0}(U_{1} \cap U_{2})
    \\[3pt]
    \uparrow  & & & & \downarrow
    \\[2pt]
    ^{\tau_{12}}K^{G}_{1}(U_{1} \cap U_{2}) & \stackrel{j^{1}_{*} \oplus j^{2}_{*}\,}{\longleftarrow} &
    ^{\tau_{1}}K^{G}_{1}(U_{1}) \oplus {^{\tau_{2}}K^{G}_{1}(U_{2})} &\stackrel{i^{1}_{*} \ominus i^{2}_{*}}{\longleftarrow}
    &^{\tau}K^{G}_1(U_{1} \cup U_{2})
    \end{array}
    \end{equation}
    where $j^{1}$ and  $j^{2}$ are the inclusions of $U_{1} \cap U_{2}$ in $U_{1}$ and  $U_{2}$
respectively and $i^{1}$ and  $i^{2}$ are the inclusions of
$U_{1}$ and  $U_{2}$ respectively in $U_{1} \cup U_{2}$. For $K$-homology with
compact support, arrows should be reversed. Again, the maps in (\ref{MV}) are
$R_G$-module maps.

Throughout this paper, the groups $H^*_G(X;A)$ denote $\rm{\check C ech}$
cohomology. In computing these cohomology groups, 
we use the relations $H^i_G(\bbR\times X;A)=H^{i}_G(X;A)$ and $H^i_G(\mathbb{T}\times X;A)
\simeq H^{i-1}_G(X;A)\oplus H^{i}_G(X;A)$ for any group $G$ (provided $G$ leaves
$\bbR$ and $\mathbb{T}$ fixed), as well as $H^i_G(X\times G/H;A)=H^i_H(X;A)$ for any
subgroup $H$ of $G$ and space $X$. $H^*_G(M\ddot{o}b;A)\simeq H^*_G(\mathbb{T};A)$ since
$M\ddot{o}b$ (the open M\"obius strip) is a deformation retract from $\mathbb{T}$. Also,
$H^1_G(pt;A)\simeq {\rm Hom}(G,A)$ for any group $G$ and any ring $A$, and
$H^2_G(pt;\bbZ)\simeq G/[G,G]$. The Schur multiplier $H^3_G(pt;\bbZ)$
for any finite subgroup of $SU(2)$ is trivial (this follows from the fact
that they have presentations (\ref{coxmos})
with the same number of generators as relations). Mayer-Vietoris
here becomes
\begin{equation}
0\rightarrow H^0_G(X;A)\rightarrow H^0_G(U_1;A)\oplus H_G^0(U_2;A)\rightarrow
H_G^0(U_1\cap U_2;A)\rightarrow H^1_G(X;A)\rightarrow\cdots\, ,
\end{equation}
for $G$-invariant open sets $U_1,U_2$ covering $X$.
We also compute some $H^*_G(X;A)$ from the spectral sequence (see eg.\
Chapter 1 of \cite{Ht})  associated to the fibration
$X\rightarrow(EG\times X)/G\rightarrow BG$; this has $E_2^{p,q}
=H^p(BG;H^q(X;A))$.

From page 226 of \cite{Ha}, we know $H^*_{SU2}(pt;\bbZ)$ is the polynomial ring
$\bbZ[w_4]$, and page 327 of \cite{Ha} says $H^*_{SU2}(pt;\bbZ_2)\simeq\bbZ_2[w_4]$, where $w_4$
has degree 4. Then \cite{Br} computes the cohomology rings
$H^*_{\bbT}(pt;\bbZ)\simeq\bbZ[w_2]$ and $H^*_{\bbT}(pt;\bbZ_2)\simeq\bbZ_2[
w_2]$, as well as
\begin{eqnarray}
H^*_{SO3}(pt;\bbZ)&\simeq&\bbZ[w_3,w_4]/(2w_3)\, ,\nonumber\\
H^*_{SO3}(pt;\bbZ_2)&\simeq&\bbZ_2[w_2,w_3]\, ,\nonumber\\
H^*_{O2}(pt;\bbZ)&\simeq& \bbZ[w_2,w_3,w_4]/(2w_2,2w_3,w_3^2-w_2w_4)\, ,\nonumber\\
H^*_{O2}(pt;\bbZ_2)&\simeq&\bbZ_2[w_1,w_2]\, ,\nonumber
\end{eqnarray}
where $w_i$ has degree $i$.

Poincar\'e duality, for $X$ a compact manifold,  says (see \cite{Tu,EEK,BMRS})
\begin{equation}
{}^\tau K_i^G(X)\simeq {}^{\tau'} K^{i+{\rm dim}(X)}_G(X)\, ,
\end{equation}
where $\tau+\tau'=(sw_1^G(X),sw_3^G(X))$ (recall that as a group, the twists form
a semi-direct product (\ref{bock}) of $H^1_G$ and $H^3_G$), where the Stiefel-Whitney class
$sw_1^G(X)=0$ iff $X$ is $G$-equivariantly orientable, and $sw_3^G(X)=0$ iff $X$ admits
a $G$-equivariant spin$_c$-structure. A useful fact is that $H^2_G(X;\bbZ_2)
=0$ implies $X$ admits an equivariant spin structure, and hence $sw_3(X)=0$.
Every compact orientable manifold of dimension $\le 4$ admits (though not
necessarily equivariantly) a spin$_c$ structure; however,
$SO(3)$ is compact and orientable and yet $sw_3^{SU2}(SO(3))\ne 0$
for the adjoint action.

In this paragraph let $G=SU(2)$.
Topologically, $G$, $G/O(2)$ and $G/\bbT$ are the 3-sphere $S^3$,
the projective plane, and the 2-sphere $S^2$; $sw_1^G(G/O(2))=1\in\bbZ_2$
but $sw_3^G(G/O(2))$, $sw_1^G(G/\bbT)$ and $sw_3^G(G/\bbT)$ all vanish.
We will also be interested in the {\it spherical manifolds} $G/\Gamma$ for
$\Gamma$ a {\it finite} subgroup of $G$. Since $H^3_G(G/\Gamma;\bbZ)\cong
H^3_\Gamma(pt;\bbZ)=0$ as mentioned above, we know $sw_3^G(G/\Gamma)=0$.
Likewise, $sw_1^G(G/\Gamma)=0$ since $G$, being a Lie group, is
$G$-equivariantly orientable for the translation action, and $G/\Gamma$ will
inherit this. Another way of seeing this is that the $G/\Gamma$ are
rational homology spheres; hence by the Lefschetz fixed point formula,
any orientation-reversing continuous map must have fixed points. Since
any $g\in G$ acts freely on $G/\Gamma$, it must preserve orientation, which
means $G/\Gamma$ is $G$-equivariantly orientable.

When $G$ fixes all of $X$, equivariant $K$-theory can be expressed
in terms of nonequivariant $K$-theory through: \cite{Se}
\begin{equation}\label{trivact}
K_G^*(X)\simeq R_G\otimes_\bbZ K^*(X) \,.
\end{equation} 
In particular, $K^0_G(pt)=R_G$. More generally,
${}^\tau K^0_G(pt)$ is the representation ring ${}^\tau R_G$
described in Definition 4.2 of \cite{FHT}, while ${}^\tau K^1_G(pt)$ is the group
${}^\tau R^1_G$ given in Proposition 4.5 of \cite{FHT}. In particular,
the torsion part of the $H^3$-component of $\tau$ concerns spinors (i.e.
representations of a central extension)
and isn't relevant to the examples considered in this paper;
the $H^1$-component of $\tau$ concerns graded representations.
A grading here is a group homomorphism $\epsilon:G\rightarrow \{\pm 1\}$.
Then $^{\epsilon}R^1_G$ can be defined to be $R_G/\mathrm{Ind}_{G^+}^G
(R_{G^+})$,
where $G^+$ is the kernel of $\epsilon$. The graded representation ring
$^\epsilon R_G$ is the collection of all finite-dimensional $\epsilon$-graded
representations, modulo the supersymmetric ones: a {\it graded representation}
is a $\bbZ_2$-graded vector space $V=V^+\oplus V^-$ carrying a $G$-action, where $G^+$ preserves,
and $G^-=G\setminus G^+$ changes, that grading; a {\it supersymmetric representation}
is a graded representation $V$ with an isomorphism $\alpha:V^\pm\rightarrow
V^\mp$ obeying $g\alpha=\pm\alpha g$ for $g\in G^{\pm}$.
In section 2.1 below we provide a novel interpretation of graded 
representations and of both $^\epsilon R_G$ and $^\epsilon R_G^1$.

For $X$ a compact manifold fixed pointwise by $G$, and $H$ a subgroup of $G$, 
we know 
${}^\tau K^*_G(X\times G/H) ={}^\tau K^*_H(X)$ (see \cite{Se}) and hence
\begin{equation}\label{cancel}
{}^\tau K_*^G(X\times G/H)={}^{\tau'}
 K_{*+{\rm dim}(G)+{\rm dim}(H)}^H(X)
 \end{equation}
  by Poincar\'e duality, for the appropriate twist $\tau'$. If $N$ is a normal 
subgroup of $G$, and $N$ acts  freely on $X$, then from the definition
(\ref{defkhom}) of equivariant K-homology, 
\begin{equation}\label{freenormal}
{}^\tau K_*^G(X)={}^\tau K_*^{G/N}(X/N)\,,
\end{equation}
where we use $H^*_G(X;A)\simeq H^*_{G/N}(X/N;A)$.

In places we will need
infinite-index induction. The usual (Mackey) induction Ind$_H^G(M)=L^2_H(G;M)$ results in
infinite multiplicities; the appropriate notion is 
Dirac induction. One special case of it  we use
is (see Theorem 13 of \cite{La} for a generalisation): \textit{If $T$ is the 
maximal torus of a connected compact simply connected group $G$, and  $\lambda$
is a dominant weight of $G$, then Dirac induction takes a $T$-character
$e^{2\pi \i \lambda}$ to the virtual $G$-representation $({sign\,{w}})
V_{w\lambda-\rho}$ if $w\lambda-\rho\in P_{+}(G)$ for some Weyl element $w$,
and to 0 if no such $w$ exists} (here, $\rho$ is the Weyl vector of $G$,
$V_\mu$ is the $G$-module with highest weight $\mu$, and $P_+$ is the set
of all dominant weights of $G$).
We describe in detail other classes of Dirac inductions in section 2.1,
when we have a better grasp on graded representation rings. 

By contrast, the (closely related) holomorphic induction of Borel-Weil theory
induces the $T$-character $e^{2\pi \i\lambda}$ to 
the $G$-representation $V_{w\lambda}$.
So eg.\ for $G=SU(2)$, Dirac induction takes $\lambda=0$ to 0, $\lambda=1,2,
\ldots$ to the $SU(2)$-representation $\sigma_{\lambda}$,
 and $\lambda=-1,-2,\ldots$ to the virtual $SU(2)$-representation
$-\sigma_{|\lambda|}$. Of course holomorphic induction sends $\lambda$ to
$\sigma_{|\lambda|+1}$.

Near the beginning of  section 1.4  many of these results are put together
into a simple example.

\subsection{Review of notions from CFT}

In this subsection we  review the basic mathematical structures of conformal field theory
(CFT) involved in
this paper. The physical interpretation of some of the following material
is given at the end of this subsection.

Choose any  compact connected simply connected Lie group $G$.
For fixed level $k$ there are finitely many positive energy representations 
$\pi_\lambda$ of the loop group $LG$, parametrised by the highest weight
$\lambda\in P_+^k(G)$. Their characters $\chi_\lambda$ 
define a finite-dimensional unitary representation of $SL(2,\bbZ)$ by
\begin{equation}\label{modtrans}
\chi_\lambda(-1/\tau)=\sum_{\mu\epsilon P_+^k(G)}S_{\lambda\mu}
\chi_\mu(\tau)\, ,\ \ \chi_\lambda(\tau+1)=\sum_{\mu\epsilon P_+^k(G)}
T_{\lambda\mu}\chi_\mu(\tau)\ .
\end{equation}
These matrices $S,T$ are called {\it modular data}, and have special properties
we won't get into. We will often abbreviate the phrase `loop group $LG$ at
level $k$' with $G_k$.

The usual tensor product of Lie algebra modules behaves additively on
the level, but it is possible (using eg. the vertex operator algebra
structure implicit here, or the Kazhdan-Lusztig coproduct) to define a
different one, usually called the fusion product, such that the fusion of
level $k$ modules is again level $k$. The resulting finite-dimensional
fusion algebra is also called the {\it Verlinde algebra} in the mathematical
literature, as it was E.\ Verlinde who 
expressed its structure constants using the matrix $S$ (see (\ref{verl})
below). The Verlinde
algebra $Ver_k(G)$ can be conveniently expressed as a quotient $R_G/I_k$ of
the representation ring $R_G$ (a polynomial algebra) by the
{\it fusion ideal} $I_k$. For example, for $G=SU(n+1)$ and $G=Sp(2n)$, $I_k$ is
the ideal of $R_G$ generated by
all representations of level exactly $k+1$ \cite{Gep}, i.e. all characters 
$ch_\lambda$ where the highest-weight $\la=\sum_{i=1}^n\la_i\Lambda_i$ satisfies
$\sum_{i=1}^n\la_i= k+1$.

Into this context we will often place
 the $r$-dimensional torus $G=T^r=\bbT^r$, but this requires
a little subtlety. The corresponding Lie
algebra, $\mathfrak{g}=\widehat{\bbC}^r$, is not Kac-Moody, and the corresponding
CFT (that of $r$ free bosons living in $\bbR^r$) is not rational. For instance,
its Verlinde algebra is infinite-dimensional. To get finite-dimensionality,
and indeed a fully rational theory to which the formalism of Freed-Hopkins-Teleman
applies, we should proceed as follows.
The role of the level $k$ is an $r\times r$ symmetric, positive-definite
integer matrix -- geometrically, it then corresponds
to the Gram matrix (with entries $b_i\cdot b_j$ for some basis $\{b_i\}$)
of an $r$-dimensional Euclidean lattice $L$. The role of $P^k_+(G)$
is played by $L^*/L$ ($L^*$ is the dual lattice Hom$(
L,\bbZ)$), and the fusion product is $[v][v']=[v+v']$. Algebraically,
this amounts to extending the {\it Heisenberg} vertex operator algebra corresponding
to $\widehat{\bbC}^r$, to the {\it lattice} vertex operator algebra corresponding
to $L$. Physically,  this amounts to bosons living in the $r$-torus
$\bbR^r/L$, at least when $L$ is even.

The modular data and  Verlinde algebra have a direct analogue in
any rational conformal field theory (RCFT) -- the highest weights
$P_+^k(G)$ and cosets $L^*/L$ become the finite set of 
chiral {\it primaries}. Another major class of examples,
in addition to the affine algebras and lattices, comes from
the doubles of finite groups (corresponding to holomorphic orbifolds
-- see eg.\ \cite{CGR}).
More generally, any CFT can be orbifolded by a finite symmetry of $G$.
The most
tractable of these, the holomorphic orbifold, recovers the representation
theory of the quantum double of $G$; its $K$-homological interpretation is developed
in \cite{Ev} and reviewed (and further developed) next subsection. Also very accessible are the 
permutation orbifolds \cite{KS} where $n$ identical RCFTs are tensored 
together and this product is orbifolded by a subgroup of the symmetric group
$S_n$. The chiral primaries of this orbifold are parametrised by pairs
$(O,\psi)$ where $O$ is a $G$-orbit of the primaries of the original theory,
and $\psi$ is an irreducible character of the stabiliser of that orbit.

The modular data and Verlinde algebra are examples of chiral data of the RCFT.
An RCFT
consists of two chiral halves spliced together. The quantity describing
this splicing is the {\it partition function} or {\it modular invariant}
${Z}(\tau)
=\sum_{\lambda,\mu}Z_{\lambda\mu}\chi_\lambda(\tau)\,{\chi_\mu}(\tau)^*$
 of the theory, as it is invariant
under the $SL(2,\bbZ)$-action of (\ref{modtrans}). In terms of
the coefficient matrix ${Z}$, the basic properties of the modular invariant are:
\begin{itemize}
\item $ZS = SZ$, $ZT = TZ\,$,
\item $Z_{\lambda\mu} \in \{0,1,2,3,\ldots\}$\,,
\item $Z_{00} = 1$\,.
\end{itemize}
The third condition comes from
uniqueness of the vacuum  --  we actually get a much richer structure by 
sometimes ignoring this normalisation constraint. In the case of the loop
group at level $k$, $0$ denotes the highest weight $k\Lambda_0$.

The  estimate \cite{BE5} $Z_{\lambda\mu}
\leq Z_{00}d_\lambda d_\mu$, where $d_\lambda = S_{0\lambda}/S_{00}>0$, shows 
that  there are at most finitely many solutions,
for a fixed modular data with a given 
representation of $SL(2,\mathbb{Z})$ and fixed $Z_{00}$.  In the case of
$G=SU(2)$, there are at most three normalised solutions for a fixed level,
according to the A-D-E classification of Cappelli, Itzykson and Zuber
\cite{CIZ}.  A Dynkin diagram is associated to each modular invariant
through the identification of diagonal terms
$\{\lambda:Z_{\lambda\lambda}\neq 0\}=\mathcal{I}$ of $Z$ with the
eigenvalues $\{S_{f\lambda}/S_{00} :\lambda \in
\mathcal{I}\}$ of the corresponding Dynkin diagram, where $f=\Lambda_1+(k-1)
\Lambda_0$ here corresponds to the
fundamental two-dimensional representation of $SU(2)$.  The ${A}_n$ modular
invariant is the
diagonal invariant at level $k=n-1$, $D_n$ is the orbifold or simple current
modular invariant at level $k=2n-4$, and ${E}_6,{E}_7,{E}_8$ are
the  exceptional modular invariants at levels 10,16,28.  For $SU(3)$ the
the analogous modular invariant classification is due to Gannon \cite{G2}.

Let $G$ be a compact connected simply connected Lie group.
Let $H$ be a connected Lie subgroup of $G$ and $\widetilde{H}$ its simply
connected universal covering. Suppose levels $k,\ell$ can be found
so that the central charge $k\,{\rm dim}(H)/(k+h^\vee)$ of $L\widetilde{H}$ at level
$k$ ($h^\vee$ is the dual Coxeter number of $H$) equals that of $LG$ at level
$\ell$. We say that ${H}_k\hookrightarrow G_\ell$ or $\widetilde{H}_k\rightarrow G_\ell$
is a {\it conformal embedding}. The point is that the restriction of
 $LG$-representations to $LH$ involves only finite multiplicities.
Because of this, given a conformal embedding
 $\widetilde{H}_k\rightarrow G_\ell$ and a choice of modular invariant for
$G_\ell$, we get a modular invariant for $\widetilde{H}_k$ by restriction of
the characters.
This is responsible for instance for the $E_6$ and $E_8$ exceptional
invariants in the $SU(2)$ classification.

All conformal embeddings have been classified \cite{BaBo,ScWa} -- eg.\
for simple $G\ne H$, the level $\ell$ must equal 1.
The easiest nontrivial example is when $G$ is simply-laced (i.e. of type A-D-E)
and $H$ is the maximal torus $T=\bbR^r/Q$ (where $r$ is the rank of $G$
and $Q$ is the root lattice).
Then the level $k$ of $LT$ is the Cartan matrix of $G$ 
(in terms of the `Hom' definition of level in section 2.2, this level 
corresponds to the natural embedding of the root lattice $Q$ in the weight
lattice $Q^*={\rm Hom}(Q,\bbZ)$).
In this case, the primaries of $LG$ level 1 are in exact one-to-one
correspondence with those of $LT$ level $k$, so the associated
modular invariant is the diagonal one $Z=I$.

In the subfactor approach to modular invariants, the Verlinde algebra
is represented by endomorphisms $\NXN$ on a III$_{1}$ factor $N$ which are non-degenerately
braided. There are two main sources
of examples -- one from loop groups and the other from quantum doubles of finite systems 
(which are not themselves necessarily
 non-degenerately braided, such as finite groups
or the Haagerup subfactor). Both are relevant for the twisted $K$-homology approach.

Examples in the III$_{1}$-setting appear from the analysis of Wassermann \cite{W}
for $\SUn$, from loop groups. 
Restricting to loops concentrated on an interval
$I\subset S^1$ (proper, i.e.\ $I\neq S^1$ and non-empty),
denote the corresponding subgroup  by:
$\LISUn=\{f\in\LSUn : f(z)=1 \,,\,\, z\notin I \} \,.$
One finds that in each positive energy representation
$\pi_\la$ the sets of operators $\pi_\la(\LISUn)$
and $\pi_\la(\LIcSUn)$ commute, where $I^\rmc$ is the
complementary interval of $I$, using that
$\SUn$ is simply connected. In turn we obtain
a subfactor:
$
\pi_\la(\LISUn)''\subset\pi_\la(\LIcSUn)'\,,$
involving hyperfinite type III$_1$ factors (see \cite{W}).
In the vacuum representation, labelled by $\la=0$,
we have Haag duality in that the inclusion collapses
to a single factor $N(I)=N(I)$.
The inclusion $
\pi_\la(\LISUn)''\subset\pi_\la(\LIcSUn)'$
 can be read as $\la(N(I))\subset N(I)$
for an endomorphism $\la$ of the local algebra $N(I)$,
yielding a system of endomorphisms $\NXN=\{\la\}$ 
 labelled by the positive energy representations.

Other examples arise from quantum doubles of finite systems.
A non-degenerate braiding in quantum double subfactors can be constructed via
three-dimensional TQFT where the crossings are represented 
with tubes. See \cite{EKaw, Iz, EP1} for details.

In either the loop group or quantum double setting, what
we have acting on a factor $N$ are 
 braided endomorphisms $\la\in\NXN$ -- these are required to commute only
up to an adjustment with a unitary $\varepsilon$ = $\varepsilon(\la, \mu)$: 
$\la\mu = Ad(\varepsilon) \mu\la\,. $
Here the family $\{\varepsilon(\la, \mu)\}$ 
can be chosen to satisfy the Yang-Baxter or braid relations, and braiding-fusion relations. 
The endomorphisms will form a system closed under composition:
$[\la][\mu] = \sum N_{\la \mu}^{\nu}[\nu]\,,$
for some multiplicities $N_{\la \mu}^{\nu}$ of positive integers (the fusion
rules). Intertwiners associated to a twist (statistics phase)
 and a Hopf link provide 
matrices   $T$ and $S$  which as  the braiding is non-degenerate, 
 gives a  representation of  the modular group $SL(2, \bbZ)$.
In the loop group setting, the fusion coefficients
$N_{\la \mu}^\nu$ of sectors match exactly the loop group fusion \cite{W}.
By a conformal spin and statistics theorem \cite{FG1,FRS2,GL2}
one can   ensure that the statistics phase (and the
modular $T$-matrix) in our subfactor context coincide
with that in conformal field theory, and hence since the Verlinde $N$ matrices
coincide, so do the modular $S$-matrices.

One modular invariant is always the trivial diagonal invariant: 
$\sum_{\lambda\epsilon \NXN}|\chi_\lambda|^2$\,.  In some sense
\cite{MS2,DV,BE4}, every `physical' modular invariant is
diagonal if looked at properly.  If we restricted our system to a subfactor
$N \subset M$, with systems of endomorphisms on both factors, 
such that the endomorphisms on $M$ decompose to endomorphisms on $N$  as
$\tau=\sum_{\lambda\epsilon\NXN}b_{\tau\lambda} {\lambda}$ according to some branching rules,
 then the diagonal modular invariant should  give an $\NXN$-modular
invariant $\sum_{\tau}|\chi_\tau|^2 =
\sum_{\tau}|\sum_{\lambda}b_{\tau\lambda}\chi_\lambda|^2$, for
$\tau\in \NXN$, $\la\in\MXM$.
In some sense, every modular invariant should look like
this or with a possible twist
$\Sigma_\tau\chi_\tau\chi^*_{\omega(\tau)},$ for a symmetry $\omega$ of
the (extended) fusion rules of $\MXM$.  The problem in
general is then to find such extensions. When there is no twist
present we have what are sometimes called type I invariants:
$Z_{\lambda\mu}=\sum_\tau b_{\tau\lambda}b_{\tau\mu}.$
These are automatically symmetric:
$Z_{\lambda\mu}=Z_{\mu\lambda}$. In the presence of a non-trivial
twist,  we have the type II invariants
$Z_{\lambda\mu}= \sum_\tau b_{\tau\lambda}b_{w(\tau)\mu}.$
These are not necessarily symmetric,  but at least there is a
symmetric vacuum coupling $Z_{0\lambda} = Z_{\lambda 0}$. Not
every modular invariant is symmetric even in this weaker sense (eg. for
$SO(16n)_1$ or for the doubles of some finite groups or
the Haagerup subfactor),
 but every known $\SUn$ modular invariant is 
symmetric in the usual stronger sense.

In practice of course, we would like to start with
the smaller system on $N$ and find an $M$ that realises a given invariant $Z$, 
i.e. inducing instead of restricting.
We induce the system on $N$ to systems on $M$, 
 using the braiding and its opposite to get two systems of endomorphisms 
on $M$, namely $\MXMp$ and $\MXMm$.
The   inclusion $N \subset M$ should be related to the original system 
$\NXN$ in the following sense.  If we consider $M$ as an $N$-$N$ bimodule
and hence as an endomorphism of $N$, then 
$\theta = {_N{M_N}}$ should decompose as a sum of sectors from $\NXN$. 
We can write this canonical endomorphism $\theta$ on $N$ as $\bar\iota\iota$,
where $\iota: N \rightarrow M$ is the inclusion and
$\bar\iota: M \rightarrow N$ its conjugate.
Then $\gamma = \iota\bar\iota$ is the dual canonical endomorphism on $M$.
  Using the braiding $\varepsilon = \varepsilon^+$ or
its opposite braiding $\varepsilon^-$, we can lift an endomorphism
$\lambda$ in $\NXN$ to those of $M$: $\alpha^\pm_\lambda =
\gamma^{-1} Ad(\varepsilon^\pm(\lambda,\theta))\lambda\gamma$. The
maps $[\lambda] \mapsto[\alpha^\pm_\lambda]$ preserve the
operations of conjugation, addition and multiplication of sectors
\cite{xu, BE1}. However, they won't necessarily be 
injective, and $\alpha^\pm_\lambda$ may be reducible. 
 What is important is their intersection $\MXMo = \MXMp \cap \MXMm$.
 When we decompose into irreducibles  we count the
number of common sectors and get a multiplicity:
\begin{equation}
\label{Zisbb}
Z_{\la\mu} = \langle \a^+_\la,\a^-_\mu \rangle \,, \qquad
\la,\mu\in \NXN \,.
\end{equation}
This matrix $Z = [Z_{\la\mu}]$ will be a modular invariant \cite{BEK1, E2}. 
We will shortly find it convenient to drop the normalisation condition
($Z_{00} = 1$), and then we must not  insist that $M$ is a factor. A modular
invariant realised by an inclusion $N \subset M$
 has vacuum multiplicity $Z_{00}$ equal to the dimension of the centre of $M$ \cite{EP1}.
The system $\MXMo$ is non-degenerately braided, and consequently also gives rise to a
 representation of the modular group $SL(2,\bbZ)$ with modular matrices
$S^{ext}$ and $T^{ext}$. The two representations of the modular group
are intertwined via the
 chiral branching coefficients $\langle
\tau,\alpha^\pm_\lambda\rangle,\tau \in \MXMo,\lambda \in \NXN\,,$
i.e. 
 $S^{ext}b^{\pm} = b^{\pm}S$
and $T^{ext}b^{\pm} = b^{\pm}T$\,.
We can decompose the modular invariant as
$Z_{\lambda\mu}=\langle\alpha^+_\lambda,
\alpha^-_\mu\rangle=\sum_{\tau\epsilon \MXMo}b_{\tau\lambda}^+
b_{\tau\mu}^{-} $ or write $Z = {b}^{+t}b^-$.

The associativity of the system of 
endomorphisms $\NXN$ on $N$ yields  a representation
$N_{\la}N_{\mu} = \sum N_{\la \mu}^{\nu} N_{\nu}$ by
commuting matrices $ N_{\la} = [N_{\la \mu}^{\nu} : \mu, \nu\in\NXN]\,,$
describing multiplication by
$\la$. Since  $N_{\bar\nu} = N_{\nu}^{t}$, they are  a family
 of commuting normal matrices and so  can be simultaneously
diagonalised:
\begin{equation}  \label{verl}
N_{\la\mu}^\nu = \sum\nolimits_\kappa
S_{\lambda\kappa}\frac{S_{\mu\kappa}}{S_{0\kappa}}
S_{\nu\kappa}^* \,.
\end{equation}
Remarkably, 
 the diagonalising matrix is the same as the $S$ matrix in the representation
  (\ref{modtrans}) of $SL(2, \bbZ)$.

The action of the $N$-$N$ system $\NXN$ 
on the $N$-$M$ endomorphisms $\NXM$ (obtained by decomposing 
$\{ \iota\lambda=\alpha^\pm_\lambda \iota : \lambda \in \NXN \}$
into irreducibles)
 gives naturally a representation of the same fusion
rules of the Verlinde algebra:
$G_{\la}G_{\mu} = \sum N_{\la \mu}^{\nu} G_{\nu}\,,$
 with matrices $ G_{\la} = [G_{\la a}^{b} : a, b \in \NXM]$.
Consequently, the matrices $G_{\la}$ will be described by the 
 same eigenvalues but with possibly different multiplicities:
\begin{equation}
\label{nimrep}
(G_\lambda)_{ab} = \sum\nolimits_\kappa
\psi_{a\kappa}\frac{S_{\lambda\kappa}}{S_{0\kappa}}
 \psi_{b\kappa}^* \,.
\end{equation}
These multiplicities are given  \cite{BEK2} exactly by the
 diagonal part of the modular invariant: 
${\rm spectrum} (G_{\lambda}) = \{S_{\lambda\kappa}/S_{0\kappa}  
: {\rm with \,multiplicity}\, Z_{\kappa\kappa} \}\,.
$
This is called a {\it nimrep} \cite{BE5} -- a non-negative integer matrix representation. 
Thus a modular invariant realised by a 
subfactor is automatically 
 equipped with a  compatible nimrep whose spectrum is described by the diagonal part of
the modular invariant.  The case of 
$\SUz$ is just the A-D-E classification of Cappelli-Itzykson-Zuber \cite{CIZ}
with the $\NXM$ system yielding the associated (unextended) Dynkin graph. 

The  complexified  finite dimensional fusion rule  algebras spanned by
$\MXMpm$ decompose as  \cite{BEK2}:
\begin{equation}
\label{chiral}
{\rm Furu}(\MXMpm) = \bigoplus_{\tau\epsilon \MXMo}\bigoplus_{\lambda\epsilon\NXN}
{\rm Mat}(b_{\tau\lambda}^\pm)\,.\end{equation} Here $b_{\tau\lambda}^\pm$
are the chiral branching coefficients $\langle
\tau,\alpha^\pm_\lambda\rangle$.
The {\it full} $M$-$M$  {\it system} $\MXM$ is obtained by decomposing 
$\{ \iota \lambda \bar\iota : \lambda \in \NXN \}$ into irreducibles, and
is generated by the  $\pm$-inductions taken together, i.e. both $\MXMpm$ when the 
$\NXN$ braiding is non-degenerate.
The
complexified fusion rule algebra of the full $M$-$M$  system $\MXM$ 
 decomposes as:
\begin{equation}
\label{fullfudec}
{\rm Furu(\MXM)} = \bigoplus_{\la,\mu\epsilon\NXN}
\Mat(Z_{\la\mu})\,,
\end{equation}
and the action of $\NXN$ on $\NXM$ (our nimrep), the Verlinde algebra of
$N$-$N$ sectors on $N$-$M$ sectors only sees the
diagonal part of this representation on:
\begin{equation}
\label{nimfudec}
 \bigoplus_{\la \epsilon\NXN}
\Mat(Z_{\la\la})\,.
\end{equation}
Counting the dimension of the space 
where this acts, we get the number of irreducible $N$-$M$ sectors:
\begin{itemize}
\item $\#\NXM$ = $\sum_{\la} Z_{\la\la}$\,.
\end{itemize}
Moreover, counting the dimensions of the $M$-$M$ sector algebras we get 
\begin{itemize}
\item $\#\MXMpm$ = $\sum_{\tau\la}(b^\pm_{\tau\lambda})^2$\,,
\item $\#\MXM$ = $\sum_{\la, \mu} Z_{\la\mu}^2$\,.
\end{itemize}
These cardinalities  can be read off as  $\tr Z$, $\tr b^{\pm t}b^\pm$ and
  $\tr ZZ^{t}$ respectively.
In the case of chiral locality where $b^+=b^-$, so that
the invariant is type I, we see that  $\#\MXMpm = \tr Z =  \#\NXM\,.$
In fact, $\MXMpm$ can be identified with the nimrep space $\NXM$ 
by mapping $\beta \in \NXN$ to $\iota\beta \in \MXM$, when chiral locality 
holds \cite{BE1}.

Now $ZZ^{t}$ is a modular invariant in its own right,
 satisfying all the axioms except possibly the normalisation.
If a modular invariant $Z$ can be realised by an inclusion $N  \subset M$,
 then there is an associated 
inclusion $N \subset \tilde{M}$, for another algebra $\tilde{M}$, which realises the modular invariant $ZZ^{t}$ \cite{EP1}
such that the full system $\MXM$ for $Z$ is identified with the classifying
CIZ system (nimrep) $\NXMt$ for $ZZ^t$.
Here $Z$ need not be normalised, and in general  $ZZ^{t}$ is certainly not normalised.
The inclusion $N \subset \tilde{M}$ is closely related to the Jones basic construction \cite{J1}
$N \subset M_{1}$ from $N \subset M$. However, 
 it cannot be precisely that as the Jones extension always yields a factor $M_{1}$, if we start from a subfactor
$N \subset M$. What is true is that  $M_{1}$ and $\tilde{M}$ yield the same $N$-$N$ sector in $\NXN$
(i.e. $ _N{M_1}_N \simeq{}_N{{\tilde{M}}_N}$ as $N$-$N$ bimodules), but 
they determine different inclusions $N \subset M_{1}$ and $N \subset \tilde{M}$. An inclusion of $N$ determines by restriction
an $N$-sector but such a sector does not necessarily or uniquely determine an inclusion.
 Taking the central decomposition of $\tilde{M} =  \oplus M_{c}$,
with $M_{c}$ as factors, then
each inclusion $N \subset \tilde{M}$ gives rise to a normalised modular invariant
 $Z_{c}$ so that $ZZ^{t} =  \sum_{c} Z_{c}$ 
decomposes into normalised modular invariants. 
In particular, both
$ZZ^{t}$ and $Z^{t}Z$ decompose into sums of  normalised modular invariants. In this way, the  CIZ graph for $ZZ^t$, namely $\NXMt \simeq \MXM$, decomposes according to
the $\NXN$ orbits.

Any $\SUz$ modular invariant can be realised by a subfactor \cite{O, xu, BE1, BEK2} and a systematic 
or unified formulation  of a subfactor which realises each is given by \cite{E1}:
\begin{equation} \label{theta}
\theta={}_{N}M_{N} = \oplus_{\la} Z_{\la\la} [\la]\,,
\end{equation}
as $N$-$N$ sectors or bimodules where the summation is over even sectors $\la$
(we identify the $SU(2)_k$ highest weight $\la=(k-\la_1)\Lambda_0+\la_1\Lambda_1$
with the Dynkin label $\la_1\in\{0,1,2,\ldots\}$).
Ocneanu \cite{ocn4} has announced that all $\SUd$ modular invariants are realised by subfactors. The
situation for $\SUz$, $\SUd$, and $\SUf$ is reviewed in 
\cite{E2}. 

We consider two $SU(2)$ modular invariants in detail, namely the $D_4$ and
$E_6$ ones. The $D_4$ modular invariant occurs at level 4: 
\begin{equation}
\label{Dfour}
{Z}_{D_4}= |\chi_0+\chi_4|^2+2|\chi_2|^2\,.
\end{equation}
Its diagonal part
$\{\lambda:Z_{\lambda\lambda}\neq0\}$ matches the
spectrum of the (unextended) Dynkin diagram $D_4$, namely
$\{S_{1\lambda}/S_{0\lambda}=2\cos\pi(\lambda+1)/6:\lambda=0,4,2,2\}.$
For this reason Cappelli, Itzykson and Zuber labelled this modular invariant
by the graph ${D}_4$. It can be realised as an orbifold (i.e. simple current)
invariant, but it is more convenient for us to view
it as a conformal embedding invariant  due to  \cite{BN}
which provides the extended system 
diagonalising the invariant. The embedding  $\SUz_4\rightarrow\SUd_1$
 means there is a (two-to-one) mapping of $\SUz$ in $\SUd$ such that
the three level 1 representations of $\SUd$
decompose into level 4 representations of $\SUz$ with finite
multiplicities.  The system $\SUd_1$ has three inequivalent
representations $\{(00),(10),(01)\}$, obeying $\bbZ_3$
 fusion rules. They decompose as
$\chi_{00}=\chi_0+\chi_4,\chi_{10}=
\chi_{01}=\chi_2\,,
$
 so that the $D_{4}$ modular
invariant for $\SUz_{4}$ arises from the diagonal invariant for
$\SUd_1$ : ${Z}_{D_4}=|\chi_{00}|^2 + |\chi_{10}|^2 + |\chi_{01}|^2.$
The conformal embedding gives
us an inclusion of factors:  $N(I) = L_I \SUz  \subset M(I) = L_I \SUd$
using the vacuum representation on $L\SUd$. 
On $N$ we have the system of endomorphisms $\SUz_{4}$ and on $M$ we have
$\SUd_1$.

The canonical endomorphism (\ref{theta}) for this $D_4$ conformal embedding
is given by the vacuum sector
$[\lambda_0] \oplus [\lambda_4]$, the chiral systems decompose as
$ \MXMpm = \{ \as 0 , \aspm 1 , \asx 21 , \asx 22 \} \,,$
and the neutral system is identified with   $\as 0$, $\asx 21$ and $\asx 22$
and obeys $\bbZ_3$ fusion rules. The full system is given by
$ \MXM = \{ \as 0, \asp 1, \asm 1, \asx 21, \asx 22, [\epsilon],
[\eta], [\eta'] \} \,, $
where 
$ \asprod 11 = [\epsilon] \oplus [\eta] \oplus [\eta']\,,$
with statistical dimensions $d_\epsilon=d_{\eta}=d_{\eta'}=1$.
The dual canonical endomorphism is
$[\gamma]=[\id_M]\oplus[\epsilon]$.
Since  $Z^2 = 2Z$ the full system $\MXM$ with cardinality 
$\tr Z^2$ = 8  decomposes as two sheets
which are copies of the Dynkin diagram $D_4$, as in Figure 2: the solid
lines denote multiplication by $\alpha^+_1$, and the dotted ones by
$\alpha^-_1$.

%%%%%%%%%%%% D_4 %%%%%%%%%%%%%
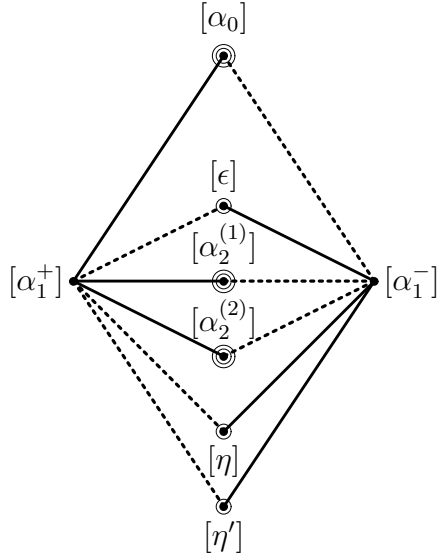
\begin{figure}[tb]
\unitlength 1.0mm
\begin{center}
\begin{picture}(60,80)
%%%%%
\thinlines 
\multiput(10,40)(40,0){2}{\circle*{1}}
\multiput(30,10)(0,10){5}{\circle*{1}}
\multiput(30,10)(0,10){5}{\circle{2}}
\multiput(30,30)(0,10){2}{\circle{3}}
\put(30,70){\circle*{1}}
\put(30,70){\circle{2}}
\put(30,70){\circle{3}}
\Thicklines 
\path(30,70)(10,40)(30,40)  \path(30,30)(10,40)
\dottedline{1.2}(30,70)(50,40)(30,40)
\dottedline{1.2}(50,40)(30,30)
\path(30,50)(50,40)(30,20)  \path(30,10)(50,40)
\dottedline{1.2}(30,50)(10,40)(30,20)
\dottedline{1.2}(10,40)(30,10)
\put(30,75){\makebox(0,0){$\as 0$}}
\put(5,40){\makebox(0,0){$\asp 1$}}
\put(55,40){\makebox(0,0){$\asm 1$}}
\put(30,54){\makebox(0,0){$[\epsilon]$}}
\put(30,45){\makebox(0,0){$\asx 21$}}
\put(30,35){\makebox(0,0){$\asx 22$}}
\put(30,16){\makebox(0,0){$[\eta]$}}
\put(30,6){\makebox(0,0){$[\eta']$}}
\end{picture}
\caption{$D_4$: Fusion graph of $\asp 1$ and $\asm 1$}
\label{D4pm}
\end{center}
\end{figure}

%%%%%%%%%%%% E_6 %%%%%%%%%%%%%
\begin{figure}[tb]
\unitlength 0.6mm
\begin{center}
\begin{picture}(100,140)
%%%%%
\thinlines 
\multiput(10,50)(0,40){2}{\circle*{2}}
\multiput(10,70)(80,0){2}{\circle{4}}
\multiput(90,50)(0,40){2}{\circle*{2}}
\multiput(50,30)(0,80){2}{\circle*{2}}
\multiput(10,70)(80,0){2}{\circle*{2}}
\multiput(50,10)(0,40){4}{\circle*{2}}
\multiput(50,10)(0,40){4}{\circle{4}}
\multiput(50,90)(0,20){3}{\circle{6}}
\Thicklines 
\path(50,130)(10,90)(10,50)(50,90) \path(50,110)(10,70) 
\dottedline{2}(50,130)(90,90)(90,50)(50,90)
\dottedline{2}(50,110)(90,70)
\path(90,90)(50.5,50)(50.5,10)(90,50) \path(90,70)(50.5,30) 
\dottedline{2}(10,90)(49.5,50)(49.5,10)(10,50)
\dottedline{2}(10,70)(49.5,30) 
\put(1,90){\makebox(0,0){$\asp 1$}}
\put(1,70){\makebox(0,0){$\asp 2$}}
\put(1,50){\makebox(0,0){$\asp 9$}}
\put(50,137){\makebox(0,0){$\as 0$}}
\put(50,118){\makebox(0,0){$\asx 31$}}
\put(50,97){\makebox(0,0){$\as {10}$}}
\put(50,57){\makebox(0,0){$[\delta]$}}
\put(56,27){\makebox(0,0){$[\zeta]$}}
\put(50,4){\makebox(0,0){$[\delta']$}}
\put(99,90){\makebox(0,0){$\asm 1$}}
\put(99,70){\makebox(0,0){$\asm 2$}}
\put(99,50){\makebox(0,0){$\asm 9$}}
\end{picture}
\caption{$E_6$: Fusion graph of $\asp 1$ and $\asm 1$}
\label{E6pm}
\end{center}
\end{figure}
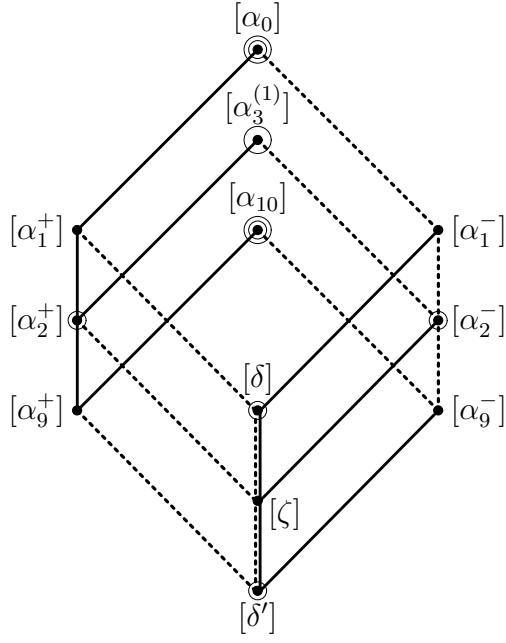

The first  exceptional modular invariant for $SU(2)$
occurs at level 10:
\begin{equation}
\label{Esix}
{Z}_{E_6}=|\chi_0+\chi_6|^2+|\chi_4+\chi_{10}|^2+|\chi_3+\chi_7|^2.
\end{equation}
Its diagonal part $\{\lambda:Z_{\lambda\lambda}\neq0\}$ matches the
spectrum of the Dynkin diagram $E_6$, namely
$\{S_{1\lambda}/S_{0\lambda}=2\cos\pi(\lambda+1)/12:\lambda=0,6,4,10,3,7\}.$
 This is obtained from the conformal embedding   $SU(2)_{10} \rightarrow
Sp(4)_1$. The system $Sp(4)_1$ again has three inequivalent
representations: the vacuum (00), vector $(01)$ and spinor $(10)$; they
reproduce the Ising fusion rules. Restricting from $Sp(4)$ to $SU(2)$, they
decompose as
$\chi_{00}=\chi_0+\chi_6,\chi_{01}=\chi_4+\chi_{10},
\chi_{10}=\chi_3+\chi_7
$
 so that the $E_{6}$ modular
invariant for $SU(2)_{10}$ arises from the diagonal invariant for
$Sp(4)_1$ : ${Z}_{{E}_6}=|\chi_{00}|^2 + |\chi_{01}|^2 + |\chi_{10}|^2.$
The conformal embedding gives
us an inclusion of factors:  $ N(I) = L_I SU(2)  \subset M(I) = L_I Sp(4)$
using the vacuum representation on $LSp(4)$. 
On $N$ we have the system of endomorphisms $SU(2)_{10}$ and on $M$ we have
$Sp(4)_1$.

We have \cite{BE1} the chiral systems 
$ \MXMpm = \{ \as 0 , \aspm 1 , \aspm 2 ,
\asx 31 , \aspm 9 , \as {10} \} \,,$
where 
$\as 3 = \asx 3 1 \oplus \as 9 \,, \
\as 4 = \as 2 \oplus \as {10} \,, \
\as 5 = \as 1 \oplus \as 9\,, \
\as 6 = \as 0 \oplus \as 2 \,, 
\as 7 = \as 1 \oplus \asx 3 1 \,.$ The neutral system
$\MXMo = \{ \as 0 ,  
\asx 31 , \as {10} \}$
with its Ising fusion rules are identified with the vacuum, spinor
and vector representations of $Sp(4)$ at level 1 respectively.  The full system is:
\[ \MXM  = \{ \as 0 , \asp 1 , \asm 1 ,  \asp 2 , \asm 2 ,
\asx 31 , \asp 9 , \asm 9 , \as {10} , [\delta] ,
[\zeta], [\delta'] \} \,, \]
where $[\delta]=\asprod 11$, $[\zeta]=\asprod 12 = \asprod 21$
and $[\delta']=\asprod 91 = \asprod 19$. 
The dual canonical endomorphism decomposes as
$[\gamma]=[\id_M]\oplus\asprod 11$, whilst $[\theta]=[\lambda_0]\oplus[\lambda_6]$.
Since  $Z^2 = 2Z$, the full system $\MXM$ with cardinality 
 $\tr Z^2$ = 12  decomposes as two sheets
which are copies of the Dynkin diagram $E_6$, as in Figure 3.
%Sigma restriction $\sigma_\beta = \bar\iota\beta{\iota}$ to
%take $M$-$M$-sectors to $N$-$N$ sectors namely
% $\sigma_b=\theta= \lambda_0+ \lambda_6
%,\sigma_v = \lambda_4 + \lambda_{10}, \sigma_s =
%\lambda_3+\lambda_7$

Physically \cite{PZ,E3}, $\NXN$ concerns the chiral bulk data (eg. Verlinde algebra), 
$\NXM$ the boundary data (eg. nimrep$=$annulus partition function),
and the full system $\MXM$ 
the defects. In particular, the endomorphisms $\lambda\in\NXN$ label the
primaries, i.e. the irreducible modules of the chiral algebra $\mathcal{A}$
of the
theory; the $a\in\NXM$ label the boundary states; and the $\alpha_i\in\MXM$
label defect lines. The endomorphisms of the neutral system $\MXM^0=\MXM^+\cap
\MXM^-$
label irreducible modules of the chiral subalgebra preserved by the boundary
conditions. The matrix $\psi$ diagonalising the nimrep (\ref{nimrep})
relates boundary states to Ishibashi states. For the special case of
modular invariants of $LG$ (or its affine algebra $\mathfrak{g}$) associated
to symmetries $\omega$ of the corresponding unextended Dynkin diagram (eg. charge
conjugation), this data has a clear  Lie theoretic interpretation \cite{GaGa}:
boundary states are labelled by integral highest-weights for the twisted
affine algebra $\mathfrak{g}^\omega$, or equivalently $\omega$-twisted 
$\mathfrak{g}$-representations;
the nimrep coefficients are twisted fusion coefficients; $\psi$ describes
how $\mathfrak{g}^\omega$-characters transform under $\tau\mapsto -1/\tau$; and
the exponents are highest-weights of another twisted algebra, called the
orbit algebra. The categorification of bulk and boundary conformal field theory
(see eg. the review article \cite{RFS}) owes much to the subfactor picture.
In particular,  the starting point there is the category $\mathcal{C}$ of
$\mathcal{A}$-modules, together with the object $A$ corresponding to
the canonical endomorphism $\theta$ of   (\ref{theta}) -- $A$ will be a
special symmetric Frobenius algebra of $\mathcal{C}$. In this language the
boundary conditions are $A$-modules and the defect lines are $A$-$A$-bimodules.
The nimrep and modular invariant are constructed from $A$ using the analogue there of
$\alpha^\pm$-induction.

\subsection{Review of the K-homological approach to CFT}

This subsection reviews the Freed-Hopkins-Teleman realisation of the Verlinde 
algebra $Ver_k(G)$, for $G$ a Lie group. It then reviews the analogous
construction (and extension) for finite groups, described in \cite{Ev}, and
concludes by describing the analogue for finite groups of conformal embeddings.
The success of this  finite group story is crucial motiviation for this paper.

Let $G$ be simple and simply connected, and $k$ any integral level. The
main result of Freed-Hopkins-Teleman \cite{FHT,FHTi,FHTii,FHTlg} is that the 
Verlinde algebra $Ver_k(G)$
can be realised as the $K$-homology group ${}^{k+h^\vee}K_0^G(G)$, where $G$ acts
on $G$ adjointly and $k+h^\vee\in\bbZ\simeq H^3_G(G;\bbZ)$ (here, $H^1_G(G;
\bbZ_2)=0$). An elegant proof of this is given in \cite{M}. 
The twist is crucial for finite-dimensionality: eg. \cite{Bry} computes that
the untwisted $K^G_*(G)$ is a free $R_G$-module of rank $2^{{\rm rank}\,G}$. 

For example, consider $G=SU(2)$ on $SU(2)$. Then eg.\ by spectral sequences
we obtain $H^1_G(G;\bbZ_2)=0$ and
$H^3_G(G;\bbZ)\simeq \bbZ$, and we identify the twist $\tau$ with the shifted
level $k+2$. The orbits of $G$ on $G$ come in two kinds: the fixed points
$\{\pm I\}$, and the generic points $gen=\bbR\times G/\bbT$ with stabiliser
$\bbT$. The six-term relation (\ref{six-term}) tells us how to glue together
the $K$-homology of the fixed points to those of $gen$: it becomes
\begin{equation}\label{su2su2}
\begin{array}{ccccc}
    ^{\tau^{\prime}}K^{G}_0(\bbR\times G/\bbT) & 
    {\longleftarrow} & ^{k+2}K^{G}_{0}(G) &
    {\longleftarrow} & ^{\tau^{\prime\prime}}K^{G}_{0}(\pm I)
    \\[3pt]
    \downarrow  & & & & \uparrow\beta
    \\[2pt]
    ^{\tau^{\prime\prime}}K^{G}_1(\pm I) & {\lori} &
    ^{k+2}K^{G}_1(G) &{\lori} &^{\tau^{\prime}}K^{G}_1(\bbR\times G/\bbT)
    \end{array}    \ .
    \end{equation}
Using the simple results on equivariant cohomology collected at the end of
section 1.2, we immediately find that the relevant cohomology groups on
the fixed points and $gen$ all vanish: eg.\
 $H^3_G(\bbR\times G/\bbT;\bbZ)\simeq H^3_{\bbT}(pt;\bbZ)=0$. This means that the twists
 $\tau',\tau''$ in (\ref{su2su2}) vanish, and the level $k+2$ can only
 appear in the maps. Of course $K^G_0(\pm I)=K_G^0(pt)\oplus K_G^0(pt)=
 R_G\oplus R_G$, while $K^G_1(\pm I)=0$. Likewise, $K_*^G(\bbR\times G/\bbT)
 \simeq K_*^G(G/\bbT)\simeq K_*^\bbT(pt)$, using (\ref{cancel}), and hence
$K_1^G(gen)= R_\bbT$ and $K_0^G(gen)=0$. So all that remains is to
identify the map $\beta:R_\bbT\rightarrow R_G\oplus R_G$, which we know
should involve $k$. The answer is:
 $\beta$ will send the polynomial $p(a)\in R_\bbT$ to D-Ind$_{\bbT}^G
 (p(a),a^{k+2}p(a))$. Eq.(\ref{su2su2}) says
 $^\tau K_1^G(G)$ is the kernel of  $\beta$ while $^\tau K_0^G(G)$ is its
 cokernel.

The presence of Dirac induction in $\beta$ is clear, but it may be hard to
anticipate this prefactor $a^{k+2}$, without
knowing the underlying bundle (which we describe below in section 2.2).
Locally about both $\pm I$ the bundle is trivial; $a^{k+2}$ is the
relative twist picked up when comparing these trivialisations.
This simple example is a baby version of the other much more complicated
calculations we do elsewhere in this paper.
This same example was worked out in eg. Example 1.7 of \cite{FHTi},
using Mayer-Vietoris and $K_{*,cs}$, with the same result (the answers must
agree since the space $G$ is compact). A very explicit yet elegant calculation of
$K_*^G(G)$ for any compact simple $G$ was done in \cite{M} using the
spectral sequence of \cite{Se2}.

When $G$ is not
simply connected, the situation is a little more complicated: there will be
torsion in both $H^3_G(G;\bbZ)$ and $H^1_G(G;\bbZ_2)$. The calculation for
$G=SO(3)$, and all classes of twists, was  worked out in \cite{FHT}; for
the appropriate twist $\tau$, $^\tau K_0^G(G)$ is again a Verlinde algebra,
namely the extended Verlinde algebra of the type $I$ $SU(2)$ modular invariants
of ${D}$-type. However, for other values of the twist $\tau$, there
is nontrivial $K$-homology $^\tau K_*^{G}(G)$ which doesn't have an obvious
CFT interpretation. We return to this in section 7.

A concrete calculation is given 
in Proposition B.5 of \cite{FHT}, where the extended Verlinde algebra for the
simple current modular invariant at $SU(2)$-level $k=4n$ is realised as
${}^{(-,-,1+2n)}K^1_{\overline{G}}(\overline{G})\simeq
{}^{++}R(k)\oplus\bbZ$ where $R(k)$ is the Verlinde algebra
of $SU(2)$ at level $k$, spanned by the irreducible representations $\sigma_i$,
$1\le i\le k+1$, ${}^\pm R(k)$ means nonspinors/spinors,
and ${}^{+\pm}R(k)$ means to identify weights in the same $J$-orbit in
${}^\pm R(k)$ (i.e.\ ${}^{++}R(k)=\bbZ[\sigma_1]\oplus[\sigma_3]\oplus
\cdots\oplus \bbZ[\sigma_{k-1}]$ and
${}^{+-}R(k)=\bbZ [\sigma_2]\oplus[\sigma_4]\oplus\cdots\oplus[\sigma_{n}]$). The extra $\bbZ$ comes from the graded
representation $\delta-1$, and corresponds to `resolving' the fixed point
$[\rho_{n+1}]$. The $K$-groups ${}^{(-,\pm,n+1)}K^0_{\overline{G}}(\overline{G})$ are both trivial, while
${}^{(-,+,n+1)}K^1_{\overline{G}}(\overline{G})={}^{+-}R(k)$.

A major clue as to extensions of Freed-Hopkins-Teleman is provided by
considering finite groups $G$.  \cite{Ev} provides the 
 $K$-homological description for modular invariants associated to
the modular data arising from the quantum double of $G$.
We'll review this description in the next few paragraphs.

Take a $G$-kernel on a
factor $M$, that is, a homomorphism from $G$ into the outer 
automorphism 
group ${Out}(M)$ of $M$  (namely the automorphism group of $M$ 
modulo the inner automorphisms). If $\nu_{g}$ in {Aut}$(M)$ 
is a choice of representatives for each $g$ in $G$ of the $G$-kernel, 
then
$ \nu_{g}\nu_{h} = {Ad} (u(g,h)) \nu_{gh}\,,$
for some unitary $u(g,h)$ in $M$, for each pair $g,h$ in $G$. 
% We can 
% assume the normalisation $\nu_{e} = id_{M}$, $u(g,e) = u(e,g) = 
% 1_{M}$, for all $g$ in $G$, where $e$ is the unit of the group. 
By associativity of 
$\nu_{g}\nu_{h}\nu_{k}$, we have a scalar 
$\omega(g,h,k)$ in $\bbT$ such that 
% \begin{eqnarray} 
$u(g,h)u(gh,k)=w(g,h,k)\nu_g(u(h,k))u(g,hk)\,.$ 
% % \end{eqnarray} 
%  i.e. $\omega = \partial_{\nu} u$, the 
% $\nu$-coboundary 
% of $u$.
A standard computation  shows that $\omega$ is a 3-cocycle
in  $Z^3(G;\bbT)$.
% :
% \be 
% \omega(g,h,k)\omega(g,hk,l)\omega(h,k,l)=\omega(gh,k,l)\omega(g,h,kl)\,, 
% \quad g,h,k,l\in G\,.
% \ee 
% Every element of $Z^{3}(G,\bbT)$ arises in this way from some  
% $G$-kernel \cite{Suth, J1} (see also \cite{KT, W2}).
% 
% where the correction term $w(,,)$ is circle valued as $M$ has trivial center.
% Associativity means that $w$ is a 3-cocycle in $Z^{3}(G,\bbT)$.
Conversely, any 3-cocycle arises in this way for some 
$G$-kernel. One can even choose $M$ to be hyperfinite \cite{J1},
 but for our purposes
any realisation will do -- the simplest being with
free group factors \cite{Suth}. 
Now in the
tube algebra approach of Ocneanu (see \cite{EKaw}) to
 the quantum double of $G$, one considers the
space of intertwiners
${\rm 
Hom}(\nu_{h}\nu_{a} , 
\nu_{hah^{-1}}\nu_{h}) = T(a, h)\,.$ This is a line bundle, and  multiplicativity
of these line bundles means that 
$T(hah^{-1},h') \otimes T(a,h) \simeq T (a, h'h)$ \cite{Ev}. This
gives a projective representation of the {\it groupoid}  $G \rtimes G$
(not to be confused with the semi-direct product of groups)
of $G$ acting on itself by conjugation, and consequently an element of
$Z^{2}(G \rtimes G; \bbT )$, which can be identified with 
the equivariant 2-cocycles $Z^{2}_{G}(G; \bbT )$. Thus, associated to
$\omega\in Z^3(G;\bbT)$ is a cocycle in $Z_G^2(G;\bbT)$\,.

Now by definition, 
the equivariant cohomology $H^n_G(G; \bbT)$ is given by $H^n((G \times EG)/G;\bbT)$.
However a model
for the classifying space $BG= EG/G$ is given by simplices associated to $n$-tuples $(g_1, \dots ,
g_n)\,,$
with edges given by $g_1, g_2 \dots \,.$
 Similarly a model for $(G \times EG)/G$ is given by $n$-simplices associated to
 $n+1$-tuples $(g, g_1, \dots , g_n)\,,$ with $g$ associated to the origin, and $g_1$ to the first
edge
to the next vertex $g_1gg_1^{-1}$, $g_2$ to the next edge to the next vertex
$g_2g_1gg_1^{-1}g_2^{-1}$ etc. 
This allows us to identify $H^n_G(G; \bbT)$ with $H^n(G\rtimes G; \bbT)$, for
that groupoid $G \rtimes G$. Hence, given $\omega$ we get a 2-cocycle in
$H^2_G(G; \bbT) \simeq \oplus_{[t]}H^2(BC_G(t); \bbT)$, where the sum is 
over all conjugacy classes, and $C_G(t)$ is the centraliser.

Once we have an element of $H^{2}_{G}(G; \bbT) = H^{3}_{G}(G;\bbZ)$,
we can construct an equivariant bundle of compacts over $G$. 
However  the 
$K$-theory of the $C^{*}$-algebras of the space of sections, does not in general lead to the
the twisted equivariant $K$-group $^{\omega}K^{0}_{G}(G)$ where
$G$ acts on itself by conjugation. The correct formulation of this twisted $K$-theory is not through
the $C^*$-algebra of the space of sections but through the representation theory of the twisted
quantum double.
If $\omega$ is a 3-cocycle on $G$, and $\alpha$ the corresponding
3-cocycle on $G \times G$, given by the difference of the two pullbacks
of $\omega$ on the factors,
 then the Verlinde algebra
is described as the equivariant $K$-homology group $^{\omega}K_{0}^{G}(G) 
\simeq {}^{\alpha}K_{0}^{G \times G}(G \times G)$. Here in the first formulation,
$G$ acts on $G$ by conjugation, and in the second, we have diagonal actions
of $G$ on $G\times G$ on the left and right.
In the second formulation, a precise description of  an element of the Verlinde algebra is  as a
vector bundle
$V$ over $G \times G$, with left and right actions of $G$ diagonally on the base space which act on
the bundle in compatible  way according to the 3-cocycle $\alpha$: 
\begin{eqnarray}
%%% related to Karoubi finite type???
    (h_{1}h_{2})w = \overline{\alpha}(h_{1}, h_{2}, g) 
    (h_{1}(h_{2}w)) \,,  \\
    w(k_{1}k_{2}) = \alpha(g, k_{1}, k_{2})  
      (w k_{1})k_{2} \,, \\  
       h (w k) = \alpha(h , g, k) 
	   (h w)k \,, 
    \end{eqnarray}
where $h_{1},h_{2},h, k_{1},k_{2},k \in G$, and
$w = w_{g} \in V_{g}$, the fibre over $g\in G \times G$. 
The transgression map from $H^{3}(G; \bbT)$ to $H^{3}_{G}(G;\bbZ)$ can be zero, and so we need to
keep track of
where the element of $ H^{3}_{G}(G;\bbZ)$ really comes from in $H^{3}(G; \bbT)$.

The product $V\otimes_G W$ in this Verlinde algebra can be naturally found 
as follows. Given $G$-$G$ bundles $V$ and $W$, divide the tensor product 
$V \otimes W$  by the relation:
\begin{equation}
    v_{a}k \otimes w_{b} = \alpha(a,k,b)  v_{a} \otimes k w_{b}
\end{equation}
and then push-forward under the product map
$(G\times G) \times (G \times G)  \rightarrow G \times G$
 to obtain a bundle over $G\times G$ we'll denote $V \otimes_{G}W,$ with
fibres $(V \otimes_{G} W)_{g} = \oplus_{ab=g} V_{a} \otimes  W_{b}$. 
Then $V \otimes_{G} W$ becomes a $G$-$G$, $\alpha$-twisted 
bundle under the natural actions:
\begin{eqnarray}
    h(v_{a} \otimes w_{b}) = \alpha(h,a,b)  hv_{a} \otimes 
    w_{b}\,,\\
    (v_{a} \otimes w_{b})l = \overline{\alpha}(a,b,l)  v_{a} \otimes w_{b}l
\,.\end{eqnarray}
The  braiding is given by 
$   v_{(e,b)} \otimes w_{(a,e)} \mapsto  w_{(a,e)} \otimes v_{(e,b)}$
together with $G\times G$-equivariance.

%Taking the 
%$K$-theory of the $C^{*}$-algebras of the space of sections, we find
%the twisted equivariant $K$-group $^{\omega}K^{0}_{G}(G)$ where
%$G$ acts on itself by conjugation.  

By analogy with the loop group case, the parameter $\omega$ is regarded as
the level. The map $H^3(BG;\bbT)\rightarrow H^3_{G\times G}(G\times G;\bbZ)$,
$\omega\mapsto \alpha$, constructed above is just the transgression $H^4(BG;\bbZ)\rightarrow
H^3_G(G;\bbZ)$ discussed in the introduction, as $H^*(BG;\bbT)\simeq
H^{*+1}(BG;\bbZ)$ for finite groups, and $H^*_{G\times G}(G\times G;\bbZ)\simeq
H^*_G(G;\bbZ)$ for any group. To simplify the discussion now, we'll consider
the case of trivial level $\omega$.

A modular invariant is described by a subgroup $H$ of $G \times G$, and the
 simplest possible situation is when the subgroup contains the diagonal
$\Delta = \{(g,g): g \in G \}$, so $\Delta \subset H \subset G \times G$.
Then  $N = \{ab^{-1} : (a,b)
\in H \}$ is a normal subgroup of $G$ and $G \times G / H$ is identified with $G/N\,.$
If $\pi : G \rightarrow G/N = L$ is the quotient map, and $\theta : G \times G \rightarrow G$
is $ (a,b) \mapsto ab^{-1}\,,$ then $H = {\rm ker}\,   \pi \theta \subset G \times G\,.$

It was remarked cryptically in \cite{CGR} that 
the surjective homomorphism $\pi : G \rightarrow G/N$ is the finite group analogue
of the conformal embedding of Lie groups discussed in section 1.3. This
can be understood as follows. The full system is identified
with the equivariant $K$-theory $K^0_{H\times H}(G\times G)$, with 
an irreducible equivariant bundle is described  by pair consisting of a double coset $HgH$ in
$H\backslash(G \times G)/H$ and an irreducible representation of 
the stabiliser subgroup
$H \times_{g}H = \{(h,k) \in H \times H : hg=gk\}$
which is isomorphic to 
$H{\cap}\, gHg^{-1}$\,.

The neutral system $\MXMo = \MXMp \cap \MXMm\,$, where
$\MXMpm$ are the $\alpha^{\pm}$-induced systems,  can be computed directly as follows
and identified with  $K^0_{L}(L) \simeq K^0_{\Delta(L) - \Delta(L)}(L \times L)\,.$
For ease of notation, let us take $G$ abelian and consider $\alpha$-induction:
$$\alpha^{\pm} : K_G^0(G) \rightarrow K_{H \times H}^0(G \times G)\,.$$
A primary field in  $K_G^0(G)$ is labelled by $[a,\pi]$ where $a \in G, \pi \in \widehat G$
(a conjugacy class and a representation of the stabiliser).
Then $\alpha$-induction is described by 
\begin{eqnarray}
 \alpha^+ : [a,\zeta] &\mapsto& [(a \times 1)H,(\zeta \times 1)|_H]\,,\nonumber\\
\alpha^- : [b,\psi]& \mapsto& [(1 \times b)H, (1 \times \psi)|_H]\,.\nonumber
\end{eqnarray}
For $\alpha^+[a,\zeta] =  \alpha^-[b,\psi]$, then \cite{Ev} we need $ab^{-1} \in N\,,$
$\zeta = \psi$ and $\zeta|N = 1\,.$ So the primary fields of $\MXMp \cap \MXMm\,$
are described by the cosets of $G/N$ and the representations of $G/N$, i.e. the
quantum double of $G/N = L$, which is $K^0_{L}(L)\,.$ 
The classifying systems $\NXM$ and $\MXN$ are identified with 
$K^0_{\Delta \times H}(G\times G)$ and 
 $K^0_{H \times \Delta}(G\times G)$ respectively and naturally with each other and
with the induced systems $\alpha^\pm(K^0_G(G))$.

The modular invariant is given
through $\alpha$-induction as $Z_{\la \mu} = \langle \alpha^+_{\la},  \alpha^-_{\la} \rangle\,,$
or through $\sigma$-restriction as $\sum_{\tau} b_{\la\tau} b_{\mu \tau}$ where the
 branching coefficient is given as $b_{\la \tau} = \langle \alpha^{\pm}_{\la},  \tau \rangle = 
\langle \la, \sigma_{\tau} \rangle \,.$ If $\tau = [k, \phi], \ell \in L, \phi
 \in \hat {L}$
is a primary field in the neutral system  $K_{L}^0(L)$, then its 
$\sigma$-restriction is given by 
$$\sigma_\tau = \sum_{g \epsilon \pi^{-1}(\ell)} [g, \phi \pi]\,.$$ 
Alternatively, in terms of vector bundles, take $V$ in $Bun^{L}(L)$, then by using the 
morphism $\pi: G \rightarrow L$, we get an equivariant bundle $W$ in $Bun^{G}(L)\,,$
and hence the pullback $\pi^*W$ in $Bun^{G}(G)$. Then $V \mapsto
\sigma_V = \pi^*W$ is the map
$Bun^{L}(L) \rightarrow Bun^{G}(G)$, which yields the modular invariant $Z = b^tb$.

\section{Gradings and bundles}

\subsection{Gradings and induction}

An $H^1$-twist involves graded representations -- we briefly mentioned
these in section 1.2.
In this subsection we rewrite sections 4.1-4.7 of \cite{FHT}, by 
interpreting graded representation rings etc
very concretely in terms of ordinary representations of an index-2 subgroup. 
To our knowledge this interpretation, which seems
conceptually simpler and more amenable to computations than that given in
\cite{FHT}, is new. We conclude the subsection with examples of
Dirac induction.

Let $G$ be a compact group, and $H$ an index-2 subgroup. Let $g\in G\setminus
H$. If $\rho$ is an irreducible $G$-representation, then one of the following holds (see eg.\ section III.11 of \cite{Sim}):

\begin{itemize}

\item[type ${2\atop 1}$:] $\overline{\rho}:=\,$Res$_H^G\rho$ is irreducible; the character
$\chi_\rho$ is not identically 0 on $G\setminus H$; there is an irreducible
$G$-representation $\rho'$ with character $\chi_{\rho'}(h)=\chi_\rho(h)$,
$\chi_{\rho'}(gh)=-\chi_{\rho'}(gh)$ for all $h\in H$; Ind$_H^G\overline{\rho}
=\rho\oplus\rho'$.

\item[type ${1\atop 2}$:]  Res$_H^G\rho$ has irreducible decomposition $\rho_1\oplus
\rho_2$, where $\chi_{\rho_2}(k)=\chi_{\rho_1}(gkg^{-1})$ for all $k\in G$;
$\chi_\rho$ is identically 0 on $G\setminus H$; Ind$_G^H\rho_i=\rho$.

\end{itemize}

A \textit{graded irreducible $G$-representation} is an irreducible $G$-representation
$\rho$ of type ${1\atop 2}$ and a choice (the $\bbZ/2\bbZ$-grading of \cite{FHT})
of calling one of $\rho_i$ `$\rho_+$' and the other `$\rho_-$'; we denote
this graded representation $\rho_+\ominus\rho_-$. The group-homomorphism $\epsilon:G\rightarrow\{\pm 1\}$
with kernel $H$ is an element of $H^1_G(pt;\bbZ_2)$. Then $^\epsilon R_G={}
^\epsilon K^0_G(pt)$
is the span of these $\rho_+\ominus\rho_-$ (where $-(\rho_+\ominus\rho_-)=\rho_-
\ominus\rho_+$).
Similarly, $^\epsilon R^1_G={}^\epsilon K^1_G(pt)$ consists of all possible
sums of irreducible
${2\atop 1}$-representations, modulo the sums of all combinations
$\rho\oplus\rho'$ -- as $\rho$ is then equivalent to $-\rho'$, we will
write the class containing $\rho$ as the anti-symmetrisation   $\rho^-=(\rho
-\rho')/2$.
These representation rings for $H^3$-twists $\tau$
are defined analogously.

The restriction map $^\epsilon$Res$^G_H:{}^{\epsilon,\tau} R_G\rightarrow {}^\tau R_H$
takes $\rho_+\ominus\rho_-$ to $\rho_+-\rho_-$, while induction $^\epsilon
{\rm Ind}^G_H:
{}^\tau R_H\rightarrow {}^{\epsilon,\tau}R_G$ takes $\rho_+$ to $\rho_+\ominus
\rho_-$ (if the usual induction Ind$^G_H\rho_+$ is type ${1\atop 2}$) and to 0 otherwise.
Graded restriction $^\epsilon$Res$'{}_G:{}^{\epsilon,\tau}R_G^1\rightarrow
{}^\tau R_G$ takes $\rho^-$ to $\rho-\rho'$, while graded induction  
$^\epsilon {\rm Ind}'{}_G:{}^\tau R_G\rightarrow {}^{\epsilon,\tau}R_G^1$ takes
$\rho$ to $\rho^-$. Frobenius reciprocity becomes the exact
sequences (4.7) of \cite{FHT} -- those equations are special cases of
the sequences (4.14) and (4.15) of \cite{FHT} on page 14, which in turn
are a special case of the exact sequences (4.2)  of \cite{UHaag}.  
Because $^\epsilon$Res and $^\epsilon$Res$'$ are injective, the $R_G$-module
structure of $^{\epsilon,\tau}R_G$ and $^{\epsilon,\tau}R_G^1$ is obtained
by restricting to $^\tau R_H$ and $^\tau R_G$ respectively.

This connects nicely with a description of graded $K$-theory due to Pimsner.
Let $A$ be a (finite-dimensional)  graded $C^*$-algebra. Then according to 
\cite{Pimsner} the graded $K$-theory is described by graded traces. More 
precisely $K_0(A), K_1(A)$ are generated by graded traces supported on the even
and odd subspaces of $A$. A graded trace $\phi$ on $A$ is a linear map which 
vanishes on the graded commutators $[x,y] = xy - 
(-1)^{\partial x\partial y}yx$, where $x, y$ are homogeneous, and $\partial$ 
denotes the grading. Let $G$ be a finite group, and $\alpha$ an element of 
$H^1_G(pt; \bbZ_2) = {\rm Hom}(G,\bbZ_2)\,.$
Then $^{\alpha}K^*_G(pt) \simeq {}^{\alpha}K_*(C^*(G))\,,$ where $A = C^*(G)$ 
is graded by $\alpha$. Thus $^{\alpha}K^*_G(pt)$ is described by graded 
characters on $G$ supported on the degree $i$ elements in
$G_i$. A graded character of degree $i$ is a map $\chi: G_i \rightarrow 
\bbC\,,$ such that for $i =0\,$:
$\chi(x_{+}y_{+})=\chi(y_{+}x_{+})$ (for $x_{+},y_{+}  \in G_0$)
and $\chi(x_{-}y_{-})=- \chi(y_{-}x_{-})$ (for $x_{-},y_{-} \in G_1$);
and for $i=1\,$:
$ \chi(x_{+}y_{-}) = \chi(y_{-}x_{+})$ for $x_{+} \in G_0\,, y_{-} \in G_{1}\,.  $
This is reminiscent of Section 4.8, page 13 of  \cite{FHT}.
In any case for $\rho$ of type ${1\atop 2}$, take $\chi=\chi_\rho$, whereas
for type ${2\atop 1}$ take $\chi=\chi_\rho-\chi_{\rho'}$.

For example, consider $G=O(2)$ and $H=\bbT$, so $\epsilon$ is the determinant
$\delta$.
The irreducible ${1\atop 2}$-representations are precisely the $\kappa_i$, while
$1$ and $\delta$ are the ${2\atop 1}$ ones. Thus $^\epsilon R_G$ can be identified with
Span$\{a^i\ominus a^{-i}\}_{i\ge 1}$, and 
$^\epsilon R^1_G$ can be identified with $\bbZ 1^-$. In $^\epsilon R_G$,
$\delta$ acts like 1 and $\kappa_1$ takes $a\ominus a^{-1}$ to
$a^2\ominus a^{-2}$ and
$a^i\ominus a^{-i}$ (for $i>1$) to $(a^{i+1}\ominus a^{-i-1})
\oplus(a^{i-1}\ominus a^{-i+1})$. In $^\epsilon R_G^1$, $\delta$
acts like $-1$ and $\kappa_i$ like 0.

This picture is quite pretty when $G$ is a finite subgroup of $SU(2)$, in
which case $G$ is associated via the McKay correspondence to a graph of
(extended) ${A}$-${D}$-${E}$ type. A grading, i.e. a homomorphism $\varphi:G\rightarrow
\bbZ_2$, corresponds to an involution of the diagram; the node corresponding to
the trivial $G$-representation is sent to the node of a different $G$-character $\psi$.
The graph of the
kernel $G_0=\varphi^{-1}(0)$ is obtained by folding that of $G$ by that
involution, which identifies $G$-representation $\rho$ with $\rho\otimes\psi$.
A node fixed by the involution corresponds to a $G$-representation
of type ${1\atop 2}$; that node splits into the nodes of the $G_0$-representations
$\rho_1,\rho_2$. On the other hand, two nodes interchanged by the involution
correspond to the $G$-representations $\rho,\rho'$, and they collapse into
the $G_0$-node corresponding to $\overline{\rho}$.
Conversely, not all involutions correspond to gradings -- indeed, folding
by some involution of the $G_0$ graph fixing the trivial $G_0$-representation
will in some cases recover the $G$-graph.

In particular,
$\mathbb{A}_{2n-1}$ has a unique grading, given by rotation by $n$ in the
graph; the folded graph is $\mathbb{A}_{n-1}$. Similarly, $\mathbb{D}_{2n-1}$
has a unique grading, given by reflection through a horizontal mirror; the
folded graph will
be $\mathbb{A}_{4n-7}$. $\mathbb{D}_{2n}$ has two inequivalent gradings,
given by reflections through vertical or horizontal mirrors; the former folds
to $\mathbb{D}_{n+1}$ while the latter folds to $\mathbb{A}_{4n-5}$.
$\mathbb{E}_7$ has a unique grading, given by
reflection through a vertical mirror, and the folded graph  is $\mathbb{E}_6$.
The $\mathbb{D}\rightarrow\mathbb{A}$ and $\mathbb{E}_7
\rightarrow\mathbb{E}_6$ foldings are reversed by an appropriate folding.
The remaining groups, namely $\mathbb{A}_{2n}$, $\mathbb{E}_6$ and
$\mathbb{E}_8$, don't have a grading.

As we know, infinite inductions involve Dirac induction, which we've
already discussed in section 1.2. 
An independent example of Dirac induction is given in section 4.12 of
\cite{FHTlg}. The situation we need later is $G=SU(2)$ and $H=O(2)$. 
The coadjoint orbits of $G$ on $\mathfrak{g}^*$ are (see section 5.3 of \cite{Kir})
the fixed point 0 (stabiliser
$G$) and the sphere of radius $r>0$ (stabiliser $\bbT$). The obvious six-term
exact sequence identifies $R_{SU2}$ with coker$\,$Res$^{SU2}_{\bbT}$,
i.e. $\sigma_\ell\leftrightarrow [a^\ell]$ for $\ell\ge 1$.
Similarly, the coadjoint action of $O(2)$ identifies $^+R_{O2}$ with
$^-R^1_{O2}\oplus {\rm coker}\,^-$Res$^{O2}_{\bbT}$, i.e. $1\leftrightarrow
{1}^-+[1]$, $\delta\leftrightarrow -{1}^-+[1]$, and $\kappa_\ell
\leftrightarrow [a^\ell]$ for $\ell\ge 1$, whereas the graded ring
$^-R_{O2}$ is identified with coker$\,^+$Res$^{O2}_{\bbT}$, i.e.
$a^\ell\ominus a^{-\ell}\leftrightarrow[a^\ell]$ for $\ell\ge 1$;
finally, for $^-R_{O2}^1\,$, $\bbZ{1}^-$ is identified with ker$\,$Res$^{O2}_{\bbT}$,
i.e. ${1}^-\leftrightarrow 1-\delta$.

This means both Dirac restriction $^+$D-Res$^{SU2}_{O2}$ and Dirac induction
$^+$D-Ind$^{SU2}_{O2}$ interchange
 $\sigma_\ell\in R_{SU2}$ and $\kappa_\ell\in R_{O2}$, while
 $^+$D-Ind$^{SU2}_{O2}$ kills both 1 and $\delta$. Similarly, both
 $^-$D-Res$^{SU2}_{O2}$ and $^-$D-Ind$^{SU2}_{O2}$ interchange $\sigma_\ell$
 and $a^\ell\ominus a^{-\ell}$. Note that Dirac induction D-Ind$^{SU2}_{\bbT}$
 is the composition of $^+ {\rm Ind}^{O2}_{\bbT}$ with Dirac induction
 $^+$D-Ind$^{SU2}_{O2}$. 

The special case of Dirac induction between finite and Lie groups does not
seem to appear explicitly in the literature. For concreteness, consider the 
situation we will encounter later: $G=O(2)$ and $H$ a finite subgroup, eg.
a cyclic or binary dihedral group. Give $O(2)$ the grading $\epsilon$ coming 
from determinant. Then the Dirac induction from $R_H$ and $^\epsilon R_{O2}^1$
will send $\rho\in R_H$ to:
\begin{eqnarray}
\oplus_{\lambda\epsilon \mathrm{Irr}(O2)}&\![\lambda]\!&\bigl(\mathrm{dim\,Hom(
Res}^{O2}_HV_\lambda,\rho)-\mathrm{dim\,Hom(Res}^{O2}_HV_\lambda,\rho\otimes
\mathrm{Res}_H^{O2}\delta)\big)\\&=&1^-\,
(\mathrm{Mult}_1(\rho)-\mathrm{Mult}_d(\rho)\label{dindfinite}
\end{eqnarray}
where $d=\mathrm{Res}_H^{O2}(\delta)$ (see Theorem 2 of \cite{L2}).

\subsection{The geometry of adjoint actions}

In this subsection we explain how to construct the bundles we will need below.
Knowing the bundle is valuable in identifying some of the maps needed
in later sections. What we are after, for 
$^\tau K_G^*(X)$, is a bundle over $X$ with fibre  the compact operators on
a $G$-stable Hilbert space $\cH\,,$ i.e.  $\cH \simeq \cH \otimes L^2(G)$ as
$G$-spaces (though sometimes we can get away with $\cH=L^2(G)$ itself). We 
will focus on the most interesting case:
$G$ acting adjointly on itself. As explained in section 1.1, it suffices to
consider separately the $H^1_G$- and $H^3_G$-twists.

Consider first the group $G$ being $n$-torus $T=\bbR^n/L$,
where $L\subset\bbR^n$ is an $n$-dimensional lattice. Of course in this case,
the adjoint action will be trivial. By K\"unneth,
$H^1_T(T;\bbZ_2)\simeq H^1(T;\bbZ_2)\otimes\bbZ_2$ and $H^3_T(T;\bbZ)\simeq H^3(T;\bbZ)
\oplus H^2_T(pt;\bbZ)\otimes H^1(T;\bbZ)$. Consider first a trivial $H^1$-twist;
transgression implies (see section 7 of \cite{FHT}) we can ignore
$H^3(T;\bbZ)$; thus, introducing the dual lattice $L^*\simeq H^2_T(pt;\bbZ)\simeq
H^1(T;\bbZ)$, we obtain that the twist $\tau$ here (the `level') lies in 
${\rm Hom}(L,L^*)$. This level $\tau\in{\rm Hom}(L,L^*)$ can be written
in integer matrix form, once a basis $\{\beta_1,\ldots,\beta_r\}$ of $L$ is chosen,
by $k_{ij}=\tau(\beta_i)(\beta_j)\in\bbZ$. The level defines a map
$T\rightarrow T$ defined by $(t_1,\ldots,t_r)\mapsto(\sum_jt_jk_{1j},\ldots,
\sum_jt_jk_{rj})$. 

Consider now the easiest case: the 1-torus $\mathbb{T}$.
The $\mathbb{T}$-equivariant bundle $\mathcal{A}_k$ on $\mathbb{T}$ associated to level $k$ 
can be constructed as follows. Take Hilbert space $\cH=L^2(\mathbb{T})$ and let
$\mathcal{K}=\mathcal{K}(\cH)$ be the algebra of compacts. Let $U_k\in U(\cH)$
be the unitary operator corresponding to multiplication by the $\mathbb{T}$-character
$\chi_k$, so $U_k\pi U_k^*=\chi_k\pi$ where $\pi$ is the regular representation
of $\mathbb{T}$ (i.e. $U_k$ defines an equivalence $\pi\otimes \chi_k\simeq \pi$). 
Then $\mathcal{A}_k$ is the $\mathbb{T}$-bundle with fibres $\mathcal{K}$, 
whose sections $f$ are maps $f:[0,1]\rightarrow\mathcal{K}$ satisfying $f(0)=
U_kf(1)U^*_k$. Define an $\mathbb{T}$-action on $\mathcal{A}_k$ by
$(t.f)(s)=Ad(\pi(t))(f(s))=\pi(t)f(s)\pi(t)^{-1}$
(that this acts on $\mathcal{A}_k$, sending sections to sections, follows 
quickly from $U_k\pi U^*_k=\chi_k\pi$). If we were to ignore this $\mathbb{T}$-action,
then $\mathcal{A}_k$ would be trivialised by any continuous path from 1 to 
$U_k$, however as a $\mathbb{T}$-equivariant bundle it is nontrivial, for 
$k \ne 0$ (as can be seen by computing $K$-homology).
We call $U_k$ the {\it twisting unitary} for the bundle.

This bundle construction is easily generalised. Let $G$ be the torus
$T=\bbR^r/L$, for some $r$-dimensional lattice $L$ (we are most interested
in $T$ being a maximal torus of a compact Lie group, in which case $L$ is
the coroot lattice $Q$). Fix a level $k\in{\rm Hom}(L,L^*)$.
 The Hilbert space is $\cH=L^2(T)$; for any $\gamma\in L^*$ we
have a character $\chi_\gamma$ for $T=\bbR^r/L$ defined by $\chi_\gamma
(t)
=e^{2\pi \i \,\gamma(t)}$; define $U_\gamma$ as before by $U_\gamma \pi
U_\gamma^*
=\chi_\gamma \pi$. The bundle $\mathcal{A}_k$ on $T$, with
fibres the compacts $\mathcal{K}=\mathcal{K}(\cH)$, is defined using the gluing
conditions $f(t)=U_{k(\ell)} f(t+\ell) U_{k(\ell)}^*$, for all $\ell\in L$,
$t\in\bbR^r$.

An $H^1$-twist is possible for $T$ on $T$, arising from the target $T$ (as opposed to 
the group $T$), and the associated bundle is as follows. Return  for 
simplicity to $T=\mathbb{T}$ acting
adjointly on itself. Let $\widetilde{T}\simeq \mathbb{T}$ be a double-cover of $T$
(so the angle parametrising $\widetilde{T}$ is half that of $T$). Identifying
the space $L^2(T)$ with the completion of the space $\bbC[z^{\pm 1}]$
of polynomials, 
the space $\cH=L^2(\widetilde{T})$ becomes the completion of the polynomials $\bbC[z^{\pm \frac{1}{2}}]$
(half-integer powers are the spinors, and integer powers are the nonspinors). 
This nonspinor/spinor decomposition $\cH=\cH_0\oplus \cH_1$ provides a natural 
grading on the compacts
$\mathcal{K}(L^2(\widetilde{T}))$; act on the overlap of the cover of the circle
by the odd unitary $U=\left(\begin{matrix}0&z^{\frac{1}{2}}
\cr z^{-\frac{1}{2}}&0\end{matrix}\right)$ (interchanging those two subspaces) -- i.e.
as you wrap around the circle, the compact operator ${c}$ becomes $U{c}U^*$.

By contrast, consider the bundle for the orthogonal group $G=O(2)$ acting 
trivially on a point. Here the (trivial untwisted) `bundle' over that point 
consists of 
the compacts $\mathcal{K}(L^2(O(2)))$ with its obvious $O(2)$ action, and the
$H^3$-twist ($H^3_{O2}(pt;\bbZ)=\bbZ_2$) is obtained by replacing 
$L^2(O(2))$ by its spinors (the $O(2)$-spinors consist of half of the two-dimensional irreducible
representations of the double-cover $\widetilde{O}(2)\simeq O(2)$).
The untwisted bundle can be $H^1$-twisted (thanks to the group $O(2)$ being disconnected), 
essentially by doubling the point (which splits $O(2)$ into its two components,
each a copy of $\bbT$). More precisely, the graded space here will be 
$\cH=L^2(\bbT)\oplus L^2(\bbT)$, and the grading on $\mathcal{K}(\cH)$
is given by the odd unitary $U:(f,g)\mapsto (\overline{g},\overline{f})$.

Let $G$ be a compact semi-simple Lie group of rank $r$, eg. $G=SU(r+1)$. 
The orbits of $G$ acting adjointly 
on itself are of course the conjugacy classes of $G$.
A convenient way to parametrise these orbits 
uses the Stiefel diagram. Fix a maximal torus $T$ of $G$, which we can
identify with $\bbR^r/Q^\vee$ (where $Q^\vee$ is the coroot lattice).
The Stiefel diagram is an affine Weyl chamber: the affine Weyl group is
a semi-direct product of translations in the coroot lattice with the finite
Weyl group. More precisely, remove from the Cartan subalgebra $\bbR^r=\bbR\otimes_\bbZ
Q^\vee$
the hyperplanes fixed by a Weyl reflection $r_\alpha$, as well as the
translates of those hyperplanes by elements of the coroot lattice $Q^\vee$.
The Stiefel diagram $S$ is the closure of any connected component.  
Any orbit of the adjoint action intersects $S$ in one and only one point.
Points in the interior of $S$ correspond to generic (`regular') elements of $G$
and have stabiliser $T$, but points on the boundary will have larger stabiliser
(if $G$ is simply connected, the dimension of the boundary stabilisers
will be greater than that of the interior).

In the special case that $G$ is of A-D-E type, we can be more explicit.
 A natural basis for the Cartan subalgebra
${\mathfrak h}$ is provided
by the dual basis $\Lambda_i^\vee$ to the simple roots $\alpha_i\in{\mathfrak h}^*$.
The Killing form is an inner product on  ${\mathfrak h}$, and so allows us
to identify ${\mathfrak h}$ and its dual, and through this $\Lambda^\vee_i$
will be identified with the fundamental weights $\Lambda_i$. The Stiefel
diagram is the convex span of $\{0,\Lambda_1^\vee,\ldots,
\Lambda^\vee_r\}$, so any element $\xi$ in it can be written as a linear combination
$\xi=\sum_{i=1}^r \xi_i \Lambda^\vee_i$, where the Dynkin labels all satisfy 
$0\le \xi_i\le 1$.

For example, the Stiefel diagram for $G=SU(2)$ consists of a closed segment,
whose endpoints correspond to the fixed points $\pm I\in G$ with stabiliser 
$G$; we can identify those endpoints with the weights $0,\Lambda_1$ (see Figure 4).
The Stiefel diagram of the (nonsimply connected) group $G=SO(3)$ is an
interval with endpoints $I$ (stabiliser $G$ and weight 0) and ${\rm diag}(1,-I)$
(stabiliser $O(2)$ and weight $\Lambda_1$); generic points have stabiliser
$\bbT$ (see Figure 4).
The Stiefel diagram for $G=SU(3)$ is an equilateral triangle with vertices
we can identify with the weights $0,\Lambda_1$, and $\Lambda_2$
($\Lambda_i$ are the fundamental weights; exponentiated, these correspond to
the three scalar matrices in $SU(3)$);
the stabilisers at the vertices are $G$ and those on the edges are $U(2)$
(see Figure 5).

\bigskip 
\begin{figure}[tb]
\begin{center}
\epsfysize=.8in \centerline{\epsffile{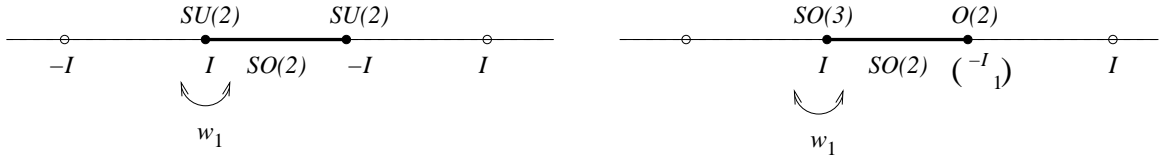}}
\caption{The Stiefel diagrams for $SU(2)$ and $SO(3)$ respectively}
\end{center}
\end{figure}

\begin{figure}[tb]
\begin{center}
\epsfysize=3in \centerline{\epsffile{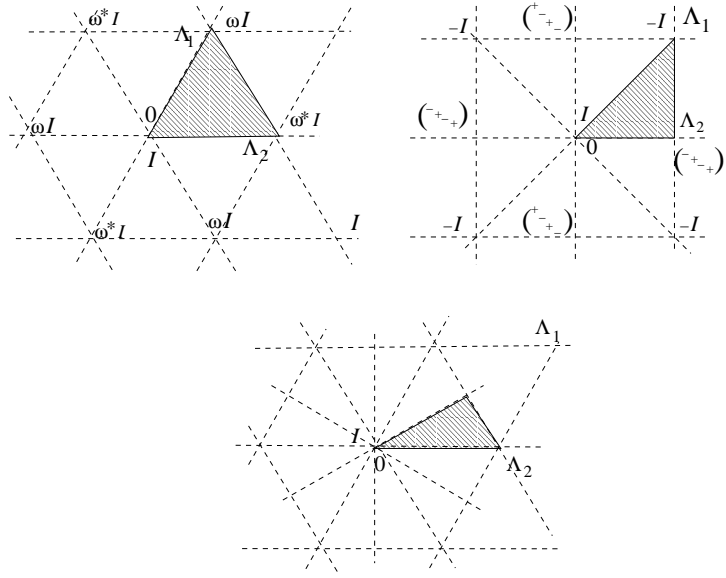}}
\caption{The Stiefel diagrams for $SU(3)$, $Sp(4)$, and
$G_2$, respectively}
\end{center}
\end{figure}

The Stiefel diagram of the symplectic group $G=Sp(4)\simeq Spin(5)$
 is a right isosceles triangle with
vertices $0,\Lambda_1$ and $\Lambda_2$ (in our labelling conventions, 
$\alpha_1$ is
the longest root, so $\Lambda_1$ is the longest fundamental weight, corresponding
to the five-dimensional representation $SO(5)$; these correspond to the
diagonal matrices $I_4$, $-I_4$, diag$(1,-1,1,-1)$
in $Sp(4)$) (see Figure 5(b)). We can take the maximal torus of $Sp(4)$ to be
the diagonal matrices diag$(\xi,\psi,\overline{\xi},\overline{\psi})$ for
complex numbers $\xi,\psi$ of modulus 1;
then the edges $0\leftrightarrow\Lambda_1$, $\Lambda_1\leftrightarrow\Lambda_2$,
$0\leftrightarrow\Lambda_2$, respectively, of the Stiefel diagram correspond to the
diagonal matrices diag$(\xi,\xi,\overline{\xi},\overline{\xi})$,
diag$(\xi,-1,\overline{\xi},-1)$ and diag$(1,\xi,1,\overline{\xi})$.
The stabilisers at the vertices are $G$, $G$ and
$SU(2)\times SU(2)$ respectively. The stabiliser at 
the edge $0 \leftrightarrow\Lambda_1$ is $U(2)$, at the edge 
$\Lambda_1\leftrightarrow\Lambda_2$ is $\bbT\times SU(2)$, and at the edge 
$0 \leftrightarrow\Lambda_2$ is $SU(2)\times \bbT$.

The level-$k$ bundle for $G=SU(2)$ can be constructed as for $T^r$, by
decomposing a representation of $G$ into weight-spaces (i.e. modules of the
maximal torus which are organised by the Weyl group). 
In particular, let $T$ be the maximal torus consisting 
of the diagonal matrices -- we can naturally identify it with the circle
$\bbR/Q$ where $Q=\sqrt{2}\bbZ$ is the (co)root lattice. The Hilbert space is 
$\cH=L^2(G)$. We want to associate a unitary $U_\gamma$ to any weight $\gamma
\in Q^*=1/\sqrt{2}\bbZ$. To do this, first fix a Stiefel diagram $S$ (here, 
half  of a fundamental domain for $T$).
 For any subrepresentation $\pi$ in $L^2(G)$, define `$\pi\otimes \gamma$'
as follows: restrict $\pi$ to $T$ (i.e. write its weight-space decomposition),
and in the Weyl-image $wS\subset T$ act like the character $\chi_{w\gamma}(e^{2\pi \i t})
=e^{2\pi \i \,w\gamma(t)}$. Then thanks to infinite-dimensionality,
$\cH\otimes \gamma\simeq\cH$ as both a representation of $T$ and the Weyl group,
so let $U_\gamma$ be the unitary defining that equivalence.
We can cover $G\simeq S^3$ with two patches: $D_1$ 
about the scalar matrix $I$ and $D_2$ about the scalar matrix $-I$.
The bundle $\mathcal{A}_k$ on $G$ (for $k$ the level), with
fibres the compacts $\mathcal{K}=\mathcal{K}(\cH)$, is defined by the 
following gluing
condition: identify $(gxg^{-1},c)$ in $D_1$  with $(gxg^{-1},Ad(\pi_g
U_k\pi_g^{-1})c)$ for any $gxg^{-1}\in D_1\cap D_2$, $c\in\mathcal{K}$. 
Again, $U_k$ is called the twisting unitary.

The consistency condition for these bundles is that when $gxg^{-1}=x$,
then $Ad(\pi_gU_k\pi_g^{-1})$ should be the identity, i.e. $\pi_gU_k\pi_g^{-1}
=\lambda_g I$ for some character (i.e. one-dimensional representation) $g\mapsto
\lambda_g$ of the stabiliser $C_G(x)$. 

When $G=SO(3)$ we have $H^1_G(G;\bbZ_2)\simeq \bbZ_2$ and $H^3_G(G;\bbZ)\simeq
\bbZ\oplus \bbZ_2$. The representation ring $R_G$ is the polynomial ring
$\bbZ[\sigma_3]$, while the spinors (corresponding to the torsion part of
$H^3_G$) form the $R_G$-module $\sigma_2\bbZ[\sigma_3]$. The nontorsion
part $\tau_{3,non}\in\bbZ$ of $H^3_G$ is done as for $SU(2)$; the torsion part 
$\tau_{3,tor}$ is done by decomposing 
the $SO(3)$-module $L^2(SU(2))$ into nonspinors (the space for $\tau_{3,tor}=0$)
and spinors (for $\tau_{3,tor}=1$). The $H^1_G$-twist is handled analogously to
that of $\bbT$, by putting a grading on the $SO(3)$-module $L^2(SU(2))=
\cH_{ns}\oplus\cH_{sp}$ and using the odd automorphism $U=\left(\begin{matrix}
0&a^{\frac{1}{2}}\cr a^{-\frac{1}{2}}&0\end{matrix}\right)$ 
(where $a^{\frac{1}{2}}
\in L^2(\widetilde{\bbT})$ at a generic point) on the overlap.

The level-$k$ bundle for $G=SU(3)$ is similar to that of $SU(2)$. Cover the Stiefel diagram
with a patch $D_0,D_1,D_2$ about each vertex. To the overlap between patch
$i$ and patch $j$, assign the twisting unitary $U_{ij}:=U_{k(\Lambda_j-\Lambda_i)}$
(where $\Lambda_0:=0$). We must check the consistency condition -- it suffices
to consider the boundary of the Stiefel diagram, say the edge diag$(z,z,z^{-2})
\subset T^2$. Which character of the stabiliser
$\{{\rm diag}(U),det(U)^{-1})\}\simeq U(2)$ restricts to the character $k\Lambda_1$ of the torus
$T^2$ of $SU(3)$? On that edge the twisting unitary acts
like $determinant^k$ (the only one-dimensional representations of the stabiliser
$U(2)$). 

The level-$k$ bundle for $G=Sp(4)$ is the same; the unitary on the
edges again corresponds to $det^k$. From these
examples it should be clear how to obtain any other bundle for $G$ acting 
adjointly on itself -- the $G_2$-bundle is explicitly described at the end of 
section 2.3.
Torsion in $H^3_G$ corresponds to groups $G$ which are nonsimply connected, 
as we explained with $O(2)$. An $H^1_G$-twist is
obtained by using a double-cover of $G$ to get a grading. 

Note that these considerations imply $H^3_G(G;\bbZ)$ contains $\bbZ$; of course
the former can be calculated by eg. spectral sequences and is found to equal
that $\bbZ$.

For $G$ compact semi-simple (say of rank $r$), Meinrenken \cite{M} found
an elegant construction of the $G$ on $G$ bundle at level 1, using the
basic representation of the associated affine algebra. Think of the Stiefel
diagram as an $r$-dimensional simplex with $r+1$ vertices, labelled say from
0 to $r$. Every nonempty
subset $I$ of $\{0,1,\ldots,r\}$ are the vertices of a subsimplex, parametrising
points in $G$ containing some stabiliser $G_I$ (though the boundary points
of this subsimplex will have larger stabiliser). The Lie algebra of $G_I$ 
is naturally identified with the Lie algebra obtained from the affine algebra
of $G$, obtained by deleting the vertices $I$ from the affine Dynkin diagram.
More precisely, we get a natural embedding of the (finite-dimensional)
Lie algebras of these stabilisers, into the (infinite-dimensional) loop
algebra. The level 1 basic representation of that affine algebra then restricts
to a coherent family of projective representations of those stabilisers,
and from this the bundle is formed (see equation (21) in \cite{M}).
By using the affine algebra representation, he obtains almost for free 
a global description of the bundle, avoiding our complicated explicit
construction of unitaries and verification of their consistency conditions.
On the other hand our construction is more general, permitting eg. $H^1_G$-twists,
is more explicit (which can help in identifying some of the maps in six-term
and Mayer-Vietoris sequences), and is inherently finite-dimensional.

In both the physics literature \cite{FFFS,MMS} and mathematics
literature (this is done very explicitly in \cite{M}), primaries are
identified with certain conjugacy classes in $G$. For example, when $G$ is of 
A-D-E type, the level $k$ primaries are naturally identified with
the conjugacy classes corresponding to points in the Stiefel diagram
with Dynkin labels $\xi_i\in(1/(k+h^\vee))\bbZ$. From our standpoint, a natural
task would be to associate to each of these conjugacy classes, a section 
of the $G$ on $G$ bundle.

This is fairly straightforward to do for $G=\bbT$ -- see the end of section
3.1 for the details for twisted $K$-theory.

\subsection{The level calculations}\bigskip

Using the  bundles constructed in subsection 2.2, we can compute the
level $k$ of the conformal embeddings $H_k\rightarrow G_1$. This provides a
nontrivial consistency check.
In this subsection we work out several examples.

Consider first the $\bbT_2\rightarrow SU(2)_1$ conformal embedding.
 The level `2' arises as the inner product 
$\alpha\cdot \alpha$: as we wrap around $T$, we traverse the Stiefel diagram
of $SU(2)$ twice, and so pass through the  overlap of the bundle cover, twice.
The first time picks up the unitary $U_{\Lambda}$, and the second picks up
the unitary $U^*_{r.\Lambda}$ where we Weyl reflect the fundamental
weight $\Lambda$ (we invert $U$ because the overlap is traversed in the 
opposite direction). The resulting unitary corresponds to weight $\Lambda-
r(\Lambda)=\alpha$, and hence to $T$-character $\alpha\cdot\alpha=2$.

The level $k\in {\rm Hom}(A_2,A_2^*)$ for $T^2_k\rightarrow SU(3)_1$
($A_2$ is the hexagonal lattice)
is recovered very similarly, and this shows how this works in general for
conformal embeddings of the maximal torus. For convenience make the patch
about the $\Lambda_2$-vertex of the $SU(3)$ Stiefel diagram very small, as
in Figure 5. This $A_2$ (co)root lattice is the span of the simple roots
$\alpha_1,\alpha_2$. Move first along the $\alpha_1$-direction: we cross
from the 0-patch to the $\Lambda_1$-patch, given by unitary $U_{r_{\alpha_2}(
\Lambda_1)}=U_{\Lambda_1}$, and back again, given by $U_{r_{\alpha_1}r_{\alpha_2}
(\Lambda_1)}^*=U_{\Lambda_1-\alpha_1}^*$. As with the $SU(2)$ calculation,
$k(\alpha_1)\in A_2^*$ will be the net weight picked up, and will equal the
difference $\Lambda_1-(\Lambda_1-\alpha_1)=\alpha_1$. Similarly, $k(\alpha_2)
=\alpha_2$, so the $LT^2$-level $k$ is given by the identity map.

More interesting is to recover from the bundles the level for the conformal embedding of
$SU(2)$ (or more precisely $SO(3)$) into $SU(3)$. This map $\mathcal{R}^{(3)}$
is given explicitly
in (\ref{SU2SU3}) below, from which we read off that the $SU(2)$ Stiefel
diagram $I_2\leftrightarrow -I_2$ embeds into the line $(2t,0,-2t)$ in the $SU(3)$ Cartan subalgebra
$(x,y,-x-y)$ as in Figure 6: the endpoints correspond to $t=0$ and 
 $t=\frac{1}{2}$ (the $-I_2$-endpoint should lie at the first coroot lattice
point on the segment after 0, since $-I_2$ lies in the kernel of 
$\mathcal{R}^{(3)}$). The $SU(2)$ simple root, in this $SU(3)$ notation,
 corresponds to $t=1$. As we move along this $SU(2)$ Stiefel diagram, we
 see we twice have to change patches in the $SU(3)$ bundle. As always,
 this is where the twisting unitaries arise: the net twist here is $\Lambda_1-
 r_{\alpha_1+\alpha_2}(\Lambda_1)=\alpha_1+\alpha_2=(1,0,-1)$. This will correspond to an
 $SU(2)$ twist of $k\Lambda_1^{su2}$ where the level $k$ is obtained by
 the inner product of the net twist $(1,0,-1)$ with the $SU(2)$ simple
 root $(2,0,-2)$. In this way we recover the value $k=4$.

\bigskip 
\begin{figure}[tb]
\begin{center}
\epsfysize=3in \centerline{\epsffile{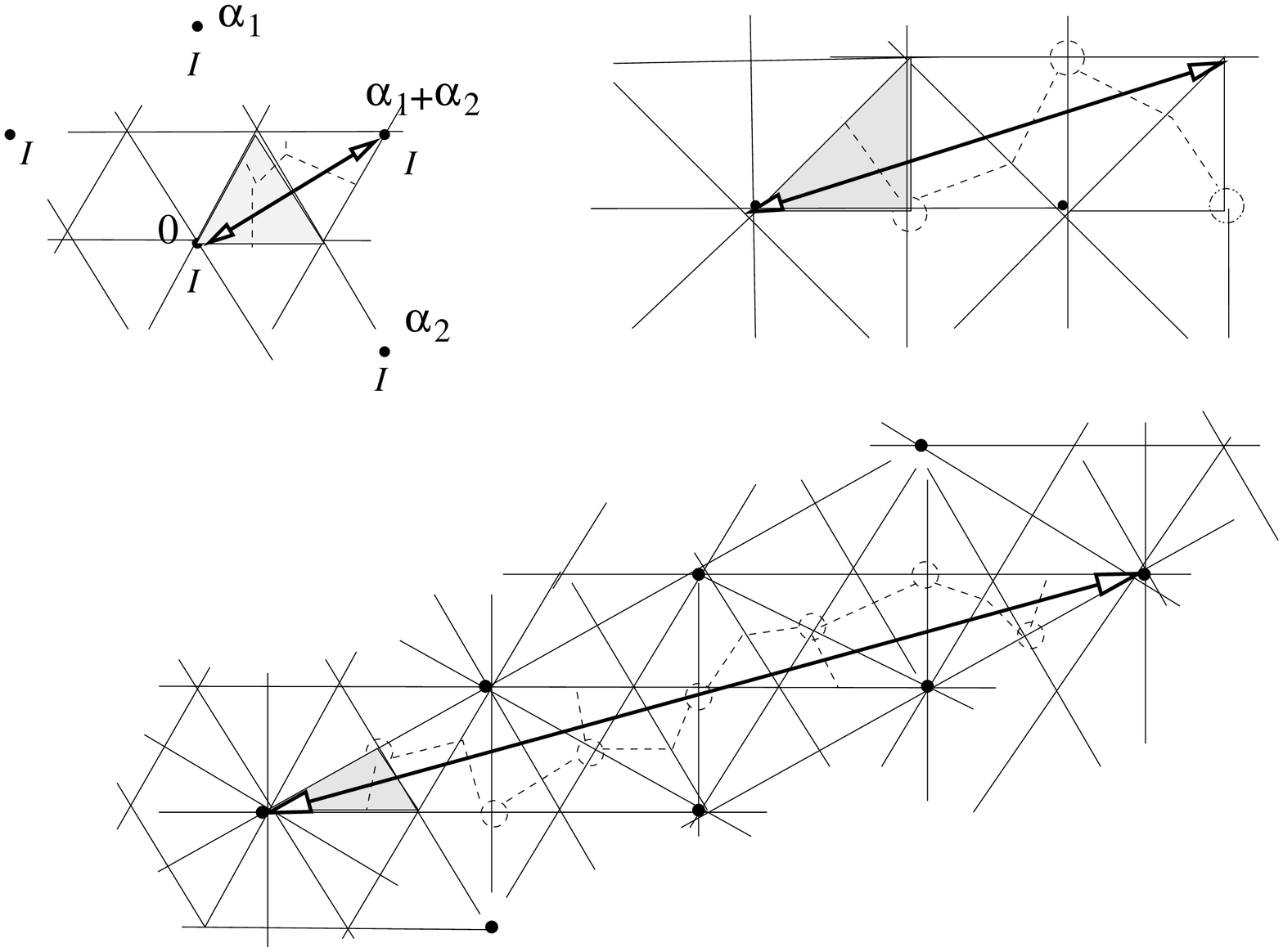}}
\caption{$SU(2)$ Stiefel diagram in $SU(3)$, $Sp(4)$, $G_2$ 
Cartan subalgebras}
\end{center}
\end{figure}

The conformal embedding of 
$SU(2)$ into $Sp(4)$, given explicitly by $\mathcal{R}^{(4)}$ in 
(\ref{su2sp4}) below, behaves similarly. The $SU(2)$ Stiefel
diagram embeds into the line $\frac{1}{\sqrt{2}}(3t,t)$ in the $Sp(4)$ Cartan 
subalgebra $(x,y)$ as in Figure 6, with endpoints at $t=0$ and 
 $t=1$; $t=2$ is the $SU(2)$ simple root (here, the simple root should 
correspond to the first coroot lattice point on the segment after 0, since
$-I_2$ does not lie in the kernel of $\mathcal{R}^{(4)}$).
As we move along this $SU(2)$ Stiefel diagram, we
 change patches three times, for a net twist of $\frac{1}{\sqrt{2}}(1,1)-
 \frac{1}{\sqrt{2}}(-1,1)+\frac{1}{\sqrt{2}}(1,1)=\frac{1}{\sqrt{2}}(3,1)$.
The level $k$ is thus $\frac{1}{\sqrt{2}}(3,1)\cdot\frac{1}{\sqrt{2}}(6,2)$.  
In this way we recover the value $k=10$.

The conformal embedding of $SU(2)$ (more precisely, $SO(3)$) into the 
compact Lie group of type $G_2$, can be analysed similarly, except we
don't give the explicit form for it (it will take the form of the
7-dimensional irreducible two-to-one representation of $SU(2)$, embedded into the
7-dimensional irreducible representation of $G_2$, whose image can be
identified with $G_2$). Choosing the realisation roots $\alpha_1=
(-1,2,-1)/\sqrt{3}$, $\alpha_2=(1,-1,0)/\sqrt{3}$, the 
Stiefel diagram can be taken to be the triangle with vertices at 0, $\Lambda_2=
(1,0,-1)/\sqrt{3}$ and $\Lambda_1/2=
(1,1,-2)/2\sqrt{3}$ (with the twisting unitaries in each patch given by the 
weights $2(\Lambda_2-0)$, $2(\Lambda_1/2-0)$, $2(\Lambda_1/2-\Lambda_2)$ -- 
the doubling is needed to get $G_2$ weights). Because we don't have the 
explicit mapping, we need help to
see how the $SU(2)$ Stiefel diagram fits inside the $G_2$ Cartan subalgebra:
Table 13 of \cite{Dyn} tells us the $SU(2)$ Cartan subalgebra is the
line $t(10\alpha_1+18\alpha_2)=\frac{t}{\sqrt{3}}(8,2,-10)$ (see Figure 6).
The endpoints of the $SU(2)$ Stiefel diagram are thus at $t=0$ and
$t=\frac{1}{2}$, and the $SU(2)$ simple root is at $t=1$.
We get 8 patch crossings, for a total twist of $2\times\bigl(
(2,0,-2)-(0,2,-2)+(-2,0,2)+(-1,0,1)\bigr)/\sqrt{3}=(-2,-4,6)/\sqrt{3}$.
Thus the $SU(2)$ level is $k=(-2,-4,6)\cdot (8,2,-10)/3=-28$.

\section{Conformal embeddings: the first examples} 

Given a conformal embedding $H_k\rightarrow G_\ell$ and the choice of diagonal
$G_\ell$ modular invariant ${Z}=\sum_\mu |\chi_\mu|^2$,
it is natural to guess that $^\tau K_*^H(G)$ recovers the full system of
the corresponding $H_k$ modular invariant -- see the end of section 1.4
for the finite group analog of this statement, which works perfectly. As we 
will find, this $K$-homological interpretation of conformal embeddings of
loop groups isn't as clean as one would like. We will give in
this section the easiest nontrivial examples, and in sections 5 and 6
give more serious examples.  

In order for this approach to conformal embeddings to work, we should have
that $H^3_H(G;\bbZ)$ contains a copy of $\bbZ$ which can be identified with
$H^3_G(G;\bbZ)\,$. But an element of $H^3_G(G;\bbZ)$ corresponds to a
$G$-equivariant bundle $K_\tau$
of compact operators on $G$, as explained in section 1.1, so restricting 
equivariance to $H$ defines the
appropriate element of $H_H^3(G;\bbZ)\,$. This can also be seen from the Borel
construction of group cohomology $H^*_G(X)=H^*((X\times EG)/G)\,$: we can
identify the universal coverings $EH$ with $EG$, so the natural projection 
$(X\times EG)/H
\rightarrow (X\times EG)/G$ becomes the map $H^*_G(X)\rightarrow
H^*_H(X)\,$.

\subsection{The Verlinde algebras for the circle}\bigskip

We consider here the $K$-homology calculations of the level $k$ Verlinde 
algebra of the circle $\mathbb{T}$, where this is understood as in section 1.2.
This falls under the Freed-Hopkins-Teleman umbrella and constitutes the
easiest example.

In fact the calculation is given in section 4 of \cite{FHTii}. 
Let $L\subset\bbR^n$ be an $n$-dimensional lattice and $L^*={\rm Hom}(L,\bbZ)$
be the dual lattice, and consider the $n$-torus $T=\bbR^n/L$.
As mentioned earlier, the twist $\tau$ here (the `level') lies in 
${\rm Hom}(L,L^*)$. 
They obtain  $^\tau K_0^T(T)\simeq
\bbZ\langle L^*/\tau(L)\rangle$, the Verlinde algebra, and $^\tau K_1^T(T)=0$.

In order to motivate the calculations given next section, it is helpful to
redo this calculation explicitly for the maximal torus $\bbR/\sqrt{2}\bbZ$
of $SU(2)$.
The orbit analysis is trivial: we have $\mathbb{T}$ acting on itself by conjugation,
but because it's abelian this action is trivial. So each point of $\bbT$ is
itself an orbit, with full stabiliser $\mathbb{T}$. Here, $L=\sqrt{2}\bbZ$
and $L^*=\frac{1}{\sqrt{2}}\bbZ$, and the twist $\tau$ can be identified 
with a nonzero even integer $k$ (so $L^*/\tau(L)\simeq \bbZ/k\bbZ$).

The $K_*$-groups  are most simply computed by Mayer-Vietoris (\ref{MV}):
\begin{equation}
   \begin{matrix}{}^{\tau}K^{\mathbb{T}}_0(\mathbb{T})&\longrightarrow&K^{
\mathbb{T}}_0(\bbR)\times 2&\longrightarrow&K^{\mathbb{T}}_0(\bbR)\times 2
   \cr\uparrow&&&&\downarrow\cr K^{\mathbb{T}}_1(\bbR)\times 2
   &{\beta\atop\longleftarrow}&K^{\mathbb{T}}_1(\bbR)\times 2&
   \longleftarrow&{}^{\tau}K^{\mathbb{T}}_1(\mathbb{T})\end{matrix}
\end{equation}
$K_1^{\mathbb{T}}(\bbR)$ is the representation ring
$R_{\mathbb{T}}=\bbZ[a^{\pm 1}]$; we've dropped the twist on those $K$-homology groups
because $H^1_{\mathbb{T}}(\bbR;\bbZ_2)$ and $H^3_{\mathbb{T}}(\bbR;\bbZ)$ both vanish.
The map $\beta:R_{\mathbb{T}}{}^2\rightarrow R_{\mathbb{T}}{}^2$ presumably sends
$(p(a),q(a))$ to $(p(a)+q(a),p(a)\pm a^k q(a))$, where `$k$'
is the $L\bbT$ level. The effect of the $H^1$-twist
would be to introduce the sign (`+' should correspond to ungraded, in order to
recover the nonequivariant $K$-homology $K_*(\mathbb{T})=\bbZ$). Then
${}^{\tau}K^{\mathbb{T}}_0(\mathbb{T})={\rm coker}\,\beta=\bbZ[a^{\pm 1}]/(1\mp a^k)$
(where we take the lower sign, i.e. `+', if we $H^1$-twist) and
${}^{\tau}K^{\mathbb{T}}_1(\mathbb{T})={\rm ker}\,\beta=0$. We can also handle the
$H^1$-twist through (\ref{lineb}), writing $^\tau K_*^{\mathbb{T}}(\mathbb{T})\simeq {}^{\tau'}
K_{*+1}^{\mathbb{T}}(M\ddot{o}b)$ where $M\ddot{o}b$ is the open M\"obius 
strip, and $\tau'\in H^3_G(G;\bbZ)$ is the nontorsion part of $\tau$.

Now,  $\bbZ[a^{\pm 1}]/(1-a^k)$ corresponds to
a $k$-dimensional ring with cyclic fusion product generated by $a$.
For $k$ odd, $\bbZ[a^{\pm}]/(1+a^k)$
is also cyclic, with generator $-a$, but for $k$ even 
 that $k$-dimensional ring is not cyclic.
This seems to suggest that we should not $H^1$-twist here. 

In summary, for this $L\bbT$ example, $K_1$
vanishes and $K_0$ gives the Verlinde algebra. We should not $H^1$-twist 
these $K_*$-groups.

Using the bundle $\mathcal{A}_k$ constructed in subsection 1.2, 
the six-term exact sequence
(\ref{six-term}) (removing a point from $\mathbb{T}$) becomes
\begin{equation}
\begin{array}{ccccc}
    0 & {\longleftarrow} & K^{0}(C_0(X;\mathcal{A}_k)\rtimes G) &
    {\longleftarrow} & \bbZ[a^{\pm 1}]
    \\[3pt]
    \downarrow  & & & & \uparrow\beta
    \\[2pt]
    0 & {\lori} &
    K^1(C_0(X;\mathcal{A}_k)\rtimes G) &{\lori} &\bbZ[a^{\pm 1}]
    \end{array}
    \end{equation}
where $\beta$ corresponds to multiplying by $1\mp a^k$, depending on the
sign of the $H^1$-twist. This recovers the previous result.

We can give a more precise description of the dual $K^*$-groups which are again
most simply computed by Mayer-Vietoris (\ref{MV}):
\begin{equation}
   \begin{matrix}{}^{\tau}K_{\mathbb{T}}^0(\mathbb{T})&\longleftarrow&K_{
\mathbb{T}}^0(\bbR)\times 2&\longleftarrow&K_{\mathbb{T}}^0(\bbR)\times 2
   \cr\downarrow&&&&\uparrow\cr K_{\mathbb{T}}^1(\bbR)\times 2
   &{\gamma\atop\longrightarrow}&K_{\mathbb{T}}^1(\bbR)\times 2&
   \longrightarrow&{}^{\tau}K_{\mathbb{T}}^1(\mathbb{T})\end{matrix}
\end{equation}

To keep track, we call the open sets $U_1$ and $U_2$ both homeomorphic to $\bbR$ in $\mathbb{T}$.
Regard this as : $\gamma$ is surjective on the first co-ordinate, so $K^1_\mathbb{T}(U_1)$ does not
contribute to $K^1_{\mathbb{T}}(\mathbb{T})$ under the map $\gamma$ -- indeed it is killed by
$\gamma$. So $K^1_\mathbb{T}(\mathbb{T})$ is described by $K^1_\mathbb{T}(U_2)$  -- under $\gamma$.
However, because of exactness there are relations imposed on $K^1_T(U_2)$ when mapped into 
$K^1_\mathbb{T}(\mathbb{T})$, namely 
$$K^1_\mathbb{T}(\mathbb{T}) \simeq K^1_\mathbb{T}(U_1)/(1-\alpha^k)K^1_\mathbb{T}(U_1) \simeq
R_{\mathbb{T}}/(1-\alpha^k)R_{\mathbb{T}}\,.$$
An alternative viewpoint is via the six-term exact sequence for $K$-theory,
for the open set $U = \bbT$  with a point removed. 
This gives $K^1_{\mathbb{T}}(\mathbb{T}) \simeq K^1_\mathbb{T}(U)/\exp( K^1_\mathbb{T}(pt))\,$,
using the exponential map.

This yields the following description of $K^1_\mathbb{T}(\mathbb{T})$ in terms of unitaries in the
(unitalisation) of the 
twisted bundle $\cA_k$ of section 2.2 of compacts on the circle. 
The sections of the bundle are maps $f$ from $[0,1]$ into $\cK(L^2(\mathbb{T}))$
 such that $f(0) = Ad(U_k)(f(1))$, where $U_k$ is the unitary associated with the twist $k$. Take
the 
 equivariant Bott map $\beta: K^0_{\mathbb{T}}(pt)  \rightarrow K^1_\mathbb{T}(U)$\,,
 where $1$ is the unit of $\mathbb(T)$
 Then $\beta(1) = w_0 =  ze_0 + 1-e_0$, where $z$ is the natural loop in $\mathbb{T}$, and
 $e_j$ is the projection in  $L^2(\mathbb{T})$ corresponding to the character $j$. 
 The action of $\mathbb{Z}_k$, identified with the
 $\{\omega^j : j\}$ where $\omega = \exp(2 \pi i /k)$
on $\mathbb{T}$ by rotation induces an action on the
 bundle $\cA_k$ and hence on $K^1_\mathbb{T}(\bbT)$. 
 (This is best seen by breaking up the bundle into an equivalent one where
 we have $k$ cuts on the circle with a jump indicated by $U_1$ at each, so that
  sections of the bundle
 are maps 
 $f$ from $[0,1]$ into $\cK(L^2(\mathbb{T}))$
 such that $f(j/k) = Ad(U_1)(f((j+1)/k))\,, i = 0,1, \dots k-1$)
 This rotation takes the generator $[w_0]$ to
 $[w_j] = [ ze_j + 1-e_j]$, compatible with the  equivariant Bott maps 
 $\beta: K^0_{\mathbb{T}}(\omega^j)  \rightarrow K^1_\mathbb{T}(U)$\,.
 Note that $[w_n] = [w_{n+k}]$ in $K^1$ due to the nature of the bundle $\cA_k$.
The labelling of the $k$ primary fields has a dual meaning in terms of the representations
of $\mathbb{T}$ or of the (conjugacy classes) of the points $\{\omega^j : j\}$
on the circle.

\subsection{The $\bbT_2\rightarrow SU(2)_1$ conformal
embedding}\bigskip

We consider next the $K$-homology calculations of the conformal embedding
$\bbT_2\rightarrow SU(2)_1$ (corresponding to `$A_2$' on
the A-D-E list of modular invariants in \cite{CIZ}). The 
$\bbT$ level is most
easily obtained  by comparing characters: the two irreducible level 1 characters of
the loop group $LSU(2)$ are theta functions divided by $\eta(\tau)$, and coincide 
with the two $L\bbT_2$ characters (so the branching rules here are
trivial). This conformal embedding, together with the diagonal $SU(2)_1$
modular invariant, yields the
diagonal $\bbT_2$ modular invariant ${Z}=I$.
The resulting full system should thus be two-dimensional.
In fact it should be identifiable with the cyclic Verlinde algebra 
$\bbZ[a^{\pm 1}]/(1-a^2)$.

The orbit analysis is easy: we have the maximal torus $\bbT$ acting on
$SU(2)$  by conjugation, where we identify $\mathbb{T}$ with the 
 diagonal matrices $\left(\begin{matrix}\alpha&0\cr 0&\overline{\alpha}
\end{matrix}\right)$.
The orbits are $\cO_f=\mathbb{T}$ (the diagonal matrices), with the full $\bbT$ as
stabiliser, and the generic points, with $C_2=\pm I$ as stabiliser
(corresponding to the centre of $SU(2)$).
The generic orbits together form $\cO_{g}=\bbR^2\times \bbT/C_2$. To see this,
parametrise $SU(2)$  with matrices $\left(\begin{matrix}\beta&\gamma\cr
-\overline{\gamma}&\overline{\beta}\end{matrix}\right)$ where $|\beta|^2+|\gamma|^2=1$. 
The generic orbits correspond
to $\gamma\ne 0$; the resulting $\bbT$ orbit will contain exactly one matrix
whose $\gamma$ entry is a positive real number. Hence each generic orbit is
 uniquely determined by its value of $\beta$, which will lie in 
the interior of the unit disc, and this is the $\bbR^2$.

Next, we need the cohomology groups $H^1_{\bbT}(\mathbb{T};\bbZ_2)=\bbZ_2$, $H^3_{\bbT}(\mathbb{T};
\bbZ)=\bbZ$, $H^1_{\bbT}(\bbR^2\times \bbT/C_2;\bbZ_2)={\rm Hom}(C_2,\bbZ_2)=\bbZ_2$
and  $H^3_{\bbT}(\bbR^2\times \bbT/C_2;\bbZ)=0$. $H^*_{\bbT}$ for $\cO_{f}$ 
was easiest to compute using K\"unneth. Also, spectral sequences immaediately 
tell us
$H^1_{\bbT}(SU(2);\bbZ_2)=0$ and $H^3_{\bbT}(SU(2);\bbZ)$ is either $\bbZ$ or
0 (hence
we know $H^3_{\bbT}(SU(2);\bbZ)=\bbZ$, since it must see the level $k$).

The obvious six-term exact sequence reads
\begin{equation}
   \begin{matrix}K^{C_2}_1(pt)&\longleftarrow&{}^\tau K^{\bbT
}_0(SU(2))&\longleftarrow
   &{}^{\tau''}K^{\bbT}_0(\mathbb{T})
   \cr\downarrow&&&&\uparrow\beta\cr ^{\tau''}K^{\bbT}_1(\mathbb{T})
   &\longrightarrow&{}^\tau K^{\bbT}_1(SU(2))&
   \longrightarrow&K^{C_2}_0(pt)\end{matrix}
\end{equation}
The $K_*$-groups for orbits $\cO_{f}$ are computed last subsection, and we obtain
${}^{\tau'}\!K^{\bbT}_0(\mathbb{T})=\bbZ[a^{\pm 1}]/(1-a^{2(k+2)})\simeq\bbZ^{2(k+2)}$
and ${}^{\tau'}K^{\bbT}_1(\mathbb{T})=0$ (since there is no $H^1$-twist here). 
Here, $a$ is the generator of the representation ring $R_{\bbT}$, and
$k+2$ is the $LSU(2)$ level shifted as usual by its dual Coxeter
number ($k=1$ corresponds to the conformal embedding). The factor there of
2 is explained in subsection 2.3.
Also, $K_0^{C_2}(pt)=R_{C_2}\simeq\bbZ^2$ and
$K_1^{C_2}(pt)=0$, and $\beta$ is an injection. We obtain $K_0^{\bbT}(SU(2))
\simeq \bbZ^{2k+2}$ and $K_1^{\bbT}(SU(2))=0$. (We compute this more elegantly
in the following subsection. This example was also computed in \cite{SN1}.)

Again, $k=1$ corresponds to the conformal embedding, but its full system
is only two-dimensional, not four. Next subsection we find that a similar
phenomenon occurs with many other conformal embeddings. There
we identify the multiplicity two occurring here with the order of the
Weyl group $C_2$ of
$SU(2)$, or if you prefer with the Euler number of the sphere $SU(2)/\bbT$. We 
discuss what this could mean in the concluding section.

\subsection{The Hodgkin spectral sequence}\bigskip

The Hodgkin spectral sequence (Thm. 6.1 of \cite{RS}) is a powerful tool for
calculating many $K$-groups $^\tau K_H^*(G)$. In particular, suppose
$G$ is a compact connected Lie group, with torsion-free fundamental group
(eg. a torus $T^n$ or a simply connected group), and $H$ is a closed
subgroup of $G$. As in section 1.1, let $X$ be a space on which $G$ acts,
so $B=C_0(X;\mathcal{K}_\tau)$ is a $C^*$-algebra carrying a $G$-action.
Then there is a spectral sequence of $R_G$-modules which strongly converges
to $K_*^H(B)={}^\tau K_H^*(X)$, with 
\begin{equation}
E^2_{p,q}={\rm Tor}_p^{R_{{G}}}(R_H,{}^\tau K^{q}_G(X)) \,.
\end{equation}
We are most interested in $X=G$, with the adjoint action of $G$, in which
case ${}^\tau K^q_G(X)=Ver_k(G)$ or 0, for a level $k$ determined from $\tau$.

Consider first the situation where $H$ is maximal (i.e. of full rank)
 in $G$. There are many
examples of this, eg. $SU(n)_1\rightarrow SU(p)_1\times SU(n-p)_1\times U(1)_{np(n-p)}$,
$G_{2,1}\rightarrow SU(2)_{3}\times SU(2)_{1}\,$, $E_{8,1}\rightarrow SU(9)_{1}$ are some
among infinitely many (\cite{SN2} used essentially this method to compute a
subset of these, namely those corresponding to Hermitian symmetric spaces).
When $H$ is of maximal rank,
\cite{Pit} (together with the validity of Serre's conjecture
that projective modules over polynomial rings over fields or PIDs are free)
tells us that $R_H$ is free over $R_G$, say $R_H\simeq (R_G)^d$. This means
that $E^2_{0,q}=Ver_1(G)\otimes_{R_G} R_H$ for $q$ odd, and all
other $E^2_{p,q}$ vanish. Hence $^\tau K^0_H(G)=0$ and
\begin{equation}
^\tau K^1_H(G)=Ver_1(G)\otimes_{R_G} R_H=R_H/{I}_1=(Ver_1(G))^d\, ,
\end{equation}
where ${I}_1$ is the level-1 fusion ideal 
given by $Ver_1(G)\simeq R_G/\mathcal{I}_1\,$. This rank $d$ 
is given by $d=|W_G|/|W_H|\,$, the Euler number of $G/H$, where $W_G,W_H$ are the
Weyl groups of $G,H$ respectively (see equation (2.9) of \cite{SN2} for the
easy derivation).

One of the easiest of these examples is $E_{8,1}\rightarrow SU(9)_1\,$: there
is only one $E_{8,1}$ modular invariant, namely the diagonal one, and it
restricts to the $SU(9)_1$ modular invariant ${Z}=|\chi_{00000000}+\chi_{00100000}
+\chi_{00000100}|^2$. The full system of ${Z}$ will be 9-dimensional.
On the other hand $Ver_1(E_8)$ is 1-dimensional, and $|W_{E8}|/|W_{SU9}|=
1920$, so $^\tau K^1_{SU9}(E_8)$ is 1920-dimensional. There certainly is
room in $^\tau K^1_{SU9}(E_8)$ for the full system, but the meaning of
the rest is unclear to us, and again is discussed in section 7.

For most conformal embeddings, $H$ does not have full rank, but at least when
$G$ has small rank, then this spectral sequence can still be very useful.
We will see this in section 6.1 below, where we compute $^\tau K_H^*(G)$ for
 $SU(2)_{10}\rightarrow Sp(4)_1$. Unfortunately, for  the 
embedding $SU(2)_4 \rightarrow SU(3)_1$ considered in section 5 (similarly
$SU(2)_{28}\rightarrow G_{2,1}$),
it is the homomorphic image $SO(3)$ and not $SU(2)$ which is embedded in
$SU(3)$, while we are interested in the $K$-groups
with respect to $SU(2)$. For those examples, the Hodgkin spectral sequence
does not seem to have a
direct use and we must dive into the orbit analysis.

\section{Permutation orbifolds} 

Let $G$ be connected compact and simply connected (although we also
take $G=\bbT$ below), and let $\pi$ be any subgroup of the symmetric group
$S_n$. Over $\bbC$, the corresponding orbifold by $\pi$  of $n$ copies of $G$ on $G$, 
will be given by the centre of the crossed-product construction
$({}^\tau K_0^{G^n}(G^n))\sdprod \pi=K^0(C(G^n;
\mathcal{K}_\tau)\otimes C^*(G^n))\sdprod \pi$, where $G^n$ acts adjointly
on $G^n$ in the obvious way, and $\pi$ acts by permuting these $n$ factors.
It is tempting to approximate this  geometrically, by guessing that 
the Verlinde algebra of this $\pi$-permutation
orbifold of $G$  is
${}^\tau K_*^{G^n\sdprod \pi}(G^n)$, where $*=0$
and $G^n$ acts adjointly on $G^n$ while $\pi$ acts on the space
$G^n$ by permuting. We need the semi-direct product $G^n\sdprod \pi$ of groups,
rather than direct product, for this to be a group action ($\pi$ will likewise 
act on the subgroup $G^n$ by permuting). For $*=1$, we'd expect a trivial 
$K$-homology group. 

For this idea to work, we would expect that $H^3_{G^n\sdprod \pi}(G^n;\bbZ)$ 
contains a $\bbZ$ which can be identified with $H^3_G(G;\bbZ)$. This can be 
seen as follows. An element of $H^3_G(G;\bbZ)$ corresponds to a $G$-equivariant
bundle $K_\tau$ of compacts on $G$. Taking the product of $n$ of these bundles,
we have a $G^n$ bundle of compacts on $G^n$. For this to make sense under
the action of the permutation group $\pi$, we require choosing the same
bundle (i.e. one element of $H_G^3(G;\bbZ)$) $n$ times. This gives a map from
$H^3_G(G;\bbZ)$ to $H^3_{G^n\sdprod \pi}(G^n;\bbZ)$, which shouldn't be the zero 
map. In section 4.2 we generalise slightly this construction, owing to an extra
large cohomology group.

To see that this is not unreasonable, consider the $\bbZ_2$-permutation
orbifold of the quantum double of finite abelian groups. Let $G$ be a finite 
abelian group, of order $n$ say, write $H=(G\times G)\sdprod\bbZ_2\,$, and
take trivial twist $\sigma\in H^4(BG;\bbZ)\,$. Then $K_1^H(G\times G)$ is
readily seen to vanish, and $K_0^H(G\times G)$ can be computed by writing
the finite set $G\times G$ as the union of the diagonal elements $D=
\{(g,g)\}$ and the off-diagonal, the latter parametrised by a set $F$ of 
$(n^2-n)/2$ $\bbZ_2$-orbit representatives $(g,h)$. Then $K_0^H(D)$ is 
$\|D\|=n$ copies of $K_0^H(pt)=R_H$, which has dimension $2n+(n^2-n)/2$
(the type ${2\atop 1}$- resp. ${1\atop 2}$-representations, using 
terminology of section 2.1). Also,
$K_0^H(F\times H/(G\times G))$ consists of $\|F\|$ copies of $K_0^{G\times G}
(pt)=R_{G\times G}\,$, which has dimension $n^2$. Thus $K_0^H(G\times G)$
has dimension $n(n^2+n)/2+(n^2-n)n^2/2=(n^4+n^2)/2\,$, which matches exactly
the number of primaries in the $\bbZ_2$-permutation orbifold of the quantum 
double of $G$. The only curious aspect of this calculation is that 
$(n^3+n^2)/2$ dimensions are associated to the diagonal $D$, whereas $2n^2$ 
primaries in the orbifold result from the doubling of the $n^2$ fixed points. 
This means the $n^2$ fixed point primaries of $\mathcal{D}(G)$
should not be identified with the $n^2$ diagonal elements of the `space' 
$G\times G$, but this should be clear since they are parametrised 
differently.

We further test this with the $S_2$-permutation orbifold of both $\bbT$ and $SU(2)$.

\subsection{The Verlinde algebra of the 2-torus}\bigskip

Before computing the $S_2$-permutation orbifold of the circle
$\bbT$, it is convenient to compute the Verlinde algebra of $T^2
$ explicitly by $K$-homology. (The result for the loop
group of any torus is quoted in section 3.1.)
Being abelian, the torus $T^2=(\bbR/\sqrt{2}\bbZ)^2$ acts trivially on itself.
The level here can be any
matrix $K=\left(\begin{matrix}k&l\cr m&n\end{matrix}\right)$ with $det\,K\ne 0$.
Judging from section 3.1, there should be no $H^1$-twist.

Let $\Delta$ be the diagonal $(x,x)$ in the torus $T^2$.
The six-term sequence will involve the $K$-holomogy of the circle $\Delta$ and 
the cylinder $T^2\setminus\Delta=\bbR\times \mathbb{T}$, and we will first compute
these from Mayer-Vietoris:
\begin{equation}
\begin{array}{ccccc}
    ^{\tau^{\prime\prime}}K^{T^2}_0(\Delta) & {\longrightarrow} & 0&
   {\longrightarrow} &0 \\[3pt]
    \uparrow  & & & & \downarrow
    \\[2pt]
    2\times K^{T^2}_{1}(\bbR) & \stackrel{\alpha\,}{\longleftarrow} &
    2\times K^{T^2}_{1}(\bbR) &{\longleftarrow} &^{\tau^{\prime\prime}}K^{T^2}_1(\Delta)
    \end{array}
    \end{equation}
\begin{equation}
\begin{array}{ccccc}
    ^{\tau^{\prime}}K^{T^2}_0(\bbR\times \mathbb{T}) & {\longrightarrow} & 
2\times K^{T^2}_{0}(\bbR^2)&
    \stackrel{\beta \,}{\longrightarrow} & 2\times K^{T^2}_{0}(\bbR^2)
    \\[3pt]
    \uparrow  & & & & \downarrow
    \\[2pt]
    0 & {\longleftarrow} &0 &{\longleftarrow}
    &^{\tau^{\prime}}K^{T^2}_1(\bbR\times \mathbb{T})
    \end{array}
    \end{equation}
Now, $K^{T^2}_{0}(\bbR^2)=K^{T^2}_{1}(\bbR)=R_{T^2}=\bbZ[a^{\pm 1},b^{\pm 1}]$.
The map $\alpha$ sends the Laurent polynomials $(p(a,b),q(a,b))$ to 
$(p+q,p+a^{k+l}b^{m+n}q)$, and so we obtain $^{\tau^{\prime\prime}}K_0^{T^2}
(\Delta)=\bbZ[a^{\pm 1},b^{\pm 1}]/(1-a^{k+l}b^{m+n})$ and
 $^{\tau^{\prime\prime}}K_1^{T^2}(\Delta)=0$. Similarly, the map $\beta$ sends 
$(p(a,b),q(a,b))$ to $(p+q,p+a^{k+l}b^{m+n}q)$ (it involves
the same cycle on $T^2$), and so we  obtain $^{\tau^{\prime}}K_1^{T^2}
(\bbR\times \mathbb{T})=\bbZ[a^{\pm 1},b^{\pm 1}]/(1-a^{k+l}b^{m+n})$ and
 $^{\tau^{\prime}}K_0^{T^2}(\bbR\times \mathbb{T})=0$. 

Thus, the six-term sequence becomes
\begin{equation}
\begin{array}{ccccc}
    0 & {\longleftarrow} & ^{K}K^{T^2}_{0}(T^2) &
    {\longleftarrow} & \bbZ[a^{\pm 1},b^{\pm 1}]/(1-a^{k+l}b^{m+n})
    \\[3pt]
    \downarrow  & & & & \uparrow\ \gamma
    \\[2pt]
    0 & {\lori} &
    ^{K}K^{T^2}_1(T^2) &{\lori} &\bbZ[a^{\pm 1},b^{\pm 1}]/(1-a^{k+l}b^{m+n})
    \end{array}
    \end{equation}
where $\gamma([r(a,b)])=[(1-a^{k}b^{m})r]$. This gives
$^{K}K^{T^2}_{1}(T^2)=0$ and
$^{K}K^{T^2}_{0}(T^2)=\bbZ[a^{\pm 1},b^{\pm 1}]/(1-a^{k+l}b^{m+n},1-a^{k}b^{m})\,$, which
as an additive group is isomorphic to the group ring of
$\bbZ^2/{\rm Span}\{\left(\begin{matrix}k+l\cr m+n\end{matrix}\right),\left(
\begin{matrix}k\cr m\end{matrix}\right)\}\,$.
 This is in agreement with Theorem 4.2(i) of 
\cite{FHTii}.

\bigskip
\subsection{The $S_2$-permutation orbifold of the circle}\bigskip

Consider for simplicity the $S_2$-permutation orbifold of $\bbT$ at level $k$.
Over $\bbC$, this will be given by the centre of the crossed-product construction
$({}^\tau K_0^{\bbT\times\bbT}(\bbT\times \bbT))\sdprod S_2=K^0(C(\bbT^2;
\mathcal{K}_\tau)\otimes C^*(\bbT^2))\sdprod S_2$. As mentioned
above, this suggests the geometric approximation 
$K^0(C(\bbT^2;\mathcal{K}_\tau)\otimes C^*(\bbT^2)\sdprod S_2)
={}^\tau K_0^{\bbT^2\sdprod S_2}(\bbT^2)$.

Let $U=T^2\sdprod S_2$. We will mimic as much as possible the previous
calculation. There are only two orbits: the diagonal circle $\Delta=\{(x,x)\}$,
which is fixed by all of $U$; the off-diagonal is the cylinder $(x,\theta)$
with free $S_2$-action $(x,\theta)\mapsto(-x,\theta+\pi)$. 
The level $K$ here should commute
with $S_2$, i.e. $K=\left(\begin{matrix}k&l\cr l&k\end{matrix}\right)$. Again,
$K$ should be invertible, i.e. $|k|\ne |l|$. From the bundle picture, it is clear
this gives an element of $H^3_U(T^2;\bbZ)$. Strictly speaking,
the permutation orbifold of $\bbT_k$ would require taking $l=0$;
nonzero $l$ would correspond to a $S_2$-orbifold of $T^2$ level $K$.

The $K$-homology of $\Delta$ can be computed from Mayer-Vietoris:
\begin{equation}
\begin{array}{ccccc}
    ^{\tau^{\prime\prime}}K^{U}_0(\Delta) & {\longrightarrow} & 0&
   {\longrightarrow} &0 \\[3pt]
    \uparrow  & & & & \downarrow
    \\[2pt]
    2\times K^{U}_{1}(\bbR) & \stackrel{\alpha'\,}{\longleftarrow} &
    2\times K^{U}_{1}(\bbR) &{\longleftarrow} &^{\tau^{\prime\prime}}K^{U}_1(\Delta)
    \end{array}
    \end{equation}
Now, $K^{T^2}_{0}(\bbR^2)=R_{T^2}=\bbZ[a^{\pm 1},b^{\pm 1}]$, while
$K^{U}_{0}(\bbR)=K^{T^2}_{1}(\bbR^2)=0$. 
The representation ring of a double-cover such as $U$, together with its
induction and restriction maps,  is described in section 2.1: 
\begin{equation}K^{U}_{1}(\bbR)=R_{U}=\bbZ[\mu_\Delta^{\pm 1},\mu,\delta]/
(\delta^2=1,\delta \mu=\mu)=\bbZ[\mu_{ij}=\mu_{ji},\delta\mu_{ii}]\,,
\end{equation} 
where $\mu_{10}=\mu$, $\mu_{\Delta}=\mu_{11}$, and $\mu_{ij}\mu_\Delta^k
=\mu_{i+k,j+k}$. Induction from $T^2$ to $U$ takes both 
$a^ib^j$ and $a^jb^i$ (for $i\ne j$) to the two-dimensional irreducible 
representation $\mu_{ij}$, and it takes the $S_2$-fixed
point $a^nb^{n}$ to $\mu_{nn}+\delta\mu_{nn}$
 ($S_2$ acts on $R_{T^2}$ by interchanging $a$ and $b$).
The map $\alpha'$ sends the polynomials $(p(\mu_\Delta,\mu,\delta),
q(\mu_{\Delta},\mu,\delta))$ to 
$(p+q,p+\mu_\Delta^{k+l} q)$, and so we obtain $^{\tau^{\prime\prime}}
K_0^{U}(\Delta)=R_U/(1-\mu_\Delta^{k+l})$ and
 $^{\tau^{\prime\prime}}K_1^{U}(\Delta)=0$. 

The $K$-homology of the cylinder $cyl=T^2\setminus\Delta$ is more
difficult. Write ${}^\tau K_*^U(cyl)=K^*_{S_2}(C(cyl;\mathcal{K}_\tau)\otimes
C^*(T^2))$ (which we can do because $T^2$ is normal in $U$), and split
these representations $C^*(T^2)\simeq C_0(\bbZ^2)$ into the 
$S_2$-fixed points $C_0(\bbZ)$ (the diagonal) and $C_0
(\bbZ^2\setminus\bbZ)$ (the upper and lower triangles). Because $S_2$ fixes
that diagonal, we may write $K^*_{S_2}(C(cyl)\otimes
C_0(\bbZ))={}^{k+l}K_*^{\bbT\times S_2}(cyl)$, and because $S_2$ is normal
in $\bbT\times S_2$ and acts freely on $cyl$, $^\tau K^*_{\bbT\sdprod S_2}(cyl)
\simeq {}^{\tau''} K_{\bbT}^*(cyl/S_2)$ (see (\ref{freenormal})). But 
$cyl/S_2$ is the open M\"obius strip $M\ddot{o}b$. Hence Mayer-Vietoris gives us:
\begin{equation}
\begin{array}{ccccc}
    ^{\tau^{\prime}}K^{\bbT}_0(M\ddot{o}b) & {\longrightarrow} & 
2\times K^{\bbT}_{0}(\bbR^2)&
    \stackrel{\beta' \,}{\longrightarrow} & 2\times K^{\bbT}_{0}(\bbR^2)
    \\[3pt]
    \uparrow  & & & & \downarrow
    \\[2pt]
    0 & {\longleftarrow} &0 &{\longleftarrow}
    &^{\tau^{\prime}}K^{\bbT}_1(M\ddot{o}b)
    \end{array}
    \end{equation}
The map $\beta'$ sends 
$(P(a),Q(a))$ to $(P+Q,P-a^{k+l}Q)$ (the
minus sign comes from the topology of $M\ddot{o}b$), and so we
obtain $^{\tau^{\prime}}K_1^{\bbT}
(M\ddot{o}b)=\bbZ[a^{\pm 1}]/(1+a^{k+l})$ and
 $^{\tau^{\prime}}K_0^{\bbT}(M\ddot{o}b)=0$. 

For the remaining `off-diagonal' part of $C^*(T^2)$, the freeness of this
$S_2$-action identifies $K^*_{S_2}(C(cyl;\mathcal{K}_\tau)\otimes 
C_0(\bbZ^2\setminus\bbZ))$ with 
$K^*(C(cyl;\mathcal{K}_\tau)\otimes C_0(tri))$, where `$tri$' denotes
(say) the lower triangle $\{(m,n)\in\bbZ^2\,|\,m>n\}$. Mayer-Vietoris
applied to the cylinder then gives us
\begin{equation}
\begin{array}{ccccc}
    K^0(C(cyl)\otimes C_0(tri)) & {\longrightarrow} & 
2\times K^{0}(C(\bbR^2)\otimes C_0(tri))&
    \stackrel{\beta'' \,}{\longrightarrow} & 2\times K^{0}(C(\bbR^2)\otimes
C_0(tri)) \\[3pt]
    \uparrow  & & & & \downarrow
    \\[2pt]
    0 & {\longleftarrow} &0 &{\longleftarrow}
    &K^1(C(cyl)\otimes C_0(tri))
    \end{array}
    \end{equation}
As usual, the twist on $cyl$ can be trivialised on the open subsets $\bbR^2$, 
so $K^{0}(C(\bbR^2)\otimes C_0(tri))$ is the off-diagonal part $\mu R_U$
of $R_U$.  The map $\beta''$ sends $(p,q)\in\mu R_U\oplus \mu R_U$
to $(p+q,p+\mu_\Delta^{k+l}q)$, by analogy with section 4.1.
Therefore we find $K^0(C(cyl;\mathcal{K}_\tau)\otimes C_0(tri))={\rm ker}\,\beta''=0$
and $K^1(C(cyl;\mathcal{K}_\tau)\otimes C_0(tri))={\rm coker}\,\beta''=\mu R_U/(
1-\mu_\Delta^{k+l})$.

The six-term sequence then tells us how to obtain $^\tau K_*^U(cyl)$:
\begin{equation}
\begin{array}{ccccc}
    0 & {\longleftarrow} & ^{\tau}K^{U}_{0}(cyl) &
    {\longleftarrow} & 0
    \\[3pt]
    \downarrow  & & & & \uparrow\ 
    \\[2pt]
\bbZ[\mu_\Delta^{\pm 1}]/(1-\mu_\Delta^{k+l}) & {\lori} &
    ^{\tau}K^{U}_1(cyl) &{\lori} &\mu R_U/(1-\mu_\Delta^{k+l})
    \end{array}
    \end{equation}
Hence $^{\tau}K^{U}_{0}(cyl)=0$ and $^{\tau}K^{U}_{1}(cyl)=\bbZ[
\mu_\Delta^{\pm 1}]/(1-\mu_\Delta^{k+l})\oplus \mu R_U/(1-
\mu_\Delta^{k+l})$, where the left summand is the submodule
and the right one is its quotient into $K_1$.

Finally, we again use the six-term to compute the desired $K$-homology:
\begin{equation}
\begin{array}{ccccc}
    0 & {\longleftarrow} & ^{\tau}K^{U}_{0}(T^2) &
    {\longleftarrow} & R_U/(1-\mu_\Delta^{k+l})
    \\[3pt]
    \downarrow  & & & & \uparrow\ \gamma'
    \\[2pt]
0
& {\lori} &    ^{\tau}K^{U}_1(T^2) &{\lori} &
\bbZ[\mu_\Delta^{\pm 1}]/(1-\mu_\Delta^{k+l}) \oplus
\mu R_U/(1-\mu_\Delta^{k+l})\mu R_U
    \end{array}
    \end{equation}
The map $\gamma'$ first projects $^\tau K_1^U(cyl)$ to $\mu R_U/(1-
\mu_\Delta^{k+l})$, killing the submodule $\bbZ[
\mu_\Delta^{\pm 1}]/(1-\mu_\Delta^{k+l})$, and then (by analogy with
section 4.1) multiplies $p\in \mu R_U$ by $(1-\mu_{k-l}\mu_\Delta^l)$,
since $\mu_{k-l}\mu_\Delta^l=\,$D-Ind$_{T^2}^U a^kb^l$. Thus we obtain
the final answer:
 \begin{eqnarray} \label{tpermorb}
 ^{\tau}K^{{U}}_{0}(T^2)&=&\,R_U/(1-\mu_\Delta^{k+l},\mu(1-\mu_{k-l}
\mu_\Delta^l))\,,\\
 ^{\tau}K^{{U}}_{1}(T^2)&=&\,\bbZ[\mu_\Delta^{\pm 1}]/(1+
\mu_\Delta^{k+l})\,.
 \end{eqnarray}

Consider for simplicity $l=0$; then this $K$-homology does not quite recover the
permutation orbifold of the circle at level $k$. In particular, RCFT tells us 
that the $k$
diagonal points $(a,a)$, being fixed by $S_2$, should get doubled, while 
$(a,b)$ and $(b,a)$ for $a\ne b$ should be identified.  
We discuss this further in section 7.

\subsection{The $S_2$-permutation orbifold of $SU(2)$} 

Now let $G=SU(2)$, $\bbT$ be the diagonal matrices, $H=(G\times G)\sdprod
\bbZ_2$, and $U=(\bbT\times \bbT)\sdprod\bbZ_2$ as before, where $\bbZ_2$ acts
on $G\times G$ and $\bbT\times\bbT$ by permuting.
The orbit analysis for $H$ on $G\times G$ is reminiscent of that of
$G$ on $G$ given at the beginning of section 1.4.
The poles $(\pm I,\pm I)\in G\times G$ form three $H$-orbits: $\mathcal{O}_H=
(I,I)\cup(-I,-I)$ is fixed by everything, while $\mathcal{O}_{G\times G}
=(I,-I)\cup(-I,I)=H/(G\times G)$ has stabiliser $G\times G$. Let $gen=SU(2)
\setminus\{\pm I\}=\bbR\times G/\bbT$ be the generic points in $G$; then
the orbits $\mathcal{O}_{G\times \bbT}=\pm I\times gen\cup gen\times \pm I=
2\times\bbR
\times H/(G\times\bbT)$ have stabiliser $G\times \bbT$. The points $(x,y)\in
gen\times gen$ are of two kinds: those with $x,y$ conjugate form the orbits
$\mathcal{O}_{U}=\bbR\times H/U$,
while the remainder form $\mathcal{O}_{T^2}=\bbR^2\times H/T^2$.
Figure 7 gives the resulting picture.

\bigskip 
\begin{figure}[tb]
\begin{center}
\epsfysize=1.5in \centerline{\epsffile{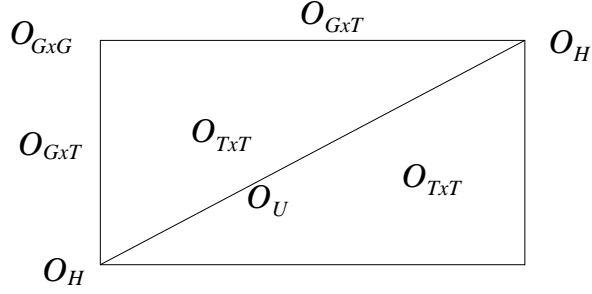}}
\caption{Figure 7: The orbits of $H$ on $G\times G$}
\end{center}
\end{figure}

Using the method of section 2.1, we easily find that the irreducible 
$H$-representations are $\rho_{ij}=\rho_{ji}$ ($i,j\ge 0$), 
$D$, and $D\rho_{ii}$, where $1=\rho_{00}$ and $D$ are one-dimensional,
and these obey $D^2=1$ and 
$D\rho_{ij}=\rho_{ij}$ for $i\ne j$. The induction Ind$_{G\times G}^H
(\sigma_i\sigma'_j)$ equals $\rho_{ij}=\rho_{ji}$ for $i\ne j$, and $\rho_{ii}
+D\rho_{ii}$ otherwise. Analogous comments for the 
$U$-representations were discussed last subsection. We will shortly need
the Dirac inductions D-Ind$_{G\times\bbT}^H={\rm Ind}_{G\times G}^H\circ\, $D-Ind$_{G\times
\bbT}^{G\times G}$ (see section 2.1 for D-Ind$_{\bbT}^G$), and D-Ind$_U^H$.
The latter sends $\mu_\Delta^j$ (resp. $\delta\mu_\Delta^j$) 
to $\rho_{|j|-1,|j|-1}$ (resp. $D\rho_{|j|-1,|j|-1}$) provided $j\ne 0$,
and kills both $1$ and $\delta$; it sends $\mu_{ij}$ for $i\ne j$
to $sgn(ij)\,\rho_{|i|-1,|j|-1}$, unless $i=-j$ in which case
it yields $-\rho_{|i|-1,|i|-1}(1+D)$.

Let $k$ be the (unshifted) level of $SU(2)$, and $\tau$ the corresponding 
element of $H_H^3(G\times G;\bbZ)$, as explained in the second paragraph of
section 4. The twist won't survive on any of these orbit spaces, so their
$K$-homology groups can be simply written down: the nonzero ones are
$K^H_0(\mathcal{O}_H)=2\times R_H$,
$K^H_0(\mathcal{O}_{G\times G})=R_{G\times G}$, $K_1^H(\mathcal{O}_{G\times
\bbT}) =2\times R_G\otimes R_\bbT$, $K_1^H(\mathcal{O}_U)=R_U$, and $K_0^H(
\mathcal{O}_{T^2})=R_{T^2}=\bbZ[a^{\pm 1},b^{\pm 1}]$.

First we use the six-term sequence to glue the compact set $\mathcal{O}_H\cup 
\mathcal{O}_{G\times G}$ to $\mathcal{O}_{G\times T}\cup \mathcal{O}_U$:
\begin{equation}
\begin{array}{ccccc}
    0 & {\longleftarrow} & ^{\tau}K^{H}_{0}(\mathcal{O}_{T^2}^c) &
    {\longleftarrow} & 2\times R_H\oplus R_{G\times G}
    \\[3pt]
    \downarrow  & & & & \uparrow\ \alpha
    \\[2pt]
0
& {\lori} &    ^{\tau}K^{H}_1(\mathcal{O}_{T^2}^c) &{\lori} &
2\times R_{G\times\bbT}\oplus R_U
    \end{array}\ ,
    \end{equation}
where $\alpha:R_{G\times \bbT}\oplus R_{G\times\bbT}\oplus R_U\rightarrow
R_H\oplus R_H\oplus R_{G\times G}$ is given by $\alpha(p_1,p_2,q)=\,$D-Ind$
(q+b^kp_1,-\mu_\Delta^{2k}q+p_2,p_1+b^kp_2)$. The sign in front of
the $\mu_\Delta$ comes from considering orientations of the edges; the
powers $k$ or $2k$ are clear from the bundle picture, where traversing an
edge crosses patches in the bundle once or twice respectively.

A tedious but straightforward calculation shows that 
 \begin{eqnarray}
{}^\tau K_0^H(\mathcal{O}_{T^2}^c)&=&\mathrm{coker}\,\alpha=R_H/
(\rho_{ij}-\rho_{2k+i,i},(1+D)\rho_{ii}-\rho_{2k+i,i},\\
&&\qquad(1-D)(\rho_{ii}-\rho_{
i+4k,i+4k}),(1-D)(\rho_{ii}-\rho_{-i+4k-2,-i+4k-2}))\,,\nonumber\\
{}^\tau K_1^H(\mathcal{O}_{T^2}^c)&=&\mathrm{ker}\,\alpha
=\mathrm{Span}\{(2\sigma_{i+k}b^j+2\sigma_{i-k}b^{j-2k},
2\sigma_{i+k}b^{-j-k}+2\sigma_{i-k}b^{-j+k},\nonumber\\
&&\qquad -(1-\delta)(\mu_{i-k,j-k}+
\mu_{-i-k,-j-k}))\}\,.\label{kercalc}
 \end{eqnarray}
The 2's in (\ref{kercalc}) are spurious for $i\ne j$, since then $\delta$ acts 
like 1. Determining (\ref{kercalc}) is the more difficult: first solve
the simpler problem by restricting $R_{G\times\bbT}$ and $R_U$ to $R^{T^2}$.
This then requires finding all $c,c'\in\mathrm{ker\,}$(D-Ind$_{T^2}^H)\,$,
$c''\in\mathrm{ker\,}$(D-Ind$_{T^2}^{G\times G})\,$, such that $a^{-k}b^{-2k}
-a^kb^{2k}$ divides $b^kc'+b^{-k}c-c''$. This can be done by 
(anti-)symmetrising with respect to $a\leftrightarrow a^{-1},b\leftrightarrow
b^{-1}$.

The desired $K$-homology is now obtained from the six-term sequence 
by gluing in $\cO_{T^2}$:
 \begin{equation}
\begin{array}{ccccc}
    R_{T^2} & {\longleftarrow} & ^{\tau}K^{H}_{0}(G\times G) &
    {\longleftarrow} &\mathrm{coker}\,\alpha 
    \\[3pt]
    \beta\ \downarrow  & & & & \uparrow
    \\[2pt]
\mathrm{ker}\,\alpha
& {\lori} &    ^{\tau}K^{H}_1(G\times G) &{\lori} & 0
    \end{array}\ .
    \end{equation}
Presumably, the map $\beta$ is an isomorphism, sending irreducible 
$T^2$-representation
$a^ib^j$ to eg. the $(i,j)$-element in (\ref{kercalc}).
This would mean we obtain
 \begin{eqnarray} \label{su2permorb}
 ^{\tau}K^{{H}}_{0}(G\times G)&=&R_H/
(\rho_{ij}-\rho_{2k+i,i},(1+D)\rho_{ii}-\rho_{2k+i,i},\nonumber\\
&&\qquad(1-D)(\rho_{ii}-\rho_{
i+4k,i+4k}),(1-D)(\rho_{ii}-\rho_{-i+4k-2}))\nonumber\\
&\simeq&\bbZ^{k(k+7)/2}\,,\\
 ^{\tau}K^{{H}}_{1}(G\times G)&=&0\,.
 \end{eqnarray}
We interpret this result in section 7.

\bigskip
\section{The ${D}_4$ modular invariant of $SU(2)$} 

We will realise this modular invariant through the conformal embedding
$SU(2)_4\rightarrow SU(3)_1$. From the finite-dimensional perspective,
this is the 
two-to-one projection of $SU(2)$ onto $SO(3)$. Explicitly, this representation
 is given as follows. We can parametrise the 3-sphere $SU(2)$ by
\begin{equation} \label{su2param}
\left(\begin{matrix}\gamma&\delta\cr -\overline{\delta}&\overline{\gamma}
\end{matrix}\right)\, ,
\end{equation}
where $\gamma,\delta\in \bbC$ satisfy $|\gamma|^2+|\delta|^2=1$. The canonical
choice of maximal torus is $\bbT$, i.e. $\delta=0$, and the other component
of the canonical $O(2)$ is $\gamma=0$. Then
\begin{equation}\label{SU2SU3}
{\cal R}^{(3)}(\gamma,\delta):=
{\cal R}^{(3)} \left(\begin{matrix}\gamma&\delta\cr -\overline{\delta}&\overline{\gamma}\end{matrix}\right)=
{\rm Re}
\left(\begin{matrix}\scriptstyle \gamma^2-\i\delta^2& \scriptstyle-\i\gamma^2+\delta^2
&\scriptstyle 2\xi\gamma\delta\cr \scriptstyle \i\gamma^2+\delta^2&\scriptstyle
\gamma^2+\i\delta^2&\scriptstyle 2\i{\xi}\gamma\delta\cr\scriptstyle
-2\xi\overline{\gamma}\delta&\scriptstyle -2\i\xi\overline{\gamma}\delta&
\scriptstyle|\gamma|^2-|\delta|^2\end{matrix}\right)\, ,
\end{equation}
where $\xi=\exp[\pi\i/4]$ and `Re' denotes the real part.
Note that ${\cal R}^{(3)}(e^{{\rm i}\theta},0)={\rm diag}(R_{2\theta},1)$
and ${\cal R}^{(3)}(0,e^{{\rm i}\theta})={\rm diag}(R'_{2\theta},1)$,
where $R_\theta,R'_\theta$ are the rotation resp.\ reflection matrices
\begin{equation}
R_\theta=\left(\begin{matrix}\cos\,\theta&\sin\,\theta\cr-\sin\,\theta&
\cos\,\theta\end{matrix}\right)\, ,\qquad
R'_\theta=\left(\begin{matrix}\sin\,\theta&\cos\,\theta\cr\cos\,\theta&
-\sin\,\theta\end{matrix}\right)\, .
\end{equation}
Of course, $R_\theta R_\phi=R_{\theta+\phi}$, $R_\theta R'_\phi=R'_{\theta+
\phi}=R'_\phi R_{-\theta}$, $R'_\theta R'_\phi=R_{\theta-\phi}$. 
Write $T={\cal R}^{(3)}(*,0)$ and $T'={\cal R}^{(3)}(0,*)$ for the
two components of the image of $O(2)$.

Throughout this section we write $G$ for $SU(2)$ and $\overline{G}$ for $SO(3)$.

Even though $G$ acts by first projecting to $SO(3)$, the groups $^\tau
K_*^G(SU(3))$ and $^\tau K_*^{SO3}(SU(3))$ are different.
Indeed, the $K$-homology ${}^\tau K_*^{SO3}(SU(3))$ is
easy to compute using the Hodgkin spectral sequence (section 3.3). In particular,
write $R$ for $R_{SU3}$; then as $R$-modules we have $R_{SO3}
=R/(\mu_{10}-\mu_{01})$  and $Ver_1(SU(3))
=R/(\mu_{20},\mu_{02})$, where here
$\mu_{ij}$ denotes the $SU(3)$-representation with highest-weight $(i,j)$.
As general facts we have
Tor$_0^R(A,B)=A\otimes_R B$ and Tor$_1^R(R/I,R/J)=(I\cap J)/(IJ)$, valid
for any ring $R$, any $R$-modules $A,B$, and any ideals $I,J$.
In this example Tor$_p$ will vanish for $p>1$, since $R_{SO3}$ has a free
resolution of length 1. We quickly find that
Tor$_0^{R}(R_{SO3},Ver_1(SU(3)))\simeq \bbZ$ and Tor$_1^R(R_{SO3},
Ver_1(SU(3)))\simeq \bbZ$, with the $R$-module structure given in both cases 
by $\mu_{10}$ and $\mu_{01}$ acting like 1. Hence ${}^\tau K^0_{SO3}
(SU(3))\simeq \bbZ\simeq {}^\tau K^1_{SO3}(SU(3))$ for the appropriate
twist $\tau$. Of course the $K$-homology then follows by Poincar\'e duality.
This differs from (\ref{d4answer}) in both the absence of torsion and
a different module structure (the generator $\sigma_3$ of $R_{SO3}$ 
acts as $+1$ in $^\tau K^{SO3}_*(SU(3))$, and as $-1$ in 
$^\tau K^{SU2}_*(SU(3))$. The $K$-theory ${}^\tau K^*_{SO3}(G_2)$ for the $E_8$
modular invariant (corresponding to the conformal embedding $SU(2)_{28}
\rightarrow G_{2,1}$) can be computed similarly. 
The relevance of this $K$-homology to conformal field theory isn't so clear
to us though.

\subsection{The orbit analysis}

\bigskip

We need to understand the orbits of the adjoint action of $\overline{G}=SO(3)$ 
on $SU(3)$. By `Stab$_{\overline{G}}{\bf B}$', we mean the set of all ${\bf A}\in SO(3)$
commuting with ${\bf B}$. It is convenient to write ${\bf A}\in \overline{G}$ and
${\bf B}\in SU(3)$ in block form as
\begin{eqnarray}
{\bf A}&=&\left(\begin{matrix}A&a\cr a'&\alpha\end{matrix}\right)\,,\label{blockA}\\
{\bf B}&=&\left(\begin{matrix}B&b\cr b'&\beta\end{matrix}\right)\,,\label{blockB}
\end{eqnarray}
where $A,B$ are $2\times 2$ matrices and $\alpha,\beta$ are numbers. A simple
observation is that if $a=0$ in (\ref{blockA}), then $|\alpha|=1$, so $a'=0$ and hence 
$A\in O(2)$, with det$\, A=\alpha=\pm 1$. 

\medskip\noindent{\bf Lemma 1.} {\it If ${\bf B}\in SU(3)$ commutes with some
${\bf A}\in \overline{G}$ with (possibly  infinite) order $n>2$, then Stab$_{\overline{G}}
({\bf B})$ contains a maximal torus of $\overline{G}$}.

\medskip\noindent{{\it Proof:}}  Without loss of generality (by conjugating
${\bf A}$ and ${\bf B}$ simultaneously by $\overline{G}$) we can take ${\bf A}\in T$,
i.e.\ ${\bf A}={\rm diag}(R_\theta,1)$ for some $\theta$. Then ${\bf A}$ commutes
with ${\bf B}$ in (\ref{blockB}) iff $R_\theta B=BR_\theta$, $R_\theta b=b$,
and $b'R_\theta=b'$. But $R_\theta\ne I$ can have no eigenvector with
eigenvalue 1, so $b=0$ and $b'=0$. Hence $B\in U(2)$. $R_\theta B=BR_\theta$
requires $B$ and $R_\theta$ to be simultaneously diagonalisable; since
$R_\theta\ne\pm I$, we get $B=e^{\i\psi}R_\phi\,$. But any such ${\bf B}$
will commute with all of $T$.\qquad QED\medskip

Therefore the finite stabilisers Stab$_{\overline{ G}}({\bf B})$ are the finite
groups of A-D-E type with exponent $\le 2$: namely, the trivial group 1,
the cyclic group $C_2$, and the dihedral group $D_2=C_2\times C_2$.
Lifting these to $G$ yields the double-covers $\bbA_1=C_2$,
$\bbA_3=C_4$, and  the quaternionic group $\bbD_4=Q_4$.

\medskip\noindent {\bf Case 1}: {\it Orbits with infinite stabiliser}.
By Lemma 1, without loss of generality (conjugating by $\overline{G}$ if
necessary) we can take ${\bf B}$ to be of form
\begin{equation}\label{Binf}
{\bf B}=\left(\begin{matrix}e^{\i\psi}R_\phi&0\cr 0&e^{-2\i\psi}\end{matrix}\right)
\end{equation}
for some angles $\psi,\phi$, so Stab$_{\overline{G}}({\bf B})$ contains $T$.
Suppose now that such a matrix commutes with some ${\bf A}\not\in T$ in
(\ref{blockA}). Then $R_\phi A=AR_\phi$, $R_\phi a=e^{-3\i\psi}a$, $a'R_\phi
=e^{-3\i\psi}a'$. If $a\ne 0$, then $a$ would be a {\em real} eigenvector
of $R_\phi$ with eigenvalue $e^{-3\i\psi}$; this would require $e^{\i\psi}$
to be a sixth root of 1. This quickly forces ${\bf B}$ to be a scalar matrix:
${\bf B}=\omega^iI$ for some $i$, where $\omega=e^{2\pi \i/3}$. These are
the three fixed points, each with stabiliser $\overline{G}$.

Otherwise, $a=0$ so ${\bf A}={\rm diag}(A,\alpha)$ where $A\in T'$ and $\alpha=
{\rm det}(A)=-1$. Then $R_\phi A=AR_\phi$ forces $R_\phi=R_{-\phi}$, i.e.\
$R_\phi=\pm I$. Therefore we may take the orbit representative to be
${\bf B}={\rm diag}(e^{\i\psi}I,e^{-2\i\psi})$. But three values of $\psi$ recover
the fixed points $\omega^iI$. The others all yield ${\bf B}$ with Stab$_{\overline
G}({\bf B})=T\cup T'\simeq O(2)$. Because only $T\cup T'$ in $\overline{G}$
stabilises $T\cup T'$, we know that each $\psi$ corresponds to a distinct
orbit.

The remaining ${\bf B}$ in (\ref{Binf}) will have Stab$_{\overline G}({\bf B})=
T$. Note that the parameter values $(\psi+\pi,\phi+\pi)$ and
$(\psi,\phi)$ correspond to identical ${\bf B}$ and should be identified.
Also, conjugating by ${\bf A}\in T'$ sends $(\psi,\phi)$ to $(\psi,-\phi)$,
so these correspond to the same $\overline{G}$-orbits and should be identified.
Since the normaliser of $T$ in $\overline{G}$ is $T\cup T'$, these are all
the redundancies we need to consider. A fundamental domain for this is
then $0\le \psi<\pi,0<\phi<\pi$; since $(0,\pi)\sim(\pi,\pi-\phi)$ the
resulting surface, parametrising the orbit representatives, is an open
M\"obius strip.

\medskip\noindent {\bf Case 2}: {\it Finite stabilisers $\ne 1$}.
By Lemma 1, any such stabiliser contains an order-2 element, which
without loss of generality we may take to be ${\bf A}_1={\rm diag}(-I,1)$. Any
${\bf B}$ commuting with ${\bf A}_1$ will look like
\begin{equation}\label{Bfin}
{\bf B}=\left( \begin{matrix}B&0\cr 0&\beta\end{matrix}\right)\,,
\end{equation}
for some $B\in U(2),\beta=\overline{{\rm det}(B)}$.

If anything else lies in the stabiliser, then it must generate $C_2\times C_2$.
Without loss of generality (up to conjugation by $T\cup T'$) it can be
taken to be ${\bf A}_2={\rm diag}(1,-1,-1)$ (the easiest way to see this, as with
any other statements we make about $\overline{G}$, is to lift to $G$).
Requiring ${\bf B}$ to commute with ${\bf A}_2$ forces ${\bf B}={\rm diag}(
e^{\i\theta_1},e^{\i\theta_2},e^{\i\theta_3})$, where $\sum_i\theta_i\equiv 0$
(mod $2\pi$). Now, the normaliser of $\langle {\bf A}_1,{\bf A}_2\rangle$ in
$\overline{G}$ is $\langle {\bf A}_1,{\bf A}_2,{\rm diag}(-1,\left(\begin{matrix}0&1\cr 1&0\end{matrix}
\right)),{\rm diag}(\left(\begin{matrix}0&1\cr 1&0\end{matrix}\right),-1)\rangle$, an extension
of $C_2\times C_2$ by $S_3=\mathrm{Aut}(C_2\times C_2)$. Conjugating by this normaliser
means that we should identify the $S_3$-permutations of $\theta_i$, since they
lie in the same $\overline{G}$-orbit. Moreover, the fixed points of this 
$S_3$-action (namely the three diagonals $\theta_i=\theta_j$) have a stabiliser
larger than $C_2\times C_2$, and hence must fall into Case 1. Deleting 
those diagonals from the $\theta_i$-torus yields six disconnected components
(namely the six triangles with vertices at $\theta_i\in\{0,2\pi/3,4\pi/3\}$).
$S_3$ permutes these six triangles, so the $\overline{G}$-orbits with stabiliser
$C_2\times C_2$ are parametrised by any one of those (open) triangular regions.

Now, consider  the orbits with stabiliser Stab$_{\overline G}({\bf B})=\langle
{\bf A}_1\rangle\simeq C_2$, where ${\bf B}$ is as in (\ref{Bfin}). The
normaliser of $\langle{\bf A}_1\rangle$ is $T\cup T'$ (as always, this is
easiest to see by lifting to $G$). Diagonalising $T$ makes the $T$-action
clearer, so parametrise ${\bf B}$ by 
\begin{equation}
{\bf B}=\left(\begin{matrix}\frac{1}{\sqrt{2}}\left(\begin{matrix}1&\i\cr \i&1
\end{matrix}\right)&0\cr
0&1\end{matrix}\right)\left(\begin{matrix}e^{\i\psi}\left(\begin{matrix}\gamma
&\delta\cr -\bar{\delta}
&\bar{\gamma}\end{matrix}\right)&0\cr 0&e^{-2\i\psi}\end{matrix}\right)
\left(\begin{matrix}\frac{1}{\sqrt{2}}
\left(\begin{matrix}1&\i\cr \i&1\end{matrix}\right)&0\cr 0&1\end{matrix}\right)
^{-1}\,,
\end{equation}
where $|\gamma|^2+|\delta|^2=1$. We must identify $(\psi+\pi,-\gamma,-\delta)$
with $(\psi,\gamma,\delta)$, since they correspond to the same ${\bf B}$.
Conjugating by $T$ tells us we can replace $\delta$ with $|\delta|$ (as it
will have the same $\overline{G}$-orbit). Finally, conjugating by $T'$ identifies
$(\psi,\overline{\gamma},|\delta|)$ with $(\psi,\gamma,|\delta|)$. We must
exclude $\delta=0$ and $\gamma\in \bbR$ as these will have an enhanced stabiliser.
Thus the $\overline{G}$-orbits with $C_2$ stabiliser are parametrised by
$(\psi,\phi,r)$ where $0\le \psi\le \pi$, $0<\phi<\pi$, $0<r<1$, and the
boundaries $(0,\phi,r)$ and $(\pi,\pi-\phi,r)$ are identified. This forms
the direct product of the open M\"obius with the interval (0,1).

\medskip\noindent{\bf Case 3}: {\it Trivial stabiliser.} Being unitary, ${\bf B}$ is diagonalisable.
First, note that Stab$_{\overline G}({\bf B})=1$ iff ${\bf B}$ has no real eigenvector,
since rotating such an eigenvector to $(0, 0, 1)^t$ 
amounts to conjugating ${\bf B}$ into a matrix of form (\ref{Bfin}). This
forces all eigenvalues of ${\bf B}$ to be distinct (a 
complex plane in $\bbC^3$ necessarily contains a nonzero 
real vector). Thus the eigenvalues of ${\bf B}$, which are constant on any 
$\overline{G}$-orbit, will lie on the torus minus the three diagonals, quotiented
by an $S_3$-action (since the eigenvalues come unordered). We know (from
case 2) this to be an open 
triangular region. ${\bf B}$ is then uniquely determined once we choose an ordered
triple $(\bbC v_1,\bbC v_2,\bbC v_3)$ of orthogonal complex lines in $\bbC^3$
(corresponding to the three eigenspaces). $\overline{G}$ acts on this triple by simultaneously
acting on each component.

We can always find in the complex line $\bbC v_1$ a vector (call it $v_1$)
with norm 1, with $\|{\rm Re}(v_1)\|\ge \|{\rm Im}(v_1)\|$, and with ${\rm Re}(v_1)$ orthogonal
to ${\rm Im}(v_1)$. Assume first that $x:=\|{\rm Re}(v_1)\|> \|{\rm Im}(v_1)\|$. Then $v_1$
is uniquely determined up to multiplication by $\pm 1$. We can then use
$SO(3)$ to simultaneously rotate ${\rm Re}(v_1)$ to $(x, 0, 0)^t$ and 
${\rm Im}(v_1)$ to $(0, y, 0)^t$ where $y=\sqrt{
1-x^2}$, so without loss of generality we can take
\begin{equation}\label{v123}
v_1=\left(\begin{matrix}x\cr {\rm i} y\cr 0\end{matrix}\right)\, ,\
v_2=\left(\begin{matrix}{\rm i}yw\cr xw\cr z\end{matrix}\right)\, ,\
v_3=\left(\begin{matrix}{\rm i}y{\bar z}\cr x{\bar z}\cr -{\bar w}\end{matrix}\right)\,,
\end{equation}
where $(w,z)\in \bbC\bbP^1$ other than $(1,0)$ or $(0,1)$. We can
uniquely rescale $(w,z)$ to $(re^{\i\theta},\sqrt{1-r^2})$ where $0<r<1$.

The remaining possibility is when $\|{\rm Re}(v_1)\|= \|{\rm Im}(v_1)\|=\frac{1}{\sqrt{2}}$.
We can again require (\ref{v123}) to hold, where $x=y=\frac{1}{\sqrt{2}}$
The difference here is that when we
simultaneously rotate the $v_i$ by diag$(R_\theta,1)$, the effect 
on $v_1$ can be undone by rescaling it by $e^{-{\rm i}\theta}$.
This remaining $\bbT$ freedom means we can take $w,z$ to both be positive.
Thus these $v_i$ are determined up to a single real parameter $0<r<1$.
Recalling the parametrisation $(x,r,\theta)$ of the previous paragraph, this 
special case corresponds  to the limit $x\rightarrow\frac{1}{\sqrt{2}}$,
and so together they form a sphere with one point removed (parametrised by
$x,\theta$) times the interval $(0,1)$ (parametrised by $r$).

\bigskip Thus the orbit spaces of the adjoint action of $G$ on $SU(3)$ are: 

\begin{itemize}

\item[$\cO_{SU2}$:] three fixed points, with stabiliser $G=SU(2)$;

\item[$\cO_{O2}$:] the space $3\times \bbR\times G/O(2)\simeq
3\times \bbR\times \bbP\bbR^2$, with 
stabiliser $O(2)$;

\item[$\cO_{\bbT}$:] $M\ddot{o}b\times G/\bbT\simeq M\ddot{o}b\times S^2$, with
stabiliser $\bbT$, where $M\ddot{o}b$ denotes the open M\"obius strip;
       
\item[$\cO_{\bbD 4}$:] $\bbR^2\times G/\bbD_4$, with stabiliser $\bbD_4\simeq Q_4\,$;

\item[$\cO_{\bbA 3}$:] $\bbR\times M\ddot{o}b\times G/\bbA_3$, with stabiliser $\bbA_3\simeq
C_4\,$;

\item[$\cO_{\bbA 1}$:] $\bbR^5\times G/\bbA_1$, with stabiliser $\bbA_1\simeq
C_2\,$.
\end{itemize}

They're placed in such an order that $\cO_{SU2}$ is compact; $\cO_{SU2}\cup\cO_{O2}$ 
is compact; $\cO_{SU2}\cup\cO_{O2}\cup\cO_{\bbT}$ is compact; etc. We need this 
for our six-term exact sequences. This compactness is clear from the 
descriptions of the orbits given above. 

For later use, we need to see more clearly some of this topology.
The three $\bbR$'s in $\cO_{O2}$ are three arcs joined at the three points of $\cO_{SU2}$, 
to form a circle $S^1$, so think of $\cO_{SU2}\cup\cO_{O2}$ as three sausages linked 
in a circle. The boundary of the M\"obius strip in $\cO_{\bbT}$ is that circle 
$S^1$, so in $\cO_{SU2}\cup\cO_{O2}\cup\cO_{\bbT}$ we imagine gluing the circle of 
$\cO_{SU2}\cup\cO_{O2}$
(doubly coiled) to the M\"obius strip, forming a closed M\"obius
strip $\overline{M\ddot{o}b}$; to each point on $M\ddot{o}b$ we place a sphere, and to
each point on $\partial M\ddot{o}b=S^1$ we place a projective plane of varying 
radius. This is depicted in Figure 8. To the boundary $S^1$ of the M\"obius strip,
we glue the disc $\bbR^2$
of $\cO_{\bbD 4}$, forming a projective plane. So the boundary of $\cO_{\bbD 4}$ is
$\cO_{SU2}\cup\cO_{O2}$. As we head in $\cO_{\bbD 4}$ to a boundary point on the three 
arcs, two of the eigenvalues of ${\bf B}$ become equal, which selects one of 
the three order-2 subgroups of $C_2\times C_2$.

\bigskip 
\begin{figure}[tb]
\begin{center}
\epsfysize=1.5in \centerline{\epsffile{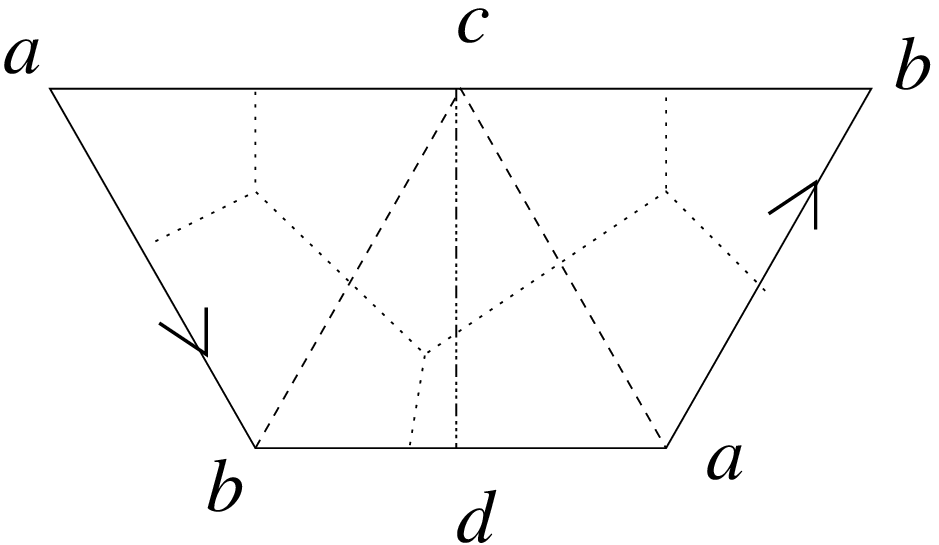}}
\caption{The closed M\"obius strip $M=\cO_G\cup\cO_{O2}\cup
\cO_{\bbT}$}
\end{center}
\end{figure}

Intriguingly, the quaternionic group $Q_4$ is given name `$\mathbb{D}_4$' in the McKay
A-D-E correspondence for finite subgroups of $SU(2)$; ${D}_4$ is also the name given
to this $SU(2)$ modular invariant, in the Cappelli-Itzykson-Zuber
A-D-E classification \cite{CIZ}.

% ${\rm \check C ech}$ cohomology (using
%12 open sets) gives $H^1_G(\cO_{SU2}\cup\cO_{O2};\bbZ_2)\simeq H^1(\cO_{SU2}\cup
%\cO_{O2};\bbZ_2)\simeq \bbZ_2^4$, while Mayer-Vietoris quickly tells us
%$H^1_G(\cO_{SU2}\cup\cO_{O2}\cup\cO_{\bbT};\bbZ_2)$ is either $\bbZ_2^4$ or $\bbZ_2^5$.

%\begin{itemize}

%\item[$\cO_{SU2}$:]  $H^1_G(3;\bbZ_2)\simeq 0$ and  $H^3_G(3;\bbZ)\simeq 0$;

%\item[$\cO_{O2}$:] $H^1_G(3\times \bbR\times G/O(2);\bbZ_2)\simeq 3\times 
%H^1_{O2}(pt;
%\bbZ_2)\simeq \bbZ_2^3$ and $H^3_G(3\times \bbR\times G/O(2);\bbZ)\simeq
%3\times H^3_{O2}(pt;\bbZ)\simeq\bbZ_2^3$;

% \item[$\cO_{\bbT}$:] $H^1_G(M\ddot{o}b\times G/\bbT;\bbZ_2)\simeq H^1_{\bbT}(\mathbb{T};\bbZ_2)
% \simeq \bbZ_2$, and $H^3_G(M\ddot{o}b\times G/\bbT;\bbZ)\simeq H^3_{\bbT}(\mathbb{T};\bbZ)
% \simeq\bbZ$;

% \item[$\cO_{\bbD 4}$:] $H^1_G(\bbR^2\times G/\bbD_4;\bbZ_2)
% \simeq H^1_{\bbD 4}(pt;\bbZ_2)\simeq \bbZ_2^2$ and $H^3_{{G}}(\bbR^2\times
% {G}/\bbD_4;\bbZ)\simeq H^3_{\bbD 4}(pt;\bbZ)\simeq 0$;

% \item[$\cO_{\bbA 3}$:] $H^1_{{G}}(\bbR\times M\ddot{o}b\times {G}/\bbA_3;\bbZ_2)\simeq H^1_{\bbA 3}(\mathbb{T};
% \bbZ_2)\simeq \bbZ_2^2$ and $H^3_{{G}}(\bbR\times M\ddot{o}b\times G/\bbA_3;\bbZ)
% \simeq H^3_{\bbA 3}(\mathbb{T};\bbZ)\simeq \bbZ_4$;

% \item[$\cO_{\bbA 1}$:] $H^1_G(\bbR^5\times G/{\bbA_1};\bbZ_2)\simeq
%H^1_{\bbA 1}(pt;\bbZ_2)=\bbZ_2$
%and $H^3_G(\bbR^5\times G/{\bbA_1};\bbZ)\simeq H^3_{\bbA 1}(pt;\bbZ)= 0$.\end{itemize}

\subsection{The K-homological calculations}

Write ${G}=SU(2)$ and $\overline{G}=SO(3)$ as before.
In this section we compute the equivariant $K$-homology
of $SU(2)$ on $SU(3)$, acting by conjugation after first projecting onto 
$SO(3)\,$ . The appropriate bundle is that of $SU(3)$ 
on $SU(3)$ at
level 1, with the group $G$ restricted to $SO(3)\subset SU(3)$ -- we let
$\tau\in H^1_G(SU(3);\bbZ_2)\oplus H^3_G(SU(3);\bbZ)$ denote the appropriate
restricted twist.
This bundle is crucial in identifying the relevant twisting
unitaries at various spots in the calculation, as we shall see.
Using spectral sequences (see section 1.2), we quickly compute that 
$H^1_G(X;\bbZ_2)=H^1(X;\bbZ_2)$ for any connected $X$ (indeed, $E^{1,0}_2=0$ 
because $G$ is connected). Thus $H^1_G(SU(3);\bbZ_2)=0$, which means there
is no global $H^1$-twist. Likewise, the spectral sequence says 
$H^3_G(SU(3);\bbZ)={\rm ker}\,{d_4}$ is $\bbZ$ or 0, and so it must be $\bbZ$,
since from the $Sp(4)$ on $Sp(4)$ bundle of section 2.2 it is clear 
$H^3_G(SU(3))$ contains $\bbZ$.
Thus $\tau\in\bbZ$. 

It is useful that all $K^G_*$-groups carry an $R_G$-module
structure. The representation rings we need are identified at the beginning
of subsection 1.2 -- recall Figure 1.

  In section 5.1 we decomposed $SU(3)$ into six spaces of $G$-orbits: $\cO_{SU2},
\ldots,\cO_{\bbA 1}$,
  where the spaces $\cO_{SU2}$, $\cO_{SU2}\cup\cO_{O2}$, $\ldots$,
  $\cO_{SU2}\cup\cO_{O2}\cup\cdots\cup\cO_{\bbA 1}=SU(3)$
  are each compact. We will use the six-term exact sequence (\ref{six-term})
a number of times, to recursively  build up $K_*^G(SU(3))$.

  \bigskip\centerline{{\bf Step 1:} ${}^\tau K_*^G(\cO_{SU2}\cup\cO_{O2}\cup
\cO_\bbT)$}
\medskip

  Orbit $\cO_{SU2}$ consists of the three fixed points, so
$K_*^G(\cO_{SU2})=3\times K_*^G(pt)=
3\times R_{SU2}\,,\, 0$.
Since the global $H^1$-twist is trivial,
the $H^1$-twist  on $\cO_{O2}=3\times \bbR\times G/O(2)$   will also be $+$.
We obtain
$$^+K_j^G(\cO_{O2})=3\times {}^+K_{j+1}^G(G/O(2))=3\times{}^-K_G^{j+1}(G/O(2))
=3\times {}^-K_{O2}^{j+1}(pt)\,,$$
where the grading $-\in H^1_{O2}(pt)$ arises from Poincar\'e duality and the 
nonorientability of the projective plane  $P^2=G/O(2)$. Hence
${}^+ K_j^{{G}}(\cO_{O2})=3\times {}^-R^1_{O2}\,,\,
3\times {}^-R_{O2}$.
These graded representation rings are explicitly described in section 2.1.

We learned in section 5.1 that the orbits with infinite stabilisers together
form a `closed M\"obius strip' $\cO_{SU2}\cup\cO_{O2}\cup\cO_\bbT$, which we'll
call $M$. It's drawn in Figure 8, and consists of three copies of the $SU(3)$ 
Stiefel diagram. The dotted lines in the figure are the `cuts', i.e. the 
overlaps of the open cover of the bundle over $M$ constructed in section 2.1
(more precisely, this bundle is that of $Sp(4)$ on $Sp(4)$, where the group
is restricted to $G$, and the space is then restricted to $M$).

The first step computes the $K$-homology of $M$ by first removing two closed
`intervals' $\overline{ab}$ and $\overline{cd}$. Here, $a,b,c$ are the three
 fixed points, while the `point' $d$ is a copy of the projective plane, with
stabiliser $O(2)$. See Figure 8.

The `interval' $\overline{ab}=1$-$S^2$-1 is the Stiefel diagram for $G$ on $G$
(see Figure 4), so 
$^\tau K_*^G(\overline{ab})$
is $^2 K_*^G(G)=Ver_0(G)$, where the shifted level $0+2$ is obtained from
the bundle, which gives the twisting unitaries attached to each cut. Thus 
$^\tau K_0^G(\overline{ab})=R_G/(\sigma)\simeq\bbZ$ and $^\tau K_1^G(\overline{ab})=0\,$.
Likewise, the `interval' $\overline{cd}=1$-$S^2$-$P^2$ is the Stiefel diagram
for $G$ on $SO(3)$ (again see Figure 4), and  $^\tau K_*^G
(\overline{cd})\simeq{}^{(+,2)} K_*^G(SO(3))\,$ (the first component of the 
(adjoint-shifted) twist $(+,2)$ is the component in 
$H_G^1(SO(3);\bbZ_2)\simeq\bbZ_2$). This $K$-homology was computed in section
4.5 of \cite{BS} to be $R_G/(\sigma)$ and $\bbZ 1^-$ (the graded representation
arises here through the application of Poincar\'e duality and the nonorientability
of the projective plane $P^2=G/O(2)$).
Both $^\tau K^G_*(\overline{ab})$ and $^\tau K^G_*(\overline{cd})$ can also
be easily constructed from first principles (and the bundle), by using the
six-term sequence and removing the endpoints as in (\ref{su2su2}). Indeed, we 
also need the $K$-homology
of the `closed interval' $P^2$-$S^2$-$P^2\,$: this is readily found to be 
$R_G/(\sigma)$ and 
$\bbZ^2 1^-\,$. Finally, $M\setminus (\overline{ab}\cup\overline{cd})$
falls into two copies of $\bbR\times(P^2$-$S^2$-$P^2)$, and the bundle is 
trivial in the $\bbR$-direction, so
$^\tau K^G_*(M\setminus(\overline{ab}\cup\overline{cd}))\simeq \bbZ^4 1^-,
\,2\times R_G/(\sigma)$.

We can compute ${}^\tau K_*^G(M\setminus\overline{ab})$ using the  six-term exact sequence:
\begin{equation}
\begin{array}{ccccc}
\bbZ^4 1^- & {\longleftarrow} & ^{\tau}K^{{G}}_{0}(M\setminus\overline{ab}) &
\stackrel{}{\longleftarrow} & R_{G}/(\sigma)
\\[3pt]
\alpha\downarrow & & & & \uparrow \beta
\\[2pt]
\bbZ 1^- & {\lori} &^{\tau}K^{{G}}_1(M\setminus\overline{ab}) &\lori &2\times
R_G/(\sigma)
\end{array}
\end{equation}
By considering the nonequivariant diagram, we see that the map $\alpha$ can be 
written  $(m_1,m_1';m_2,m_2')1^-\mapsto (m_1+m_1'+
m_2+m_2')1^-$, while $\beta$ sends $(n_1;n_2)[1]\mapsto (n_1+n_2)[1]\,$.
We obtain
 \begin{eqnarray} \label{step1}
 ^{\tau}K^{{G}}_{0}(M\setminus\overline{ab})&=&\,\bbZ^3 1^- \,,\\
 ^{\tau}K^{{G}}_{1}(M\setminus\overline{ab})&=&\,R_G/(\sigma)\,.
 \end{eqnarray}

Now glue in $\overline{ab}$:
\begin{equation}
\begin{array}{ccccc}
\bbZ^3 1^- & {\longleftarrow} & ^{\tau}K^{{G}}_{0}(M) &
\stackrel{}{\longleftarrow} & R_{G}/(\sigma)
\\[3pt]
\downarrow & & & & \uparrow \gamma
\\[2pt]
0 & {\lori} &^{\tau}K^{{G}}_1(M) &\stackrel{}{\lori} &R_G/(\sigma)
\end{array}
\end{equation}
To see what $\gamma$ is, consider the nonequivariant version of this
calculation: the closed M\"obius strip is homotopic to $S^1$, so its
nonequivariant $K$-homology $K_*=K_{*,cs}$ will be that of $S^1$, namely $\bbZ,\,
\bbZ$. Computing this from the six-term by removing a closed interval (homotopic
to $pt$) and thus opening up the strip, we likewise get a vertical map $\gamma':\bbZ\rightarrow\bbZ$
which must vanish in order to recover these $K_*$-groups. The equivariant $K$-homology
here is likewise $1+1$-dimensional, so $\gamma$
should likewise vanish. Hence we obtain:
\begin{eqnarray}
    ^{\tau}K^{{G}}_{0}(\cO_{G}\cup\cO_{O2}\cup\cO_{\bbT})&=&
R_G/(\sigma)\oplus \bbZ^3 1^-\simeq\bbZ^4\,,\\
^{\tau}K^{{G}}_{1}(\cO_{G}\cup\cO_{O2}\cup\cO_{\bbT})&=& R_G/(\sigma)
\simeq\bbZ\,.
\end{eqnarray}
As an $R_G$-module, $^\tau K_0^G(M)$ is a semi-direct sum, where  
$R_G/(\sigma)$ is the submodule, while $\bbZ^3 1^-$ 
is a homomorphic image.

\bigskip\centerline{{\bf Step 2:} $^\tau K_*^G(\cO_{SU2}\cup\cO_{O2}\cup\cO_{\bbT}
\cup\cO_{\bbD 4})$} \medskip

The $K$-homology of $\cO_{\bbD 4}$ is immediate from (\ref{cancel}) (in section
1.2 we explain why the application of Poincar\'e duality needed to relate
$K^G_*(G/\Gamma)$ to the representation ring of the finite subgroup $\Gamma<SU(2)$
 doesn't introduce any $H^1_G$- or $H^3_G$-twists):
\begin{eqnarray}
    K_0^{{G}}(\cO_{\bbD 4})&=&0\,,\\
    K_1^{{G}}(\cO_{\bbD 4})&=&R_{\bbD 4}\simeq\bbZ^5\,.
        \end{eqnarray}
The six-term exact sequence becomes
\begin{equation}\label{d4psi}
\begin{array}{ccccc}
0& {\longleftarrow} & ^{\tau}K^{{G}}_{0}(M\cup
\cO_{\bbD 4}) &{\longleftarrow} & \bbZ^3 1^- \oplus R_G/(\sigma)\\[3pt]
\downarrow & & & & \uparrow\psi 
\\[2pt]
R_G/(\sigma) & {\lori} &
^{\tau}K^{{G}}_1(M\cup\cO_{\bbD 4})&{\lori}&R_{\bbD 4}
\end{array}
\end{equation}
First let's identify the composition of $\psi$ with the projection from
$^\tau K_0(M)$ to $\bbZ^3 1^-$. We learned in section 5.1 that the $\bbR^2$ of
$\cO_{\bbD 4}$ is naturally a triangle, whose three edges can be identified
with the three $\bbR$s of $\cO_{O2}$; heading from the interior to each of
those edges gives three embeddings of $\bbD_4$ in $O(2)$, identifying $\delta$
in turn with each  $s_i$, $i\ne 0$. In particular,
each of these edges can be identified with the value of Res$_{\bbD 4}^{O2}
\,\psi$ there, which will be one of the nontrivial representations $s_i$. 
The composition of $\psi$ with the projection should be three
copies of Dirac 
induction (\ref{dindfinite}) from $R_{\bbD 4}$ to $^- R^1_{O2}$, each with
the different embedding. In particular, the kernel of this composition will be
Span$\{1+s_1+s_2+s_3,t\}$. By the $R_G$-module property of $\psi$, the
image of $\psi$ must intersect $R_G/(\sigma)$ trivially, and so that
also equals the kernel of $\psi$ itself. We obtain:
\begin{eqnarray}\label{step3}
    ^{\tau}K^{{G}}_{0}(\cO_{G}\cup\cdots\cup\cO_{\bbD 4})&=&
    R_G/(\sigma)\simeq\bbZ\,,\\
 ^{\tau}K^{{G}}_{1}(\cO_{G}\cup\cdots\cup\cO_{\bbD 4})&=&
R_G/(\sigma)\oplus \mathrm{Span}\{1+s_1+s_2+s_3,t\}\simeq\bbZ^{3}\,.
\end{eqnarray}
Here, $R_G/(\sigma)$ is a submodule of $^\tau K_1^G$.

\bigskip\centerline{{\bf Step 3:} $^{\tau}K_*^G(\cO_{SU2}\cup\cO_{O2}\cup\cO_{\bbT}
\cup\cO_{\bbD 4}\cup\cO_{\bbA 3})$} \medskip

Recall $\cO_{\bbA 3}=\bbR\times M\ddot{o}b\times G/\bbA_3$. 
We compute $K_*^G(\cO_{\bbA 3})$ by Mayer-Vietoris:
\begin{equation}\label{step4a}
\begin{array}{ccccc}
^{\tau}K^{G}_0(\cO_{\bbA 3}) & {\longrightarrow} & 2\times K^{G}_{0}(\bbR^2\times G/\bbA_3)&
{\longrightarrow} &
2\times K^{G}_{0}(\bbR^2\times G/\bbA_3)
\\[3pt]
\uparrow & & & & \downarrow
\\[2pt]
2\times K^{G}_{1}(\bbR^2\times G/\bbA_3) & \stackrel{\epsilon \,}{\longleftarrow} &
2\times K^{G}_{1}(\bbR^2\times G/\bbA_3)
&{\longleftarrow}
&^{\tau}K^{G}_1(\cO_{\bbA 3})           
\end{array}
\end{equation}
By the usual arguments $K^{G}_{1}(\bbR^2\times G/\bbA_3)=R_{\bbA 3}$ but
$K^{G}_{0}(\bbR^2\times G/\bbA_3)=0$. We take $\epsilon(f,g)=(f+g,f-g)$ (the
sign   arising because of the M\"obius).
For this map, (\ref{step4a}) says $K^G_*(\cO_{\bbA 3})=\bbZ_2\otimes R_{\bbA 3},
\,0$ and the six-term sequence immediately gives:
\begin{eqnarray}
    ^{\tau}K^{{G}}_{0}(\cO_{G}\cup\cdots\cup\cO_{\bbA 3})&=&
 (\bbZ_2\otimes R_{\bbA 3}) \oplus R_G/(\sigma)\simeq\bbZ_2^4\oplus\bbZ\,,
\label{tors}\\
 ^{\tau}K^{{G}}_{1}(\cO_{G}\cup\cdots\cup\cO_{\bbA 3})&=&
R_G/(\sigma)\oplus\mathrm{Span}\{1+s_1+s_2+s_3,t\}\simeq\bbZ^{3}\,,
\end{eqnarray}
with the first summand of both these $K$-homology groups being the submodule. Note
the torsion in $^\tau K_0^G$.

    \bigskip\centerline{{\bf Step 4: $^{\tau}K_*^G(SU(3))$} \medskip }

We get by the usual arguments that 
$K_*^{{G}}(\cO_{\bbA 1})=K_{*}^{\bbA 1}(pt)$ is $R_{\bbA 1},\,0$ for $*=0,1$
respectively. The six-term exact sequence becomes:
\begin{equation}
\begin{array}{ccccc}
R_{\bbA 1}& {\longleftarrow} & ^{\tau}K^{{G}}_{0}(SU(3)) &
{\longleftarrow} & \bbZ_2^4\oplus \bbZ
\\[3pt]
\phi\downarrow & & & & \uparrow
\\[2pt]
R_G/(\sigma)\oplus \mathrm{Span}\{1+s_1+s_2+s_3,t\} & {\lori} &
^{\tau}K^{{G}}_1(SU(3)) &{\lori} &0
\end{array}
\end{equation}
As usual, compose $\phi$ with the projection from $^\tau K_1^G(\cO_G\cup\cdots
\cup \cO_{\bbA 3})$ to the submodule of $R_{\bbD 4}$. Now, $\cO_{\bbD 4}=
\bbR^2\times G/\bbD_4$ consists of the boundary points of $\cO_{\bbA 1}=
\bbR^5\times G/\bbA_1$, as we move along two of the five $\bbR$s. So
that composition of $\phi$ with the projection should be given by
the induction from $R_{\bbA 1}$ to $R_{\bbD 4}$. Indeed, it is well-defined,
 sending $r_1''$ to
$\sum_is_i$ and $r''_{-1}$ to $t$. By the $R_G$-module property $\phi$
cannot see $R_G/(\sigma)$. The final answer is then:
\begin{eqnarray} \label{d4answer}
    ^{\tau}K^{{G}}_{0}(SU(3))&=&(\bbZ_2\otimes R_{\bbA 3})
 \oplus R_G/(\sigma)\simeq\bbZ_2^4\oplus\bbZ\,,\\
 ^{\tau}K^{{G}}_{1}(SU(3))&=& R_G/(\sigma)\simeq\bbZ\,.
\end{eqnarray}
We discuss the meaning of this $K$-homology, and in particular its
relation to the full system of the
$D_4$ modular invariant, in  the concluding  section.
The torsion in $^\tau K_0^G(SU(3))$ is mysterious though.

\section{The $E_6$ modular invariant of $SU(2)$}

Write ${G}=SU(2)$ as before. 
The `$E_6$' exceptional modular invariant of $SU(2)$ arises from
the conformal embedding of $SU(2)$ at level 10, into
$Sp(4)$ level 1. This conformal embedding belongs to an infinite series of
$Spin(n)$ level 5 into $Spin((n-1)(n+2)/2)$ at level 1 \cite{BaBo,ScWa}, where the
embedding $Spin(n)$ into $Spin((n-1)(n+2)/2)$ is given by the representation
with Dynkin labels $(2,0,0,\ldots,0)$. For us, $n=3$ and we identify $Spin(3)$
with $SU(2)$ and $Spin(5)$ with $Sp(4)$; the
doubling of the level from 5 to 10 comes from the identification of
Lie algebras $so(3)$ and $su(2)$. This doubling can be confirmed by a
conformal charge calculation, or indeed the calculation in section 2.3 above.

We will identify the symplectic group $Sp(4)$ with the set of all $4\times 4$
unitary matrices $\mathbf{B}$ commuting with $J=\bigl({0\atop -I_2}{I_2\atop 0}\bigr)$:
this commutation with $J$ is equivalent to the block form 
\begin{equation}\label{(2)}
\mathbf{B}=\left(\begin{matrix}B&C\cr D& E\end{matrix}\right)\,,
\end{equation}
 where $B^tD=D^tB$, $E^tC=C^tE$, $I=B^tE-D^tC$.

In a few spots in the orbit argument, it is very convenient to have an
explicit description of this embedding $G\hookrightarrow Sp(4)$,
 which we will call ${\cal R}^{(4)}$:
\begin{equation}\label{su2sp4}
{\cal R}^{(4)}(\gamma,\delta):=
\left(\begin{matrix}
\gamma^3&\sqrt{3}\gamma^2\delta&\delta^3&-\sqrt{3}\gamma\delta^2\cr
-\sqrt{3}\gamma^2\overline{\delta}&(3|\gamma|^2-2)\gamma&\sqrt{3}
\overline{\gamma}{\delta}^2 &(1-3|\gamma|^2)\delta\cr
-\overline{\delta}^3&\sqrt{3}\overline{\gamma}\overline{\delta}^2&\overline{\gamma}^3
&\sqrt{3}\overline{\gamma}^2\overline{\delta}\cr
-\sqrt{3}\gamma\overline{\delta}^2&(3|\gamma|^2-1)\overline{\delta}&-\sqrt{3}
\overline{\gamma}^2\delta&(3|\gamma|^2-2)\overline{\gamma}
\end{matrix}\right)\end{equation}
(recall (\ref{su2param})).
This can be written in block form $\mathbf{A}=\bigl({A\atop-\overline{A}'}
{A'\atop \overline{A}}\bigr)$. Useful special cases of (\ref{su2sp4})
are ${\cal R}^{(4)}(e^{\i t},0)={\rm diag}(e^{3\i t},e^{\i t},e^{-3\i
  t},e^{-\i t})$ and $\mathcal{R}^{(4)}(0,1)=\bigl({0\atop -I}{I\atop
  0}\bigr)$. Write $T={\cal R}^{(4)}(*,0),T'={\cal R}^{(4)}(0,*)$ for the
two circles coming from the image of $O(2)\subset G$. It will be useful in
section 6.2 to note that the normaliser of $T$ in ${\cal R}^{(4)}(G)$ is
$T\cup T'\simeq O(2)$.

\subsection{The K-homology groups}\bigskip

The $K$-homology here can be elegantly computed by the spectral sequence
methods of section 3.3.

The representation rings $R_{SU2}$ and $R_{Sp4}$ are identified with
$\bbZ[\sigma]$, and $\bbZ[s,v]$ respectively. Here if $\sigma_d$ is the $d$-dimensional
representation  of $\SUz$, then $\sigma = \sigma_2$, and the restriction of the spinor
representation $s$ from $Sp(4)$ to  $\SUz$ is $\sigma_4 = \sigma^3 - 2\sigma$,
and the restriction of the vector representation $v$ is $\sigma_5 = \sigma^4 - 3\sigma^2 +1$\,.
Restriction makes $R_{SU2}$ into an $R_{Sp4}$-module.
At level $1$, the fusion rules of $Sp(4)$ coincide with those of the Ising 
model, and are described by the fusion ideal $I_1$
generated by $s^2 - v - 1, v^2 - 1, vs - s $, however $v^2 -1 $ is redundant as
$v^2 - 1 = (1-v)(s^2 - v - 1) + s(vs - s)$\,. If $G = Sp(4)$, we have a free 
resolution of the Verlinde algebra at level $1$ by
$$ 0 \rightarrow R_G \stackrel{f}{\rightarrow} R^2_G \stackrel{g}{\rightarrow}
 R_G\rightarrow R_G/I_1\rightarrow 0\, ,$$
where $f(a) = ((ns - s)a, (s^2 - n - 1)b)$ and
$g(a,b) = (s^2 - n - 1)a + (ns - s)b$. This resolution has the same length as that
of Meinrenken \cite{M} but smaller degree.
To compute Tor$_i^{R_G}(R_{SU2},Ver_1(Sp(4)))$, we ignore the last term in the free resolution and tensor with $R_{SU2}$ to get a complex
$$ 0 \rightarrow ( R_G\otimes_{R_G} R_{SU2} \simeq R_{SU2}) \stackrel{\partial_2}{\rightarrow} (R^2_G  \otimes_{R_G} 
R_{SU2} \simeq R_{SU2}^2) \stackrel{\partial_1}{\rightarrow}(R_G\otimes_{R_G} R_{SU2}
\simeq R_{SU2})\rightarrow 0\,.$$
Here, if $p, q \in \bbZ[\sigma]$, then 
\begin{eqnarray}
\partial_1(p,q)&=& (s^2 - v - 1)p + (vs - s)q = (\sigma_4^2 - \sigma_5 - 1)p + 
(\sigma_5 -1)\sigma_4q\nonumber\\
&=& (\sigma^2 -2 )[(\sigma^4 - 3\sigma^2 + 1)p + \sigma^3(\sigma^2 -3)q]\,.
\end{eqnarray}
  Hence the image of $\partial_1$ is $(\sigma^2-2)\bbZ[\sigma]$,
and Tor$_0^{R_G}(R_{SU2},Ver_1(Sp(4)))=H_0 = {\rm ker}(\partial_0)/{\rm im}(\partial_1) = 
\bbZ[\sigma]/(\sigma^2-2)$ is two-dimensional.
Similarly 
\begin{equation}
\partial_2(p) = ( (vs - s)p, (s^2 - v - 1)p)  = (\sigma^2-2)( 
\sigma^3(\sigma^2 -3)p, -(\sigma^4 - 3\sigma^2 + 1)p)\,,
\end{equation}
 and so im$(\partial_2)
= ( \sigma^3(\sigma^2 -3), -(\sigma^4 - 3\sigma^2 + 1))(\sigma^2-2)\bbZ[\sigma] 
\subset {\rm ker}(\partial_1) = ( \sigma^3(\sigma^2 -3), -(\sigma^4 - 3\sigma^2 + 
1))(\sigma^2-2)\bbZ[\sigma]$
and $H_1 = {\rm ker}(\partial_1)/{\rm im}(\partial_2) = \bbZ[\sigma]/(\sigma^2-2)$ is again
two-dimensional.

Thus in the Hodgkin spectral sequence (section 3.3), $E^2_{p,q\,even}$ is 
two-dimensional for $p=0,1$ while all other $E^2_{p,q}$ vanish. Therefore for
$r>1$, all maps $d_r:E^r_{p,q}\rightarrow E^r_{p-r,q+r-1}$ in this spectral 
sequence will be trivial (having either trivial domain or range), so
$E^\infty_{p,q}=E^2_{p,q}$ and
\begin{equation}\label{e6answer}
^\tau K_{SU2}^p(Sp(4))\simeq\bbZ[\sigma]/(\sigma^2-2)
\end{equation}
is two-dimensional for any $p$.
By comparison, $^\tau K_{Sp4}^*(Sp(4))$ is $\bbZ^3,\,0$ (recovering the
Ising fusions); in (\ref{e6answer}) the spinor $s\in Ver_1(Sp(4))$ is sent
to 0 and the vector $v\in Ver_1(Sp(4))$ goes to 1.

The full system of the $SU(2)_{10}$ `$E_6$' modular invariant is 12-dimensional,
built out of two copies of the (unextended) $E_6$ Dynkin diagram, as in
Figure 3.
We don't know yet how to reconcile  this with (\ref{e6answer}) (see the concluding
section for some thoughts in this direction), but
based on similar calculations earlier in this paper, we may hope that the
full system arises by having $SU(2)$ act instead on some closed submanifold
of $Sp(4)$. For this purpose, we now proceed to work out the $SU(2)$-orbits
in $Sp(4)$ and recompute (\ref{e6answer}) the long way.

\subsection{The orbit analysis}\bigskip

We want the orbits of the conjugate action of ${G}$ on $Sp(4)$,
using the embedding ${\cal R}^{(4)}$. 

We'll be finding the orbits in inverse order of the size of their
stabilisers, by taking an element of maximal order and diagonalising it. 
The simplest way to verify that we're not counting some orbit
twice, i.e.\ that what is written for the stabiliser is the full stabiliser and
not merely a subgroup of it, seems to be to diagonalise the different generators
of the stabiliser, and confirm 
visually that
the resulting expression for $\mathbf{B}$ doesn't fall into a different orbit.

First of all, recall that
any element $\widetilde{{\mathbf A}}\in{G}$ lies in a maximal torus, i.e.\
there is some matrix $\widetilde{P}\in{G}$ such that
$\widetilde{P}{}^{-1}\widetilde{{\mathbf A}}\widetilde{P}=\left(\begin{matrix}e^{\i\theta}&0\cr
0&e^{-\i\theta}\end{matrix}\right)$. It is clear from this that if
$\widetilde{{\mathbf A}}$ has order $n$, 
then so will $\mathbf{A}={\cal R}^{(4)}(\widetilde{\mathbf{A}})$.

\smallskip\noindent{\bf Lemma 2.} {\it If $\mathbf{B}\in Sp(4)$ has finite stabiliser,
then any
$\mathbf{A}\in {\cal R}^{(4)}({G})$ in its stabiliser has order $\le 4$ or 6.}
\smallskip

\noindent{\it Proof:} Suppose $\mathbf{A}$ has finite order $n> 6$ or $n=5$. 
Diagonalising it, without loss of generality
we can consider $\mathbf{A}={\rm diag}(\xi_n^3,\xi_n,\xi_n^{-3},\xi_n^{-1})$ 
where $\xi_n=\exp[2\pi\i/n]$. Commuting with $\mathbf{B}$ in (\ref{(2)}), we 
get that $B$ and $E$ are both diagonal and $C=D=0$. But such a $\mathbf{B}$
commutes with the full maximal torus $T$, and so would have infinite 
stabiliser. {QED}\smallskip

Note that ${\cal R}^{(4)}(-1,0)=-I$ is always in the stabiliser.
The only finite subgroups of $SU(2)$ of even order, whose elements all have 
order $\le 4$ or 6, are the cyclic groups $\bbA_1=C_2,\bbA_3=C_4$, and 
$\bbA_5=C_6$, the binary dihedral groups $\bbD_4=B(C_2\times C_2)=Q_4$ and
$\bbD_5=BS_3$, as well as the tetrahedral group $\bbE_6=BA_4$.
The character tables
and other information about these groups are given in eg.\ \cite{IN}.

\medskip\noindent{\bf Case 1}: {\it Orbits with infinite stabiliser}\medskip

The only orbits with infinite stabiliser are
\begin{eqnarray}
\cO_G&=&\pm I\,,\\
\cO_{O2}&=&2\times {G}/{O(2)}\,,\\
\cO_{\bbT}&=&(S^2-4)\times {G}/\bbT\,.\label{tetra}
\end{eqnarray}
The four points in $\cO_G\cup\cO_{O2}$ are the punctures of $S^2$,
so the union $\cO_G\cup\cO_{O2}\cup\cO_{\bbT}$ is compact, as is both $\cO_G$ 
and $\cO_{O2}$. 
So we can take care of the issue of infinite-dimensionality, with a single
six-term exact sequence.

The reason only $\pm I\in Sp(4)$ have stabiliser ${G}$, 
is Schur's Lemma: ${\cal R}^{(4)}$ is an irreducible
representation of ${G}$, so the only matrices which can commute with
all ${\cal R}^{(4)}(G)$ are the scalar matrices.

Now suppose the stabiliser contains a maximal torus, 
which without loss of generality we can take to be the `canonical' one
${\bbT}$. Write the block-forms $\mathbf{A}={\cal R}^{(4)}(
\alpha,0)={\rm diag}(A,\overline{A})$ and $\mathbf{B}$ in (\ref{(2)}).
$\mathbf{A}\mathbf{B}=\mathbf{B}\mathbf{A}$ requires $B$ and $E$ to be diagonal,
and $C=D=0$. Since $\mathbf{B}\in Sp(4)$, we also have $E=\overline{B}$.
The set of
such $\mathbf{B}$ form a torus ${T}^2$.
It is elementary to confirm that
such a matrix $\mathbf{B}$ commutes with  some $\mathbf{A}\not\in {T}\cup T'$
iff $\mathbf{B}=\pm I$, and commutes with some $\mathbf{A}\in T'$ iff
$\mathbf{B}={\rm diag}(\pm 1,\pm 1,\pm 1,\pm 1)$, where the first and third,
and second and fourth, signs must be equal for $\mathbf{B}\in Sp(4)$.

The final ingredient is the Weyl group of ${G}$. More precisely, 
conjugating by any reflection in $T'$ will amount to
a nontrivial involution of $\cO_{\bbT}$, sending $\mathbf{B}$ to its complex
conjugate. This simultaneous complex-conjugation of $\mathbb{T}\times
\mathbb{T}$ fixes the four points making up $\cO_G$ and $\cO_{O2}$ (as it must).
We should identify points identified by this involution, so this gives
$(T^2-4)/2$, which naturally folds  to a tetrahedron with its vertices
removed, i.e.\ is homeomorphic to the sphere with four punctures.

\medskip\noindent{\bf Case 2}: {\it All finite stabilisers containing elements of order $6$}\medskip

These orbits are
\begin{eqnarray}
\cO_{\mathbb{A}5}&=&(B^2\setminus 3\mathbb{R})\times G/\mathbb{A}_5\,,\label{A5orb}\\
\cO_{\mathbb{D}5}&=&2\times \mathbb{R}\times{G}/\mathbb{D}_5\,,\\
\cO_{\mathbb{E}6}&=&\mathbb{R}\times {G}/\mathbb{E}_6\,,
\end{eqnarray}
where $B^2$ denotes an open solid ball (whose boundary is the tetrahedron of $\cO_{\mathbb{T}}$),
and the 3 $\mathbb{R}$'s of $\cO_{\mathbb{A}5}$ are the $2+1$ $\bbR$'s of
$\cO_{\mathbb{D}5}$ and $\cO_{\mathbb{E}6}$. The 
$\mathbb{R}$'s of $\mathcal{O}_{\mathbb{D}5}$ are chords with endpoints
$\pm I\in\mathcal{O}_G$ and $\mathrm{diag}(\mp 1,\pm 1,\mp 1,\pm 1)\times G/
O(2)\in\mathcal{O}_{O2}$. The $\bbR$ of $\cO_{\mathbb{E}6}$ is a chord with
endpoints at the orbits $\pm I$ in $\cO_G\,$. Of course $\mathbb{E}_6$ is the 
symmetry group of
the tetrahedron, and we find that it is the largest finite stabiliser for this 
action of $G$ on $Sp(4)$.

Let $\mu$ be an order 6 element in the stabiliser of $\mathbf{B}\in Sp(4)$.
Without loss of generality diagonalise $\mu$, so $\mathbf{A}=
{\cal R}^{(4)}(\xi_6,0)={\mathrm{diag}}(-1,\xi_6,-1,\overline{\xi_6})$.
Then $\mathbf{B}$ commuting with $\mathbf{A}$ requires
\begin{equation}\label{order6}
\mathbf{B}=\left(\begin{matrix}b&&c&\cr&b'&&\cr d&&e&\cr &&&e'\end{matrix}\right)
\end{equation}
(we avoid writing the 0's), and $\mathbf{B}\in Sp(4)$ then forces
$b'e'=1\,$, $|b'|=1\,$, $e=\overline{b}\,$, $d=-\overline{c}\,$, and $1=|b|^2+|c|^2$.
Because $\mathbf{B}$ has finite stabiliser, $c\ne 0$. 

The finite stabilisers containing an order 6 element are $\mathbb{A}_5\,$,
$\mathbb{D}_5$ and $\mathbb{E}_6\,$. Which of the $\mathbf{B}$ in (\ref{order6})
have stabiliser $\mathbb{D}_5$? All subgroups of $G$ isomorphic to 
$\mathbb{D}_5$ are conjugate in $G$ to $\langle \mathrm{diag}
(\xi_6,\overline{\xi_6}),\bigl({0\atop-1}{1\atop 0}\bigr)
\rangle$ (this fails if $G$ is replaced with $U(2)$). The matrices $\mathbf{B}$ of (\ref{order6})
which commute with $\mathcal{R}^{(4)}(0,1)$ have $c,b,b'\in\bbR\,$, so $b=e$
and $b'=e'=\pm 1$. Conjugating everything by $\mathcal{R}^{(4)}(\xi_{12},0)$
(which normalises this $\mathbb{D}_5$), we see we can take $c>0$. 
Thus the orbits with $\mathbb{D}_5$ stabiliser form two circular arcs:
$(b,b',c)=(\cos\,\theta,\pm 1,\sin\,\theta)$ for $0<\theta<\pi$.

Likewise, all subgroups of $G$ isomorphic to $\mathbb{E}_6$ are 
conjugate in $G$ to $\langle \mathrm{diag}(\xi_6,\overline{\xi_6}),\tau\rangle$
for $\tau=\frac{1}{\sqrt{3}}\bigl({-\i\atop (-1-\i)\overline{\xi_{24}}}{(1+\i)\xi_{24}
\atop\i}\bigr)$ (again this fails for $U(2)$). The matrices $\mathbf{B}$
of (\ref{order6}) which commute with $\tau$ have $c=-2\sqrt{2}y$, 
$b'=x-3y\i$ where $b=x+\i y$. Conjugating by $\mathcal{R}^{(4)}(0,\xi_{12})$
(which normalises $\mathbb{E}_6$) again shows we can restrict to $c>0$.
We thus get one elliptical arc, where $1=x^2+9y^2$.

The complement of these three chords will have stabiliser $\mathbb{A}_5$.
Conjugating any such $\mathbf{B}$
by $\mathbf{T}\subset G$, we see that again we may take  $c$ real and positive.
The parameter space so far is a solid
torus -- the interior of the torus of case 1. As was the case there,
there is a final folding that can be done: by eg. $\mathcal{R}^{(4)}(0,1)$.
This replaces the torus with the tetrahedron as before.

\medskip\noindent{\bf Case 3}: {\it All finite stabilisers containing elements of order $4$ but not higher}\medskip

These orbits are
\begin{eqnarray}
\cO_{\mathbb{A}3}&=&(B^2\setminus 4\mathbb{R})\times G/\mathbb{A}_3\,,\label{orba3}\\
\cO_{\mathbb{D}4}&=&\mathbb{R}\times{G}/\mathbb{D}_4\,.
\end{eqnarray}
This open ball $B^2$ also has boundary the tetrahedron $S^2$ of {\bf Case 1} (though
of course it is disjoint from the open ball of {\bf Case 2}). 
The endpoints of the chord $\mathbb{R}$ of $\cO_{\mathbb{D}4}$ lie in distinct
$\cO_{O2}$-orbits. The
four $\bbR$'s in $\cO_{\mathbb{A}3}$ come from overlaps with $\cO_{\mathbb{D}4}$,
$\cO_{\mathbb{D}5}$ and $\cO_{\mathbb{E}6}$, and form a square with vertices
on the boundary $S^2$. 

Without loss of generality we can diagonalise an order 4 element in the
stabiliser. We find that ${\mathbf B}$ commutes with ${\mathcal R}^{(4)}(\i,0)$
iff it is of the form
\begin{equation}\label{pre4}
{\mathbf B}=\left(\begin{matrix}b&&&c\cr &b'&c'&\cr &d&e&\cr d'&&&e'\end{matrix}\right)\,.
\end{equation}
Such a matrix ${\mathbf B}$ lies in $Sp(4)$ iff it is in fact of the form
\begin{equation}\label{stab4}
{\mathbf B}=\left(\begin{matrix}b&&&c\cr &\overline{c}b/\overline{c}'&c'&\cr 
&-\overline{c}'&\overline{b}&\cr -\overline{c}&&&c\overline{b}/c'\end{matrix}\right)\,,
\end{equation}
where $|c|=|c'|\ne 0$ (they must be nonzero, otherwise the stabiliser will
be infinite) and $|b|^2+|c|^2=1$.

The finite stabilisers containing an order-4 element are $\mathbb{D}_5$ and
$\mathbb{E}_6$ (both already dealt with) and $\mathbb{A}_3$ and $\mathbb{D}_4$.
Consider first $\mathbb{D}_4$. Any such subgroup of $G$ will be conjugate to
$\langle \bigl({\i\atop 0}{0\atop-\i}\bigr),\bigl({0\atop -1}{1\atop 0}\bigr)
\rangle$. Requiring ${\mathbf B}$ in (\ref{pre4}) to also commute with
${\mathcal R}^{(4)}(0,1)$ means
\begin{equation}
{\mathbf B}=\left(\begin{matrix}b&&&c\cr &b'&c'&\cr &-c&b&\cr -c'&&&b'\end{matrix}\right)
\end{equation}
and imposing (\ref{stab4}) then tells us $b,c$ are both real, $b^2+c^2=1$,
$c'=\pm c$ and $b'=\pm b$ (same sign). But we could have diagonalised any of
the other order-4 elements in the stabiliser; conjugating by the order-6
element $\mu=\frac{1}{\sqrt{2}}
\bigl({\overline{\xi}\atop -\xi}{\overline{\xi}\atop \xi}\bigr)\in G$,
or by the order-4 element $\nu=\bigl({0\atop -\overline{\xi}}{\xi\atop 0}\bigr)$ (where 
$\xi=\xi_8$) cyclically permutes these order-4 elements, so we should identify
$\mathbf{B}$ with its conjugates by ${\mathcal R}^{(4)}(\mu)$ and ${\mathcal R}^{(4)}
(\nu)$. Doing this,
we see that $c=+c'$ actually commutes with ${\mathcal R}^{(4)}(\mu)$, so these
$\mathbf{B}$ actually have stabiliser containing $\langle \mathbb{D}_4,\mu\rangle=\mathbb{E}_6$
and so have been considered already in {\bf Case 1}. Thus it suffices to
consider the other sign, $c=-c'$. The action of $\mu$ on $b,c$ can then be
described by the order-3 rotation $b+\i c\mapsto \xi_3\,(b+\i c)$ on the
unit circle, while that of $\nu$ is complex-conjugation: $b+\i c\mapsto b-\i c$. 
The six points $e^{\ell\pi i/6}$ on the unit circle lie in
the orbits $\cO_{O2}\cup \cO_G$, so should be removed. The result is a single arc worth
of $\mathbb{D}_4$-orbits.

To identify the $\mathbb{A}_3$ orbits, return to (\ref{stab4}). Conjugating by
arbitrary elements of $T$, we see that we can require $c>0$. Hence the $\mathbb{A}_3$
orbits are parametrised by the value of $b$ (in the unit disc) and the value of
Arg$(c')$ (in the circle), i.e. a solid torus. We must 
eliminate the case where $b$ and $c'$ are both real, as their stabiliser
would be at least $\mathbb{D}_4$. We must also eliminate the case with 
stabiliser at least $\mathbb{D}_5$. Finally, we must identify conjugates 
by $\mathcal{R}^{(4)}(0,1)$, which sends $(b,c,c')$ to $(\overline{b},c,\overline{c}')$.
The result is described above.

\subsection{The K-homological calculations}\bigskip

We consider here the $K$-homology calculations of the $E_6$ conformal embedding
for $SU(2)$ at level 10. Recall the representation rings in section 1.2.

First, we need the cohomology groups $H^1_{{G}}(Sp(4);\bbZ_2)$
and $H^3_{{G}}(Sp(4);\bbZ)$. For the former, $E_\infty^{1,0}=
{\rm Hom}(SU(2),\bbZ_2)=0$ since $SU(2)$ is connected. We have
$H^2_{SU2}(pt;\bbZ_2)=0$ so the desired  $H^1_{{G}}(Sp(4);\bbZ_2)$
equals $E_\infty^{0,1}=H^1(Sp(4);\bbZ_2)\simeq{\mathrm{Hom}}(H_1(Sp(4)),\bbZ_2)=0$
(see page 291 of \cite{Ha}). This means that there is no possibility for
a global $H^1$-grading here. 

Computing $H^3_{{G}}(Sp(4);\bbZ)$ by spectral sequences requires
knowing $H^q(Sp(4);\bbZ)=\bbZ,0,0,\bbZ,0$ for $q=0,1,2,3,4$ resp. (see
page 434 of \cite{Ha}), as well as $H^p_{G}(pt;\bbZ_2)=\bbZ_2$ or 0, and
 $H^p_{G}(pt;\bbZ)=\bbZ$ or 0, both depending on
whether or not 4 divides $p$ (see section 1.2).
Then 
$H^3_{{G}}(Sp(4);\bbZ)={\rm ker}\,d_4$ for
$d_4:H^3(Sp(4);\bbZ)\rightarrow H^4_{G}(pt;\bbZ)$.
From the bundle picture, $H^3_G(Sp(4);\bbZ)$ contains at least $\bbZ$ 
(associated to the level of the $Sp(4)$ Verlinde algebra), and so it
must equal $\bbZ$.

\bigskip
\centerline{{\bf Step 1:} {\it The infinite stabilisers}}\medskip

This is where finite-dimensionality is won (or lost) -- assuming the
parameter spaces of the orbit spaces $\cO_G,\cO_{O2},\ldots,\cO_{\bbA 1}$ are
sufficiently nice (which they are).
The tetrahedron $Tet=S^2$ of (\ref{tetra}) is drawn in Figure 9(a);
the vertices are $a,b\in\cO_G$ and $c,d\in\cO_{O2}$. In Figure 9(b) the closed
edges $\overline{ad}$ and $\overline{bc}$ are removed, resulting in an open 
cylinder $cyl$. Four copies (the triangles in (b)) of
the $Sp(4)$-Stiefel diagram tile this tetrahedron, although they
coincide with only half of the faces of the tetrahedron (the triangles in
(a)). The dotted lines in (b) are the cuts (i.e. pairwise intersections) of
the open cover of the tetrahedron coming from the bundle constructed in
section 2.2.

\bigskip 
\begin{figure}[tb]
\begin{center}
\epsfysize=2in \centerline{\epsffile{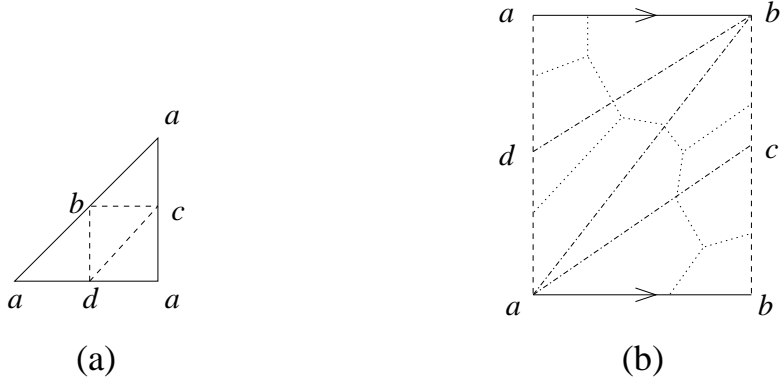}}
\caption{(a) The tetrahedron $Tet$\qquad (b) The cylinder $cyl$}
\end{center}
\end{figure}

We compute the $K$-homology of $Tet$ in two steps: first work out
the cylinder $cyl$, then use the six-term exact sequence to glue in the two edges.
The $K$-homology $^\tau K^G_*(cyl\times G/\bbT)$ collapses to $^2 K_{*+1}^\bbT
(S^1)$, where the twist $2\in H^3_{\bbT}(S^1)$ is determined from Figure
9(b) and the bundle of section 2.2 (which gives the twisting unitary across 
each cut -- the bundle on $cyl=S^1\times\bbR$ splits as a product and is trivial on
$\bbR$). Thus the $K$-homology of $cyl$ (recall section 3.1) is
\begin{eqnarray}
^\tau K_0^G(S^2\setminus(\overline{ad}\cup\overline{bc}))&=&\,0\,,\\
^\tau K_1^G(S^2\setminus(\overline{ad}\cup\overline{bc}))&=&\,\bbZ[a^{\pm 1}]/(1-a^2)\,.
\end{eqnarray}

Note that $^\tau K^G_*(\overline{ad})$ and ${}^\tau K^G_*(\overline{bc})$
can both be identified with $^{\tau'} K^G_*(SO(3))$ for some level $\tau'$
(which is found in the usual way to be 1 for both intervals). This $K$-homology
is computed in section 4.1 of \cite{BS}:
\begin{eqnarray}
^\tau K_0^G(\overline{ad})\simeq {}^\tau K_0^G(\overline{bc})&=&\,0\,,\label{cap0}\\
^\tau K_1^G(\overline{ad})\simeq {}^\tau K_1^G(\overline{bc})&=&\,\bbZ 1^-\,.
\end{eqnarray}
 Here, $1^-\in {}^-R^1_{O2}$ (recall section 2.1); the reason for the $H^1$-twist
is the nonorientability of the projective plane $G/O(2)$.

From this and (\ref{six-term}), the $K$-homology of the tetrahedron $Tet$ 
is immediate:
\begin{eqnarray}
^\tau K_0^G(\cO_G\cup\cO_{O2}\cup\cO_{\bbT})&=&\,0\,,\\
^\tau K_1^G(\cO_G\cup\cO_{O2}\cup\cO_{\bbT})&=&\,\bbZ^2 1^-\oplus 
R_\bbT/(1-a^2)\simeq\bbZ^4\,.\label{K1tet}
\end{eqnarray}
This is really a semi-direct sum of $R_G$-modules; $\bbZ^2 1^-$ is a 
submodule of $^\tau K_1^G$, while $R_\bbT/(1-a^2)$ is a homomorphic image 
(quotient).

\bigskip
\centerline{{\bf Step 2:} {\it The  $\bbA_5$-ball}}\medskip

Recall the `tetrahedron' $Tet=\cO_G\cup \cO_{O2}\cup\cO_{\bbT}$ of
Step 1. The boundary of the genus-3 volume $B^2\setminus 3\bbR$ of (\ref{A5orb})
is $Tet$ together with the three chords $\cO_{\bbD 5}\cup\cO_{\bbE 6}$. First
let's glue those chords to $Tet$ using the six-term sequence:
\begin{equation}\label{e6alpha}
  \begin{matrix} 2\times R_{\bbD 5}\oplus R_{\bbE 6}&\longleftarrow&{}^\tau K^{
{G}}_0(Tet\cup \cO_{\bbD 5}\cup\cO_{\bbE 6})&\longleftarrow
   &0\cr \alpha\downarrow&&&&\uparrow\cr 
  \bbZ^2 1^-\oplus R_\bbT/(1-a^2) 
   &\longrightarrow&{}^\tau K^{{G}}_1(Tet\cup \cO_{\bbD 5}\cup\cO_{\bbE 6})
   &\longrightarrow&   0\end{matrix}
\end{equation}
As explained in section 1.2, the representation rings in the upper-left entry 
are untwisted. We claim $\alpha\equiv 0$. Neither the $\bbD_5$- nor the 
$\bbE_6$-`chords' touch $\cO_\bbT$, so $\alpha$ ignores the $R_\bbT/(1-a^2)$ 
summand. Likewise, $\bbE_6$ can't be related to $1^-\in{}^-R^1_{O2}$. 
The easiest way to see $\alpha$ must vanish is because the
alternative would mess up (\ref{A3betagamma}) below. We'll return to this 
shortly, but assume $\alpha\equiv 0$ for now. We then get:
\begin{eqnarray}
^\tau K_0^G(Tet\cup \cO_{\bbD 5}\cup\cO_{\bbE 6})&=&\,
R^2_{\bbD 5}\oplus R_{\bbE 6}\simeq\bbZ^{19}\,,\\
^\tau K_1^G(Tet\cup \cO_{\bbD 5}\cup\cO_{\bbE 6})&=&\,\bbZ^2 1^-\oplus 
R_\bbT/(1-a^2)\simeq\bbZ^4\,.
\end{eqnarray}
As with (\ref{K1tet}), $^\tau K_1^G$ is a semi-direct sum with $\bbZ^2 1^-$
the submodule.

Note from (\ref{trivact}) that $K^*_G((B^2\setminus n\times\bbR)\times G/H)
\simeq R_H\otimes_\bbZ K^*(T_n)$ where $T_n$ is the open solid torus of genus 
$n$; its nonequivariant $K$-homology is easily found by induction to be 
$K^0(T_n)\simeq\bbZ^n$, $K^1(T_n)\simeq\bbZ$. Thus $K_*^G(\cO_{\bbA 5})\simeq 
R_{\bbA 5}, R_{\bbA 5}^3$ (the dimension-shift comes from factoring off an
implicit $\bbR^3$ before applying Poincar\'e duality).
Writing $B_{\bbA 5}$ for the closed ball $Tet\cup \cO_{\bbD 5}\cup\cO_{\bbE 6}
\cup\cO_{\bbA 5}$, we obtain:
\begin{equation}\label{A5betagamma}
  \begin{matrix} R_{\bbA 5}&\stackrel{\delta\,}{\longleftarrow}&{}^\tau K^{
{G}}_0(B_{\bbA 5})&\longleftarrow &R_{\bbD 5}^2\oplus R_{\bbE 6}
   \cr \beta\downarrow&&&&\uparrow\,\gamma\cr 
  \bbZ^2 1^-\oplus R_\bbT/(1-a^2) 
   &\longrightarrow&{}^\tau K^{{G}}_1(B_{\bbA 5})
   &\longrightarrow&   R_{\bbA 5}^3\end{matrix}
\end{equation}
The map $\gamma$ is clear from the nonequivariant calculation: it will be
the diagonal inductions $(\rho_1,\rho_2,\rho_3)\mapsto(\mathrm{Ind}_{\bbA 5}^{\bbD 5}
\rho_1,\mathrm{Ind}_{\bbA 5}^{\bbD 5}\rho_2,\mathrm{Ind}_{\bbA 5}^{\bbE 6}\rho_3)$
(these inductions are explicitly described in section 1.2). We claim $\beta$ 
must be the 0-map. This can be seen by calculating $B_{\bbA 5}$ in a different
 order, as follows.
We can compute the $K$-homology of $cyl\cup\cO_{\bbA 5}$ using Mayer-Vietoris:
choosing our open cover $U,V$ appropriately, so both open sets are 
$G$-homeomorphic to $\bbR^2$ times $L:=(\bbR\times G/\bbA_5)\cup G/\bbT$,
we get
\begin{equation}\label{e6st2}
   \begin{matrix}{}^{\tau}K^{G}_0(cyl\cup\cO_{\bbA 5})&\longrightarrow&K^G_0
(L)\times 2&\longrightarrow&K^{G}_0(L)\times 2
\oplus R_{\bbA 5}\times 2  \cr\uparrow&&&&\downarrow\cr 
K^{G}_1(L)\times 2
   &{\longleftarrow}&K^{G}_1(L)\times 2&
   \longleftarrow&{}^{\tau}K^{G}_1(cyl\cup\cO_{\bbA 5})\end{matrix}
\end{equation}
The groups $K_*^G(L)$ are easily determined from the six-term sequence to
be $R_{\bbT}\oplus R_{\bbA 5},\,0$. Putting these into (\ref{e6st2}),
we obtain in the usual way
$^{\tau}K^{G}_*(cyl\cup\cO_{\bbA 5})=R_{\bbA 5},\,R_{\bbT}/(1-a^2)\oplus
R_{\bbA 5}^3$. We can cap this bounded cylinder by gluing in $\overline{ad}$
and $\overline{bc}$ as in Step 1, and we find that the $K$-homology of
$Tet\cup\cO_{\bbA 5}$ is $R_{\bbA 5},\,
(R_{\bbT}/(1-a^2)\oplus R_{\bbA 5}^3)\oplus \bbZ^21^1$ (the induction
(\ref{dindfinite}) from $R_{\bbA 5}$ to $^-R^1_{O2}$ will vanish). Finally,
gluing in $\cO_{\bbD 5}\cup\cO_{\bbE 6}$, we recover $^\tau K_*^G(B_{\bbA 5})$,
in particular obtaining that the map $^\tau K_0(B_{\bbA 5})\rightarrow
{}^\tau K_0(Tet\cup\cO_{\bbA 5})$ appearing in this final six-term is manifestly
surjective. Now, we can calculate $^\tau K_0^G(\cO_{\bbA 5})$ in exactly
a parallel way as $Tet\cup\cO_{\bbA 5}$, and we find that they are naturally
isomorphic. Therefore the map $\delta$ in (\ref{A5betagamma}) is likewise
surjective, and hence  $\beta$ must be identically 0.
(What we need in Step 4 below are $\beta,\gamma$ -- it is unnecessary to
determine the $K$-homology $^\tau K_*^G(B_{\bbA 5})$, though this is
now immediate.)

\bigskip
\centerline{{\bf Step 3:} {\it The  $\bbA_3$-ball}}\medskip

We want to repeat the above analysis, for the $\bbA_3$-ball.
The boundary of the genus-4 volume $B^2\setminus 4\bbR$ of (\ref{orba3})
is $Tet$ together with the four chords $\cO_{\bbD 5}\cup\cO_{\bbE 6}\cup\cO_{\bbD 4}$.
As in step 2, we know $K_*^G(\cO_{\bbA 3})=R_{\bbA 3},R_{\bbA 3}^4$.
First, we fill the $\cO_{\bbD 4}$ hole in $\cO_{\bbA 3}$, in order to obtain:
\begin{eqnarray}
^\tau K_0^G(\cO_{\bbA 3}\cup\cO_{\bbD 4})&=&\,\mathrm{coker\ (Ind}_{\bbA 3}^{\bbD 4})\oplus K_0^G(\cO_{\bbA 3})\nonumber\\&=&
(\mathrm{Span}\{[s_0-s_2],[s_1-s_3]\})\oplus R_{\bbA 3}\,,\\
^\tau K_1^G(\cO_{\bbA 3}\cup\cO_{\bbD 4})&=&\,\mathrm{ker\ (Ind}_{\bbA 3}^{\bbD 4})\oplus R_{\bbA 3}^3=\bbZ (r_\i'-r_{-\i}')\oplus R_{\bbA 3}^3
\simeq\bbZ^{13}\,.
\end{eqnarray}

Write $B_{\bbA 3}$ for the closed ball $Tet\cup \cO_{\bbD 5}\cup\cO_{\bbE 6}
\cup\cO_{\bbA 3}\cup\cO_{\bbD 4}$, we obtain the analogue of 
(\ref{A5betagamma}):
\begin{equation}\label{A3betagamma}
  \begin{matrix}\mathrm{Span}\{[s_0-s_2],[s_1-s_3]\}\oplus R_{\bbA 3} &\longleftarrow&{}^\tau K^{
{G}}_0(B_{\bbA 3})&\longleftarrow &R_{\bbD 5}^2\oplus R_{\bbE 6}
   \cr \beta'\downarrow&&&&\uparrow\,\gamma'\cr 
  \bbZ^2 1^-\oplus R_\bbT/(1-a^2) 
   &\longrightarrow&{}^\tau K^{{G}}_1(B_{\bbA 3})
   &\longrightarrow&  \bbZ(r_\i'-r_{-\i}')\oplus R_{\bbA 5}^3\end{matrix}
\end{equation}
The map $\gamma'$ will again be
the diagonal inductions, sending  $(x;\rho_1,\rho_2,\rho_3)$ to
$(\mathrm{Ind}_{\bbA 3}^{\bbD 5}
\rho_1, \mathrm{Ind}_{\bbA 3}^{\bbD 5}\rho_2, \mathrm{Ind}_{\bbA 3}^{\bbE 6}\rho_3)$.
Incidentally, this is why $\alpha$ in (\ref{e6alpha}) must vanish: the alternative
would be the (nontrivial) induction (\ref{dindfinite}) from $R_{\bbD 5}$ to
$^-R^1_{O2}$, which would  kill some of the $R^2_{\bbD 5}$ and $\gamma'$
here would no longer be defined everywhere in $R_{\bbA 5}^3$.
Again, we will have $\beta'\equiv 0$. To see this, recalculate the $K$-homology
of $B_{\bbA 3}$ by gluing the four chords $\cO_{\bbD 4}\cup\cO_{\bbD 5}\cup
\cO_{\bbE 6}$ simultaneously into $Tet\cup\cO_{\bbA 3}$: the analog of the
map $\gamma$ must as usual be diagonal inductions; the analog of the map
$\beta$ goes from $R_{\bbA 3}$ to $\bbZ^2 1^-\oplus R_{\bbT}/(1-a^2)$,
and must be 0 for easy reasons: (\ref{dindfinite}) is 0, and the $R_G$-module
property says $R_{\bbA 3}$ can't see $R_{\bbT}/(1-a^2)$. Thus $\beta'$ in
(\ref{A3betagamma}) is required to vanish identically.
(As in the previous step, all that we are after
here is to identify $\beta',\gamma'$ -- we don't need 
$^\tau K_*^G(B_{\bbA 3})$.)

\bigskip
\centerline{{\bf Step 4:} {\it Gluing together the  $\bbA_5$- and $\bbA_3$-balls}}\medskip

The final step is just the six-term, gluing $Tet\cup\cO_{\bbD 5}\cup\cO_{\bbE 6}$
to $\cO_{\bbA 5}\cup(\cO_{\bbA 3}\cup\cO_{\bbD 4})$:
\begin{equation}
  \begin{matrix}R_{\bbA 5}\oplus R_{\bbA 3}\oplus\bbZ
%\mathrm{Span}\{[s_0-s_2],[s_1-s_3]\} 
&\longleftarrow&{}^\tau K^{
{G}}_0(B_{\bbA 5}\cup\cO_{\bbA 3}\cup\cO_{\bbD 4})&\longleftarrow &
R_{\bbD 5}^2\oplus R_{\bbE 6}\cr
\beta+\beta'\downarrow&&&&\uparrow\,\gamma+\gamma'\cr 
  \bbZ^2 1^-\oplus  R_\bbT/(1-a^2) 
   &\longrightarrow&{}^\tau K^{{G}}_1(B_{\bbA 5}\cup\cO_{\bbA 3}\cup\cO_{\bbD 4})
   &\longrightarrow& R^3_{\bbA 5}\oplus  \bbZ %(r_\i'-r_{-\i}')
\oplus R_{\bbA 3}^3
\end{matrix}
\end{equation}
where the vertical maps (using obvious notation) are given explicitly in Steps 
2 and 3 above. We thus obtain the final answer, for the $K$-homology of the
complement of the generic orbits $\cO_{\bbA 1}$ in $Sp(4)$:
\begin{eqnarray}
^\tau K_0^G(Sp(4)\setminus gen)&=&\,
\mathrm{coker}\,(\gamma+\gamma')\oplus \mathrm{ker}\,(\beta+\beta')\nonumber\\&=&
R_{\bbA 5}\oplus R_{\bbA 3}\oplus\mathrm{Span}\{[s_0-s_2],
[s_1-s_3]\}\simeq\bbZ^{12}\,,\\ 
^\tau K_1^G(Sp(4)\setminus gen)&=&\,
\mathrm{coker}\,(\beta+\beta')\oplus \mathrm{ker}\,(\gamma+\gamma')\nonumber\\&=&\,
(\bbZ^2 1^-\oplus  \bbZ[a^{\pm 1}]/(1-a^2))\oplus\bbZ(r_\i'-r_{-\i}') 
\nonumber\\&&\oplus\,\mathrm{Span}^2\{r_\omega\!-\!r_{\omega^2},r_{-\omega}\!-\!
r_{-\omega^2},
r_1\!+\!2r_\omega\!-\!r'_1\!-\!r'_{-1},r_{-1}\!+\!2r_{-\omega}\!-\!r'_\i\!-\!
r'_{-i}\}\nonumber\\&&\oplus\,
\mathrm{Span}\{r'_\i\!-\!r'_{-\i},r_1\!+\!r_\omega\!+\!r_{\omega^2}\!-\!r_1'\!-\!r'_{-1},
r_{-1}\!+\!r_{-\omega}\!+\!r_{-\omega^2}\!-\!2r'_{\i}\}\nonumber\\
&&\simeq\bbZ^{16}\,.
\end{eqnarray}
Presumably, gluing in the generic orbits (i.e. those with stabiliser $\bbA_1$)
reduces this $\bbZ^{12},\bbZ^{16}$ to the $\bbZ^2$ we obtained in
section 6.1. We interpret this in the concluding section.

\section{Interpretations, Questions  and Speculations}

This paper is the first of a series devoted to deepening the connection between
twisted equivariant $K$-homology and conformal field theory. We constructed the relevant
bundles and provided several detailed calculations of $K$-homology which should
be relevant to the full systems, nimreps, branching coefficients,...
 of conformal embeddings and orbifolds.
This concluding section suggests some preliminary interpretations of these.

It should be remarked that although deep connections between $K$-theory and
conformal field theory/string theory have been known for some time, their precise
relationship is often very subtle. An old example of this is the D-brane 
charge group, which for space-time $X$ is identified with a twisted
nonequivariant $K$-theory $^\delta K^*(X)$, with twist $\delta\in H^3(X;\bbZ)$ 
given by the $H$-flux (determined locally by the $B$-field). This charge group 
can also be calculated independently from conformal field theory, at least
when $X$ is a Lie group $G$, and the answers agree
apart from a multiplicity $2^{{\rm rank}\,G}$ appearing in the $K$-homology. There 
still is no noncontroversial explanation
of this multiplicity (see \cite{Mon,GGR,Bra} and references 
therein for this story). A more recent example is ${}^\tau K_G^*(G)$ when $G$ 
is compact but nonsimply connected, as mentioned in section 1.4. 
The role of the $H^1$-twist in that $SO(3)$ example is to choose between
the identifications $\sigma_i\sim\sigma_{k+2-i}$ and
$\sigma_i\sim-\sigma_{k+2-i}$; the former holds for the standard Wess-Zumino-Witten
$SO(3)$ theory, and a possible physical realisation of the latter is proposed 
in \cite{F}. In this spirit, $H^1$-twists of $^\tau K^*_G(G)$ for finite $G$ is
often possible, although its conformal field theoretic meaning it seems
has never been explored. 

In sections 3.2 and 3.3 we computed the $K$-homology ${}^\tau K_G^*(H)$ of 
conformal embeddings of equal rank, and obtained a result of higher dimension 
than we would have naively expected, which would have been the dimension of the
corresponding full system. By comparison, the $K$-homology for the conformal embeddings
of sections 5 and 6 was smaller than the corresponding full systems.
What does all this mean? Consider for concreteness the $\bbT_2\rightarrow SU(2)_1$ 
example of section 3.2, with $k=1$. Perhaps we should
have stopped that calculation at the generic orbits $\bbR^2\times\bbT/C_2$, 
with $C_2$ stabiliser. Perhaps including the
other orbits (with stabiliser $\bbT$) incorrectly doubles the answer.
So more generally the full system of $H_k\rightarrow G_\ell$ should perhaps
be obtained as $K^H_0(S)$ for some $H$-invariant submanifold,
 and only rarely will $S=G$. We expect this to be the
correct explanation, as it seems to be in line with the discussion on
conjugacy classes given below. Another possible interpretation
of the section 3.2 calculation is that the $K$-homology calculation sees
two different versions of this conformal embedding, and adds both answers together. 
The two embeddings
of $\bbT$ in $SU(2)$ would be distinguished by their orientation; more generally,
the Weyl group would permute these conformal embeddings. This could tie in
with the appearance of the Weyl group in section 3.3, and with the aforementioned 
explanation of the D-brane charge group multiplicities.

The holomorphic orbifolds by finite groups (\cite{Ev} and section 1.4 above)
performs beautifully. Likewise, permutation orbifolds of quantum doubles of
finite abelian groups also works out (see the beginning of section 4). 
Another class of accessible and important orbifolds are the $\bbZ_2$-orbifolds
of lattice theories -- we will consider these in future work. By contrast,
the orbifold calculations in sections 4.2 and 4.3 also don't quite match what
happens with orbifolds in conformal field theory \cite{KS}: the primaries fixed
by the orbifold group $\bbZ_2$ would be doubled (`fixed point resolution') and 
the remaining primaries replaced by their $\bbZ_2$-orbits.
However in section 4.2, $K_1$ in equation (\ref{tpermorb}) is $k$-dimensional 
while $K_0$ is 
$k(k+1)$-dimensional, so we recover the doubling of the fixed points, but
not the folding of the remainder. In section 4.3 the reverse happens.
In particular, the primaries of the $S_2$-permutation orbifold of the $SU(2)$ level $k$ 
theory are parametrised by pairs $(i,j)$ where $0\le i<j\le k$ (these are the
$S_2$-orbits of nonfixed points), as well as double multiplicities of the 
fixed points $(i,i)$.
In section 4.3 we get the correct folding for the nonfixed
points, but not the correct doubling of fixed points. 
As mentioned at the beginning of section 4, it is clear that this $K$-homology
should only be an approximation;
it is tempting to guess that there is
a natural `symmetrising map' from the $K$-homology computed in those subsections
 to the groups (namely the centre of the crossed-product construction)
 exactly capturing this permutation orbifold; in one example
this map would be surjective and in the other it would be injective.

In sections 5 and 6 we study in detail the $SU(2)_4\rightarrow SU(3)_1$ and
$SU(2)_{10}\rightarrow Sp(4)_1$ conformal embeddings, which give rise to
the modular invariants called $D_4$ and $E_6$ respectively in the
$SU(2)_k$ list of Cappelli-Itzykson-Zuber \cite{CIZ}.
Perhaps the most interesting observation to come out of this analysis is
that the largest finite stabiliser in this action of $SU(2)$ on $SU(3)$
resp. $Sp(4)$, is called $\bbD_4$ resp. $\bbE_6$ on McKay's list \cite{McK}.
We expect this pattern to continue with the `$E_8$' conformal embedding
$SU(2)_{28}\rightarrow G_{2,1}$. We also expect the $\bbE_7$ group to arise in
this way, using the realisation of the `$E_7$' modular invariant by a $\bbZ_2$-orbifold
of the `$D_{10}$' modular invariant. To our knowledge the only other direct
relation between the A-D-E of Cappelli-Itzykson-Zuber and McKay's A-D-E of
finite subgroups of $SU(2)$ are some speculative remarks near the end of
\cite{HH} relating orbifolds of certain supersymmetric gauge theories with
$SU(n)_k$ modular invariants; see also \cite{HS}. (For a fairly direct construction of the
$SU(2)_k$ modular invariants from the {\it Lie} groups of A-D-E type,
see \cite{Na}.) Will the largest finite stabilisers in the conformal embeddings
for $SU(3)_k$, agree with the $SU(3)$-subgroups which \cite{HH} associate to
those modular invariants? These conformal embeddings, namely $SU(3)_3\rightarrow
SO(8)_1$, $SU(3)_5\rightarrow SU(8)_1$, $SU(3)_8\rightarrow E_{6,1}$, and 
$SU(3)_{21}\rightarrow E_{7,1}$, will be studied in the sequel to this paper.

The conformal embedding $SU(2)_4\rightarrow SU(3)_1$ considered in section 5
yields the $SU(2)_4$ modular invariant
$|\chi_0+\chi_4|^2+2|\chi_2|^2$, called $D_4$ in \cite{CIZ}. The full
system is 8-dimensional, consisting of two copies of the $D_4$
 diagram, as in Figure 2. To what extent can we see these $D_4$'s in the
 $K$-homological groups of section 5?
The distinguishing feature of the (unextended) $D_4$ diagram is the $S_3$
symmetry of the three endpoints, which fixes the central vertex. This
$S_3$ symmetry exists throughout section 5: on the full space $SU(3)$
it is generated by multiplying any orbit by a scalar matrix $\omega^iI$
(these form the centre of $SU(3)$),
and by complex conjugation -- all of these commute with the $SU(2)$ action.
This is realised by arbitrarily permuting the three connected components
of $\cO_{SU2}$ and $\cO_{O2}$, as well as $S_3=\mathrm{Aut}(\bbD_4/\langle-1
\rangle)$ (responsible
for permuting the three non-trivial one-dimensional representations $s_i$
of $\bbD_4$).

The final $K$-homological groups (\ref{d4answer}) of $SU(2)$ on $SU(3)$
are only one-dimensional. However we can see copies of the $D_4$ diagram
in $^\tau K_0^G(\cO_{SU2}\cup\cO_{O2}\cup\cO_{\bbT})$ and
$K_1^G(\cO_{\bbD 4})$. The map $\psi$ of (\ref{d4psi}) identifies these
two $D_4$'s. As with the previous examples, we would expect that there
is an $SU(2)$-invariant submanifold of $SU(3)$, whose $K$-homology
consists of two copies of $D_4$, one of which is in $K_0$ and defines the
nimrep; the three endpoints of that $D_4$ should be a copy of $Ver_1(SU(3))$.

Unfortunately, these considerations are made more difficult because the
$R_{SU2}$-module structure obtained in section 5.2 
from $K$-homology does not agree with that of the full system (which should  
come from $\alpha$-induction). Is there a  meaning in CFT of this other
$R_{SU2}$-module structure, manifest in the $K$-homology of section 5?
The presence of torsion in (\ref{d4answer}) suggests that we may need to
use $K$-theory over $\bbQ$, rather than $\bbZ$.

The conformal embedding $SU(2)_{10}\rightarrow Sp(4)_1$ considered in section 6
yields the $SU(2)_{10}$ exceptional modular invariant (\ref{Esix}),
called $E_6$ in \cite{CIZ}. The full
system is 12-dimensional, consisting of two copies of the $E_6$
 diagram, as in Figure 3. Again, the final $K$-homology groups (\ref{e6answer})
are likewise too small to contain the full system. However, the $\bbZ_2$
symmetry of each diagram presumably comes from the centre $\bbZ_2$ of
$Sp(4)$ -- this corresponds for instance to permuting the two connected
components in each of $\cO_{SU2},\cO_{O2},\cO_{\bbD 5}$. The (unextended)
$E_6$ diagram can be built from two $\bbA_5$ and one $\bbA_3$ groups,
glued together at the midpoint $-1$ (see (\ref{coxmos})).

Clearly, a key (but difficult) question is to see $\alpha$-induction directly in
these examples. This is responsible for the `correct' $R_G$-module structures
in the full system; they will differ from the `obvious' $R_G$-module
structure inherited directly from the $K$-homology, because there are two
$\alpha$'s (namely $\alpha^{\pm}$) 
but only one `obvious'. Incidentally, the obvious one is the
$R_G$-structure in both $^\tau K_G^*(G)$ \cite{FHT} and $^\tau K_{SU2}^*(SO3)$
\cite{BS}. For the finite groups the $\alpha^{\pm}$ are found in \cite{Ev} 
(see also section 1.4 of this paper). \cite{TX} finds
the natural ring structure on the $K$-groups of twisted equivariant $K$-theory
-- it is essentially the external Kasparov product in equivariant $KK$-theory.
In particular, for twisted equivariant $K$-homology (see Remark 4.30 of \cite{TX}),
there will be a graded product $^\tau K_i^G(X)\times {}^\tau K_j^G(X)
\rightarrow 
{}^\tau K_{i+j}^G(X)$ (at least when the twist is transgressed).
This should agree with the algebra of the full system, and from this we can
obtain the braiding etc.

Two distinct A-D-E graphs can be associated to a
given finite subgroup $\Gamma$ of $SU(2)$: in the {\it McKay} or {\it 
cohomological} picture \cite{McK,IN}, the vertices of the {\it extended}
Dynkin diagram are labelled with irreducible $\Gamma$-representations; and
in the {\it Du Val} or {\it homological} picture \cite{IR,Bry}, the vertices
of the corresponding {\it unextended} Dynkin diagram are labelled with
{\it nontrivial} conjugacy classes of $\Gamma$. The former is described in
section 1.2; the latter can be illustrated nicely in the $\bbE_8$ case. Indeed,
the binary icosahedral group $\bbE_8=\langle a,b,c\,|\,a^5=b^3=c^2=abc=-1\rangle$ 
has 8 conjugacy classes other than 1: we can take as representatives $a^i$ for
$i=1,2,3,4$, $b^i$ for $i=1,2$, and $c$, as well as $a^5=b^3=c^2=-1$. Then the
unextended $E_8$ Dynkin diagram is obtained by identifying an endpoint
(corresponding to $-1$) of a 5-chain (corresponding to the $a^i$), a 3-chain
(the $b^i$) and a 2-chain (the $c^i$). This really is the {\it unextended} 
diagram: 
the missing conjugacy class, namely 1, should be the affine vertex, but
its edges would be all wrong. (This homological picture is very reminiscent
of the generators of eg. the $E_8$ Lie group representation ring found in \cite{Ad}: denote
by $a,b,c$ the fundamental $E_8$-representations of dimension 248, 3875 and 
147250 respectively; then the exterior powers $\bigwedge^i a$ for $i\le 4$,
$\bigwedge^ib$ for $i\le 2$, and $c$, together with either $\bigwedge^5a$ or
$\bigwedge^3b$ or $\bigwedge^2c$, generates the polynomial ring $R_{E8}$.)

This cohomological picture is the one used exclusively in our calculations,
and indeed conformal field theory identifies the primaries (i.e. the basis  
of the Verlinde algebra $R_G/I_k$) directly with representations.
However, D-branes in Wess-Zumino-Witten models (the conformal field theories
of primary interest in this paper) have long been associated to conjugacy
classes in $G$ (see eg. \cite{FFFS,MMS}). More fundamentally, Verlinde's
formula (\ref{verl}) identifies all characters (one-dimensional representations)
$\lambda\mapsto S_{\lambda\kappa}/S_{0\kappa}$ of the Verlinde algebra,
labelling them with primaries $\kappa$. In the case of Wess-Zumino-Witten
models on a Lie group $G$, $S_{\lambda\kappa}/S_{0\kappa}$ equals the
character of the $\lambda$-representation of $G$, evaluated at the
(conjugacy class of the) element of finite order $\exp(2\pi\i\,(\kappa+\rho)/(k+h^\vee))$
in $G$. This is precisely the conjugacy class which eg. \cite{M} associates to
the primary $\kappa$.

It is interesting that \cite{M}, who like us works in the language of $K$-homology,
interprets primaries as special conjugacy classes. Indeed the representation
rings appearing in our $K$-homology groups arise through the (unnatural) application
of Poincar\'e duality. The Dynkin diagrams which both conformal field theory and
subfactors attach to the $SU(2)_k$ modular invariants and full systems are
{\it unextended}. This suggests that a more natural treatment of these $K$-homology
calculations could be to directly involve conjugacy classes rather than representations.
(See the end of section 3.1 for an explanation of how conjugacy classes arise 
in the $K$-theoretic treatment of the Verlinde algebra of a circle.)

String theory is well-defined on singular spaces, and many relations between
it and the McKay correspondence and resolution of quotient singularities,
have been explored. Let's briefly describe one which bears some formal
similarity to the considerations of section 2.2. In type IIB string theory,
the 10-dimensional background comes factored into $\bbR^{3,1}\times I^6$,
where the {\it transverse} or {\it internal} space $I^6$ is a Calabi-Yau
3-fold. For ease of calculation it is common to locally model $I^6$ with
a (non-compact) toric variety $Y$, eg. $\bbC^3/\Gamma$ where $\Gamma$ is
a finite abelian subgroup of $SU(3)$. These $Y$ have a singularity at the fixed point
$(0,0,0)$; probe that singularity with $N$ D3-branes. These branes will
fill space-time $\bbR^{3,1}$ but be localised to a single point (namely 0)
in the transverse space $Y$. We are supposed to live in this $\bbR^{3,1}$,
and the optimistic among us may hope that the low energy effective gauge
theory of this string theory (for the vacuum corresponding to $0\in Y$)
would be that of a supersymmetric extension of the Standard Model of
particle physics. The method of {\it brane tilings} (see eg. \cite{HV})
is a successful way of determining that effective theory. Periodic quivers
are polygonal tilings of a 2-torus, just as in Figure 5 (eg. the $Sp(4)$-Stiefel
diagram would correespond to the square tiling of the conifold); their dual 
(brane tilings) is what we would call the graph of cuts associated to our covers
of the bundles of section 2.2. The vertices of the periodic quiver label
irreducible $\Gamma$-representations. $\Gamma$ breaks the $U(N)$ gauge
symmetry of the $N$ coincident branes down to some $\prod_i SU(N_i)$
(one factor $SU(N_i)$ for each vertex) -- perhaps these gauge groups could
correspond to some of our stabilisers. To the directed edges of the quiver
graph are associated the matter fields of the theory (called `bi-fundamentals'
because they carry the $SU(N_i)\times SU(N_j)$ representation $\phi_i\otimes
\overline{\phi_j}$ of $N_i$-, $N_j$-dimensional fundamental representations).
These resemble the unitaries of section 2.2. The final ingredient needed
to identify the effective theory is the superpotential, a polynomial in
the bi-fundamentals, each term of which is determined from the faces of the
periodic quiver. For us these terms would be constants. These dimer models
are two-dimensional (so for us correspond to rank-2 Lie groups), but a three-dimensional
generalisation (corresponding to 3-dimensional superCFT, toric Calabi-Yau
4-fold singularities, and hence rank-3 groups) has been recently proposed
-- see eg. \cite{LLP}. We don't know if
there is anything deep underlying the formal similarities of these seemingly
independent pictures.

The $H^3$-untwisted `Verlinde algebra' for $SU(2)$ (so $k=\infty$)
is the representation ring $R_{SU2}$.
This can be realised through the $K$-homology of the fixed point algebra of the 
product action on the Pauli algebra of the infinite tensor product of $2 \times 2$ matrices
\cite{Wassthes}: $K_0(( \otimes _{\bbN} M_2)^{SU2}) \simeq R_{SU2}$. If we identify
$R_{SU2}$ with $\bbZ[t]$, the polynomials in an indeterminate $t$, then the 
non-zero elements of the positive cone
are identified with $\{P : P(t) > 0, t \in (0,1/4] \}$\,.
The fixed point algebra $(\otimes _{\bbN} M_2)^{SU2}$ is the generic Temperley-Lieb algebra.
Indeed if we deform this situation with a quantum group $SU_q(2)$ 
 to  $( \otimes _{\bbN} M_2)^{SU_q(2)}$ 
then we have the Temperley-Lieb algebra at say a root of unity $q = \exp(\i \pi (k +2))$.
The $K$-groups $K_0(( \otimes _{\bbN} M_2)^{SU_q(2)})$ of these algebras
can be identified with $\bbZ[t]/  \langle P_k \rangle$, the corresponding 
Verlinde algebra at level $k$ \cite{EvGould}. Here $P_i$ are the polynomials
defined by $P_i = P_{i-2} - tP_{i-1}, P_0 = 1, P_{-1} =0$, and the non-zero positive elements of
 the $K$-group are $\{ Q + \langle P_k\rangle \}$ where $Q((2 \cos(\pi /(k+2)))^{-2}) > 0\,.$
The ring structure on the Verlinde algebra is induced by the multiplication map
$( \otimes _{\bbN} M_2) \times ( \otimes _{\bbN} M_2)
\rightarrow \otimes _{\bbN} M_2 $ on the Pauli algebra. 
The $K$-group $K_0(( \otimes _{\bbN} M_2)^{SU_q(2)})$  should in turn be identified
with the equivariant $K$-group $K_0^{SU_q(2)}( \otimes _{\bbN} M_2)$ \cite{NeTu}.
This should generalise to the other groups. It would be interesting to 
pursue the considerations of this 
paper (conformal embeddings, orbifolds, etc) from this quantum group context.

 In these remaining lines,
we'll give  a small taste of the work in progress.
There is more to Freed-Hopkins-Teleman than writing the Verlinde algebra as a
$K$-group. They have bundles/equivalence$\, = \,$Verlinde algebra.
If we only look at equivalence classes, then we never see the braiding and
hence the associated representation of
the modular group $SL(2,\bbZ)$,
as the Verlinde algebra is commutative. But there is here a special choice of
isomorphism of bundle products
$V\times W\simeq W\times V$ which gives the braiding.
Similarly, it is useful to think in terms of concrete algebras  --
i.e. a graded equivariant bundle of compact operators over a space $X$,
with appropriate $H^3$ and $H^1$ invariants,
such that the $K$-theory of the $C^*$-algebra $A$ of its sections is really
$K(X)$. We should think in terms
of objects realising $K$-theory and $K$-homology, rather then just their
equivalence classes. The analysis of Freed-Hopkins-Teleman of the Verlinde algebra
already realises the primary fields as supersymmetric operators. Mickelsson 
wrote this out explicitly for the case of $\SUz$. If $\cA$ are the smooth 
$su(2)$-valued vector potentials 
on $\bbT$, write $Q_A = Q + \hat A ,\, A \in \cA$, where $Q$ is the free supercharge 
on $H_f \otimes H_b$ satisfying $Q^2 = L_0$, and $H_f$ are fermions with a level 2 represention
and $H_b$ are bosons with an irreducible level $k$ representation of $\LSUz $ and 
$\hat A $ is an interaction term. Then $Q_A$ is a family of self adjoint Fredholm operators,
equivariant with respect
to a central extension of the loop group $\LSUz $, and $\exp(-\i \pi\,
\mathrm{ sgn}(Q_A))$ basically defines 
an element in the $K$-group $K_{SU2}^1(SU(2))$. (The index computations of
Jaffe-Lesniewski-Weitsman \cite{Connes book} produce from supercharge operators $Q$ spectral
triples giving
elements in $K$-homology.) However it is not just the primary fields which need to be explicitly
realised in this way, but all the associated objects of a modular invariant, such as the boundary
$\NXM$
and the full system $\MXM$, including the Dirac-like canonical inclusion 
$\iota \in \NXN$, the canonical endomorphism $\theta = \bar \iota \iota$ in the
Verlinde algebra and the dual canonical
endomorphism $\gamma = \iota \bar \iota \in \MXM$ as spectral objects via Fredholm modules, Dirac
operators and
spectral triples. Indeed going beyond this, the maps between these $K$-groups, such as the modular
invariant itself, branching coefficients, sigma-restriction and alpha-induction should have
$KK$-theoretic interpretations.

In (\ref{extE}) we interpreted $K$-homology as classifying certain extensions.
More generally, the extensions
\begin{equation}\label{**}
0 \rightarrow K \otimes A \rightarrow E \rightarrow B \rightarrow 0
\end{equation}
together with suspensions, yield the Kasparov groups $KK_{*}(A,B)$ (page 118 of \cite{EKaw}). Now
by a Universal Coefficient Theorem 
there is an exact sequence
 $KK_{*}(A,B) \rightarrow {\rm Hom}(K_{*}(A),K_{*}(B)) \rightarrow 0$ 
as on page 120 of \cite{EKaw}.
In particular, taking $A = B$ to be the object giving the Verlinde algebra,
 a modular invariant is just an element of Hom$(K(A), K(A))$
and so gives rise to an element of $KK_1(A,A)$.
Hence a modular invariant will give rise to very special $KK$-elements, as would sigma restriction,
alpha induction,
 which should be analysed via spectral triples, Fredholm modules and Dirac 
operators.\footnote{Constantin Teleman recently informed us that he and Dan 
Freed have 
independently been trying to reconstruct the braiding from the Dirac families.}
This should be relatively straightforward for finite groups.

\newcommand\biba[7]   {\bibitem{#1} {#2:} {\sl #3.} {\rm #4} {\bf #5,}
                    {#6 } {#7}}
                    \newcommand\bibx[4]   {\bibitem{#1} {#2:} {\sl #3} {\rm #4}}

\def\ASENS            {Ann. Sci. \'Ec. Norm. Sup.}
\def\AM   {Acta Math.}
   \def\AnM              {Ann. Math.}
   \def\CMP              {Commun.\ Math.\ Phys.}
   \def\IJM              {Internat.\ J. Math.}
   \def\JAMS             {J. Amer. Math. Soc.}
\def\JFA              {J.\ Funct.\ Anal.}
\def\JMP              {J.\ Math.\ Phys.}
\def\JRA              {J. Reine Angew. Math.}
\def\JSP              {J.\ Stat.\ Physics}
\def\LMP              {Lett.\ Math.\ Phys.}
\def\RMP              {Rev.\ Math.\ Phys.}
\def\RNM              {Res.\ Notes\ Math.}
\def\RIMS             {Publ.\ RIMS.\ Kyoto Univ.}
\def\Inv              {Invent.\ Math.}
\def\npbp             {Nucl.\ Phys.\ {\bf B} (Proc.\ Suppl.)}
\def\nupb             {Nucl.\ Phys.\ {\bf B}}
\def\nup              {Nucl.\ Phys. }
\def\nupp             {Nucl.\ Phys.\ (Proc.\ Suppl.) }
\def\adma             {Adv.\ Math.}
\def\coma             {Con\-temp.\ Math.}
\def\PAMS             {Proc. Amer. Math. Soc.}
\def\PJM              {Pacific J. Math.}
\def\ijmp             {Int.\ J.\ Mod.\ Phys.\ {\bf A}}
\def\jpa              {J.\ Phys.\ {\bf A}}
\def\PLB              {Phys.\ Lett.\ {\bf B}}
\def\RIMS             {Publ.\ RIMS, Kyoto Univ.}
\def\Top               {Topology}
\def\TAMS             {Trans.\ Amer.\ Math.\ Soc.}

\def\Duke              {Duke Math.\ J.}
\def\K                 {K-theory}
\def\JOP               {J.\ Oper.\ Theory}

\vspace{0.2cm}\addtolength{\baselineskip}{-2pt}
\begin{footnotesize}
\noindent{\it Acknowledgement.}

The authors thank the University of Alberta, BIRS, Cardiff University,
Dublin IAS, ETH (Z\"urich), the Fields Institute, Universit\"at
Hamburg, the Newton Institute, and Swansea University for generous
hospitality while researching this
paper. Their research was supported in part by CTS (ETH), EPSRC
GR/S80592/01,
EU-NCG Research Training Network: MRTN-CT-2006 031962, Humboldt-Stiftung,
and NSERC. Help from John Jones, Ryszard Nest and especially Constantin
Teleman was
invaluable and much appreciated; we thank Yang-Hui He and especially
Matthias Gaberdiel for several useful comments on the manuscript.

\end{footnotesize}

\end{document}